\documentclass[dvips,11pt]{article}
\usepackage{amsmath,amssymb,amsfonts,amscd,amsthm,pstricks,pst-node,bm,xspace}
\def\simarrow{\mathrel{\raise -0.5mm\hbox{$\sim$}}\hspace{-1.8mm}{\rightarrow} } 

\def\bsimarrow{\leftarrow\hspace{-0.7mm}\mathrel{\raise -0.5mm\hbox{$\backsim$}} }

\def\bt{\begin{tabular}}
\def\te{\end{tabular}}

\def\lettrine#1#2#3{\noindent\hangindent#1\hangafter-#2
\hskip-#1\smash{\hbox to #1{#3\hfill}}\ignorespaces}

\newcommand{\To}[1]{\mathop{\to}\limits_{#1}}

\def\BM{\begin{pmatrix}}
\def\EM{\end{pmatrix}}

\def\pt{\mathrel{\scriptstyle \bullet}}

\def\ds{\displaystyle}

\def\d=f{\buildrel\hbox{\scriptsize d\'{e}f}\over \Longleftrightarrow}

\def\square{\hfill\hbox{\vrule height .9ex width .8ex depth -.1ex}}
\def\rit{\text{\it I\hskip -2pt  R}}

\def\nit{\text{\it I\hskip -2pt  N}}

\def\rl {\rit^{\hskip 1pt\ell}}

\def\Bd{{\text B}}

\def\Cs{{\cal C}}

\def\Ed{{\text E}}

\def\Fs{{\cal F}}

\def\Hs{{\cal H}}

\def\Os{{\cal O}}

\def\be{\begin{equation}}
\def\ee{\end{equation}}
\def\beqn{\begin{eqnarray}}
\def\eeqn{\end{eqnarray}}
\def\nobeqn{\begin{eqnarray*}}
\def\noeeqn{\end{eqnarray*}}
\def\ba{\left(\begin{array}}
\def\ea{\end{array} \right) }

\def\bpr{\paragraph{Proof.}}

\def\epr{\square\vskip 6pt}

\def\o{\overline}
\def\and{\; \mbox{and} \;}

\newcommand{\half}{\frac{1}{2}}

\def\hfl#1#2{\smash{\mathop{\hbox to 12mm{\rightarrowfill}}
\limits^{\scriptstyle #1}_{\scriptstyle #2}}}

\def\mod{\mathop{\rm mod}\nolimits}

\def\Be{\begin{enumerate}}
\def\Ee{\end{enumerate}}

\def\Bena{\begin{enumerate}
\def\labelenumi{\theenumi)}
\def\theenumi{\arabic{enumi}}
\def\labelenumii{\theenumii)}
\def\theenumii{\alph{enumii}}}

\def\Bean{\begin{enumerate}
\def\labelenumii{\theenumii)}
\def\theenumii{\arabic{enumii}}
\def\labelenumi{\theenumi)}
\def\theenumi{\alph{enumi}}}

\def\Bero{\begin{enumerate}
\def\labelenumii{\theenumii)}
\def\theenumii{\arabic{enumii}}
\def\labelenumi{(\theenumi)}
\def\theenumi{\roman{enumi}}}

\def\BeRo{\begin{enumerate}
\def\labelenumii{\theenumii)}
\def\theenumii{\arabic{enumii}}
\def\labelenumi{(\theenumi)}
\def\theenumi{\Roman{enumi}}}

\def\Bi{\vskip 11pt\begin{itemize}\itemsep=18pt}
\def\Ei{\end{itemize}\vskip 11pt}

\def\Bd{\begin{description}}
\def\Ed{\end{description}}

\def\R{\right}
\def\L{\left}
\def\F{\frac}

\usepackage[latin1]{inputenc}

\def\bigoplus{\mathop{\oplus}\limits}

\def\prod{\mathop{\Pi}\limits}
\def\sum{\mathop{\Sigma}\limits}

\def\bbf{\boldmath\bfseries}
\def\o{\overline}

\def\Bi{\begin{itemize}}
\def\Ei{\end{itemize}}

\newcommand{\ZZ}{\mathbb{Z}\,}
\newcommand{\NN}{\mathbb{N}\,}

\newcommand{\Aa}{\mathbb{A}\,}

\newcommand{\RR}{\mathbb{R}\,}
\renewcommand{\rit}{\RR}
\newcommand{\QQ}{\mathbb{Q}\,}
\newcommand{\FF}{\mathbb{F}\,}

\renewcommand{\SS}{\mathbb{S}\,}

\newcommand{\DD}{\mathfrak{D}}

\newcommand{\BB}{\mathsf{B}}
\renewcommand{\BB}{\mathfrak{B}}

\def\tr{\operatorname{trace}}

\def\Gal{\operatorname{Gal}}
\def\Frob{\operatorname{Frob}}
\def\End{\operatorname{End}}
\def\mod{\operatorname{mod}}

\def\Im{\operatorname{im}}

\def\cusp{\operatorname{cusp}}
\def\Irr{\operatorname{Irr}}

\def\Rep{\operatorname{Rep}}

\def\Repsp{\operatorname{Repsp}}
\def\FRepsp{\operatorname{FRepsp}}
\def\FREPsp{\operatorname{FREPsp}}
\def\Nu{\operatorname{Nu}}
\def\Em{\operatorname{em}}
\def\INF{\operatorname{INF}}

\def\ELLIP{\operatorname{ELLIP}}
\def\Ellip{\operatorname{ELLIP}}
\def\ellip{\operatorname{ellip}}

\def\GL{\operatorname{GL}}
\def\gl{\operatorname{gl}}

\def\bM{\begin{matrix}}
\def\eM{\end{matrix}}

\def\pt{\bullet}
\def\lr{left (resp. right) }
\def\rl{right (resp. left) }

\def\bp{\beta ^+_{p_r}}
\def\bmp{\beta ^-_{p_r}}

\def\resp#1{(resp. #1)}
\def\rresp#1{\qquad \mbox{(resp.} \quad #1\ )}

\def\bpr{\noindent{\bf{Proof\/}}:\;\;}
\def\To{\begin{CD} @>>>\end{CD}}

\def\RL{_{R\times L}}
\def\Rl{_{R,L}}

\def\wt{\widetilde}

\begin{document}

\setcounter{page}{0}
{\pagestyle{empty}
\null\vfill
\begin{center}
{\LARGE The Goldbach conjecture resulting from global-local cuspidal representations and deformations of Galois representations}
\vfill
{Christian {\sc Pierre\/}}
\vskip 11pt


\vfill

\begin{abstract}
In the basic general frame of the Langlands global program, a local $p$-adic elliptic semimodule corresponding to a local (left) cuspidal form is constructed from it global equivalent covered by $p^\ell$ roots.

In the same context, global and local bilinear deformations of Galois representations inducing the invariance of their respective residue fields are introduced as well as global and local bilinear quantum deformations leaving invariant the orders of the inertia subgroups.  More particularly, the inverse quantum deformation of  a closed curve responsible for its splitting directly leads to the Goldbach conjecture.

\end{abstract}
\end{center}
\vfill

Mathematics subject classification (2000): 11F03, 11P32, 11R37, 14B12.
\eject

\tableofcontents
\vfill\eject}
\setcounter{page}{1}
\def\thepage{\arabic{page}}
\parindent=0pt 
\section{Introduction}

In spite of its apparent simplicity, the Goldbach conjecture asserting that ``every integer superior or equal to 4 is the sum of two prime numbers'' has resisted since 1742 \cite{Gol} to a convincing proof or disproof.

As for the Wiles demonstration of Fermat's last theorem \cite{Wil}, it it perhaps a sign that {\bbf these striking problems of number theory\/} cannot be directly solved by means of this unique field but {\bbf have to be tackled in a more general context\/} for example from the program of Langlands.

The aim of this paper thus consists in approaching the conjecture of Goldbach on the basis of new breakthroughs in the (bilinear) global-local cuspidal representations and in the deformations of Galois representations.
\vskip 11pt

{\bbf The main results of this work deal with:}
\Bean
\item {\bbf the generation of a local $p$-adic cuspidal form \cite{Lan2} from its global equivalent covered by $p^\ell$ roots: this leads directly to the Serre-conjecture on $p$-adic Galois representations associated with modular forms and to the Shimura-Taniyama-Weil conjecture\/}.

\item {\bbf global and local bilinear deformations of Galois representations\/} inducing the invariance of their respective residue fields {\bbf and the introduction of global and local bilinear quantum deformations\/} leaving invariant the orders of the inertia subgroups.
\Ee
\vskip 11pt

{\bbf The general basic mathematical frame\/} used in this work {\bbf is that of the Langlands global program\/} \cite{Pie1} dealing with bijections between the equivalence classes of the $n^2$-dimensional representation of the bilinear global Weil group given by the bilinear algebraic semigroup $\GL_n(L_{\o v}\times L_v)$ and the corresponding conjugacy classes of the cuspidal representation of $\GL_n(L_{\o v}\times L_v)$ where $L_v$ \resp{$L_{\o v}$} denotes the complete set of \lr archimedean pseudo-ramified completions.
\vskip 11pt

{\bbf Non abelian global class field concepts\/} are reviewed in chapter~2.  They are based on a set of increasing finite symmetric splitting semifields characterized by Galois extension degrees which are integers modulo $N$.  The corresponding completions, resulting from isomorphisms of compactification, are infinite archimedean pseudo-ramified completions, defining archimedean pseudo-ramified real places.

{\bbf Bilinear algebraic semigroups\/} $\GL_n(L_{\o v}\times L_v)=T^t_n(L_{\o v})\times T_n(L_v)$, over the product of the set $L_v$ of left archimedean pseudo-ramified completions by the symmetric set $L_{\o v}$ of right completions, {\bbf are expressed according to the Gauss bilinear decomposition\/}, i.e. by means of the product of the group $T^t_n(L_{\o v})$ of lower triangular matrices with entries in $L_{\o v}$ and referring to the lower half space by the group $T_n(L_v)$ of upper triangular matrices with entries in $L_v$ and referring to the upper half space.

The algebraic representation space $\Repsp(\GL_n(L_{\o v}\times L_v))$ of the bilinear algebraic semigroup $\GL_n(L_{\o v}\times L_v)$ decomposes according to its conjugacy class representatives
$g^{(n)}\RL[i,m_i]$, $1\le i\le t\le \infty $, where $m_i$ refers to their multiplicities.

On the toroidal compactification 
\[G^{(n)}(L^T_{\o v}\times L^T_v)\equiv \Repsp(\GL_n(L^T_{\o v}\times L^T_v))\]
 of
 $ \Repsp(\GL_n(L_{\o v}\times L_v))$, the bisemisheaf $\Phi (G^{(n)}
(L_{\o v}\times L_v))$ of differentiable smooth bifunctions
\[ \phi \RL(g^{(n)}_{T\RL}[i,m_i])=
 \phi(g^{(n)}_{T_R}[i,m_i])\times
\phi(g^{(n)}_{T_L}[i,m_i])\]
is defined where
$g^{(n)}_{T_R}[i,m_i]$ \resp{$g^{(n)}_{T_L}[i,m_i]$} is a toroidal compactified \rl conjugacy class representative, i.e. a $n$-dimensional semitorus restricted to the lower \resp{upper} half space.

On the set of bisections of $\Phi (G^{(n)}
(L^T_{\o v}\times L^T_v))$, {\bbf a global elliptic
$\Gamma (\Phi (G^{(n)}(L^T_{\o v}\times L^T_v)))$-bisemimodule
$\ELLIP_R(n,i,m_i)\otimes_D \ELLIP_L(n,i,m_i)$ is constructed in such a way that:\/}
\Bean
\item {\bbf $\ELLIP_R(n,i,m_i) \otimes_D
\ELLIP_L(n,i,m_i)$ covers the corresponding cuspidal form
$f_R(z)\otimes_D f_L(z)$\/} as introduced in \cite{Pie1} and in \cite{Pie2}.

\item $
\begin{aligned}[t]
\ELLIP_L(n,i,m_i)&= \bigoplus^t_{i=1}\bigoplus_{m_i}
\lambda ^{\half}(n,i,m_i)\ e^{2\pi i(i)x}\; , \quad x\in\rit^n\\
\rresp{\ELLIP_R(n,i,m_i)&=  \bigoplus^t_{i=1}\bigoplus_{m_i}
\lambda ^{\half}(n,i,m_i)\ e^{-2\pi i(i)x}}, \end{aligned}$\\
is the sum of smooth differentiable functions (resp. cofunctions) on the conjugacy class representatives
$g^{(n)}_{T_L}[i,m_i]$ \resp{$g^{(n)}_{T_R}[i,m_i]$}, where $\lambda ^\half(n,i,m_i)$ is the square root of the considered Hecke character.
\Ee

{\bbf $\ELLIP_R(n,i,m_i) \otimes_D
\ELLIP_L(n,i,m_i)$ then constitutes a cuspidal representation space of the algebraic bilinear semigroup $\GL_n(L_{\o v}\times L_v)$\/} as required by the global program of Langlands and recalled in section 3.1.
\vskip 11pt

Let now $L_{[v_p]}$
\resp{$L_{[\o v _p]}$} denote the truncated set of \lr archimedean pseudo-ramified completions superior and equal to the $p$-th infinite place, where $p$ is a prime integer.

{\bbf A $n$-dimensional global elliptic
$\Gamma (\Phi (G^{(n)}(L^T_{[\o v_p]}\times L^T_{[\o v_p]})))$-bisemimodule
$\ELLIP_R(n,\linebreak i\ge p,m_i)\otimes
\ELLIP_L(n,i\ge p,m_i)$ can then be envisaged as well as its covering\/} global elliptic bisemimodule
$\ELLIP_R(n,i\ge p,m_i,\o p^{(\ell)})\otimes
\ELLIP_L(n,i\ge p,m_i,\o p^{(\ell)})$ {\bbf by $p^\ell$ roots\/} in such a way that each term
$\ellip_L(n,p+h,m_{p+h},\o p^\ell)\simeq 
\ds\prod^n_{c=1}r_c(p+h,m_{p+h},\o p^\ell)^{p^\ell}\ e^{2\pi i(p^\ell)x'_c}$, $x'_c\in\rit$, of 
$\ELLIP_L(n,i\ge p,m_i,\o p^{(\ell)})$ is the covering by $p^\ell$ roots (``$\ell$'' varying from one term to another) of the corresponding term
$\ellip_L(n,[p+h],m_{p+h})$ of $\ELLIP_L(n,i\ge p,m_i)$.

It is then proved that every $n$-dimensional semitorus $T^n_L[p+h,m_{p+h},\o p^\ell]\simeq \ellip_L(n,[p+h],m_{p+h})$ is a discrete valuation (semi)ring of which uniformizing element is
$r(p+h,m_{p+h},\o p^\ell)^p
=\ds\prod^n_{c=1}r_c(p+h,m_{p+h},\o p^\ell)^p$
and units are the invertible elements 
$e^{2\pi  i(p^\ell)x'}=\ds\prod^n_{c=1}\ e^{2\pi  i(p^\ell)x'_c}$, $x'\in\rit^n$.

The Kronecker-Weber theorem, expressing that every finite abelian extension of $\QQ$ is contained in a cyclotomic extension of $\QQ$, directly follows from the precedent considerations.
\vskip 11pt

{\bbf A local $p$-adic elliptic semimodule corresponding to a local left cuspidal form is then constructed\/} in section 3.3 {\bbf from its global equivalent covered by $p^\ell$ roots\/}.

To this end, it is shown that a set $\{\wt L_{v_{p+h}}\}_{h}$ of finite increasing global subsemifields ``above $p$'' can be covered in a etale way by (a) $p$-adic finite extension field(s) leading to a global($\leftrightarrow $)local one-to-one correspondence if the number of global and local elements correspond, i.e. if the number of ``global'' algebraic points is a power of $p$.

Starting with the two-dimensional global left elliptic
$\Gamma (\Phi (G^{(2)}(L^T_{[v_p]})))$-semimodule
$\ELLIP_L(2,\linebreak i\ge p,m_i,\o p^{(\ell)})$, covered by $p^\ell$ roots, the local elliptic left
$\End(G^{(2)}(K^+_p))$-semimodule\linebreak
$\ELLIP(2,x,K^+_p)=\ds\bigoplus_r\lambda ^\half_p(2,r,m_r)(x)\ f(\mu ^{q_r\cdot r})$ is constructed in one-to-one correspondence where $x$ is a closed point of the finite Galois extension of the non archimedean $p$-adic left semifield $L^+_p$ introduced in section 3.3.1\\
in such a way that:
\Bean
\item the Frobenius substitution $\mu \to\mu ^{q_r}$ on every local Frobenius endomorphism $\mu :x\to x^p$ be considered where $q_r=\sum_r f_r\cdot e_r$ with $f_r$ the ``local'' residue degree of the $r$-th prime ideal of the considered Galois extension and $e_r$ the corresponding ramification index.

\item $e_{L_{[v_p]}\to k^+_p}:\lambda (2,p+h,m_{p+h})\to \lambda _p(2,r,m_r)$ be the embedding of the product $\lambda (2,p+h,m_{p+h})$, right by left, of Hecke characters over $L_{[v_p]}$ into its equivalent $\lambda _p(2,r,m_r)$ over $K^+_p$.
\Ee

This condition corresponds to the embedding $i(a_\ell)=\tr(\Frob_\ell)$ into $\QQ_p$ of the ring of integers of a finite extension $E_f$ (i.e. the ring of the coefficients of the cuspidal form $f$) of $\QQ$, $a_\ell$ being the coefficient of the cuspidal form.
\vskip 11pt

{\bbf The Serre\/} (Eichler, Deligne, Shimura) {\bbf conjecture\/} \cite{C-F-T}, \cite{D-S}, \cite{Swi}, {\bbf asserting that Galois representations $\rho =G_\QQ\to\GL_2(\QQ_p)$ can be associated to modular forms directly results from this construction of a local $p$-adic cuspidal form\/}.

On this basis, {\bbf two kinds of explicit deformations of $n$-dimensional representations of Galois or Weil groups given by bilinear algebraic semigroups over complete global and local Noetherian bisemirings are considered in chapter 4\/}.
\vskip 11pt

First, {\bbf local $p$-adic coefficient semiring homomorphisms are envisaged in such a way that they induce an isomorphism on their residue semifields\/} leading to a base change in the considered finite Galois extensions.  Similarly, {\bbf global coefficient semiring homomorphisms are defined in such a way that they induce an isomorphism on their global residue semifields\/}.

It is then proved that the inverse image of the homomorphism 
$h_{L'_{L_p}\to L_{L_p}}$ between global coefficient semirings
$L'_{L_p}\equiv L'_{[v_p]}$ and
$L_{L_p}\equiv L_{[v_p]}$
 is isomorphic to the inverse image of the homomorphism $h_{B^{'+}_p\to B^+_p}$ between local coefficient semirings $B^{'+}_p$ and $B^+_p$ if the number of elements of the global kernel $K(h_{L'_{L_p}\to L_{L_p}})$ is equal to the number of elements of the local kernel $K(h_{B^{'+}_p\to B^+_p})$, i.e. is a power of $p$.
\vskip 11pt

{\bbf A $n$-dimensional global bilinear deformation of\/}
\[
\rho _L:\qquad \Gal({\dot{\wt L}}_{R_p}/k)\times
\Gal({\dot{\wt L}}_{L_p}/k)\To\GL_n(L_{R_p}\times L_{L_P})\]
{\bbf is then an equivalence class of liftings\/}
\[\rho _{L'_c}=\rho _L+\delta \rho _{L'_c}\;, \qquad 1\le c\le\infty \;,\]
 where $\delta \rho _{L'_c}$ refers to the kernel of $\rho _L$, in such a way that the kernels of two deformed algebraic bilinear semigroups 
$\GL_n(L'_{R_{p_{c_1}}} \times L'_{L_{p_{c_1}}})$ and
$\GL_n(L'_{R_{p_{c_2}}} \times L'_{L_{p_{c_2}}})$
differ by powers of orders of their inertia bilinear subgroups.
\vskip 11pt

Similarly, {\bbf a $n$-dimensional local $p$-adic bilinear deformation of\/}
\[\rho _K:\Gal(K^-_p/L^-_p)\times \Gal(K^+_p/L^+_p))\to
\GL_n(K^-_p\times K^+_p)\;,\]
 in the sense of Mazur \cite{Maz2}, {\bbf is an equivalence class of liftings\/}
 \[ \rho _{K'_d}=\rho _K+\delta \rho _{K'_d}\;, \qquad 1\le d\le\infty \;,\]
 where $\delta \rho _{K'_d}$ refers to the kernel of $\rho _K$, in such a way that the kernels of the two deformed algebraic bilinear semigroups 
$\GL_n(K^{'-}_{p_{d_1}}\times K^{'+}_{p_{d_1}})$ and
$\GL_n(K^{'-}_{p_{d_2}}\times K^{'+}_{p_{d_2}})$ differ by powers of their ramification indices.

A second type of deformations of Galois representations, called {\bbf quantum deformations\/}, is envisaged on the basis of global and local coefficient semiring quantum homomorphisms.

{\bbf A uniform quantum homomorphism\/}
\[Qh_{L_{L_{p+j}}\to L_{L_p}}: \qquad L_{L_{p+j}}\To L_{L_p}\]
{\bbf between two global compactified coefficient semirings\/}
is such that:
\Bi
\item it induces an isomorphism on their global inertia subgroups;
\item it increases the global residue semifield $L_{L_p}$ by an increment of $j$ quanta, i.e. $j$ irreducible closed algebraic subsets of degree $N$, on every completion.
\Ei

Similarly, {\bbf a quantum homomorphism\/}
\[ Qh_{B^+_{p_t}\to B^+_{p_r}}: \qquad B^+_{p_t}\To B^+_{p_r}\;, \quad t=r+s\;, \]
{\bbf between two local coefficient semirings\/} is such that:
\Bi
\item it induces an isomorphism on their ``local'' inertia subgroups (having thus the same ramification index);
\item it increases the residue degrees of the $r$ residue subsemifields of $B^+_{p_r}$ by a same integer increment.
\Ei
\vskip 11pt

{\bbf A $n$-dimensional global bilinear quantum deformation of $\rho _L$\/}, defined above, is then {\bbf an equivalence class of liftings\/}
\[ \rho _{L_j}=\rho _L+\delta \rho _{L_j}\;, \qquad 1\le j\le \infty \;, \]
where $\delta \rho _{L_j}$ refers to the kernel of $\rho _L$, in such a way that the kernels of the two deformed algebraic bilinear semigroups 
$\GL_n(L_{R_{p+j_1}} \times L_{L_{p+j_1}})$ and
$\GL_n(L_{R_{p+j_2}} \times L_{L_{p+j_2}})$ differ by powers of their global residue degrees.

Similarly, {\bbf a $n$-dimensional local bilinear ``quantum'' deformation of\/}
\[
\rho _{K_{p_r}} : \quad \Gal(K^-_{p_r}/L^-_p) \times
\Gal(K^+_{p_r}/L^+_p)  \To \GL_n( K^-_{p_r}\times K^+_{p_r} )\]
in an equivalence class of liftings:
\[ \rho _{K_{p_t}}=\rho _{K_{p_r}}+\delta \rho _{K_{p_t}}\;, \]
where $\delta \rho _{K_{p_t}}$ refers to the kernel of $\rho _{K_{p_r}}$, in such a way that the kernels of the two deformed algebraic bilinear semigroups
$\GL_n(K^-_{p_{t_1}}\times K^+_{p_{t_1}})$ and
$\GL_n(K^-_{p_{t_2}}\times K^+_{p_{t_2}})$ differ by powers of their local residue degrees.

Taking into account the Langlands global correspondences
\[
\ELLIP\FREPsp(\GL_n(L{\RL}_p)): \qquad
\GL_n(L_{R_p}\times L_{L_p}) \To
\ELLIP\RL(n,i\ge p,m_i)\;, \]
between the bilinear algebraic semigroup
$\GL_n(L_{R_p}\times L_{L_p})$ and the $n$-dimensional global elliptic bisemimodule
$\ELLIP\RL(n,i\ge p,m_i)$ introduced above,
{\bbf a $n$-dimensional global elliptic bilinear quantum deformation of
$\rho ^{\ELLIP}_L=\ELLIP\FREPsp(\GL_n(L{\RL}_p)\circ \rho _L$ is an equivalence class of liftings\/}
\[ \rho ^{\ELLIP}_{L_j}=
\rho ^{\ELLIP}_L+\delta \rho ^{\ELLIP}_{L_j}\]
inducing the injective morphism
\[\DD^{\{p\}\to\{p+j\}}\RL (n): \qquad
\ELLIP\RL(n,i\ge p,m_i)\To
\ELLIP\RL(n,i\ge p+j,m_i)\]
which is quantum deformation of $\ELLIP\RL(n,i\ge p,m_i)$
increasing the global residue degree ``$i$'' of each left and right term $i\ge p$ of
$\ELLIP\RL(n,i\ge p,m_i)$ by an amount of an integer ``$j$''.

{\bbf The injective morphism\/}
\[ \DD\RL^{[p]\to[p+j]}(n) : \qquad
\ellip\RL(n,[p],m_p) \To \ellip\RL(n,[p+j],m_{p+j})\]
{\bbf restricted to the $(p,m_p)$-th conjugacy class representative
$\ellip\RL(n,[p],m_p)$ of $\ELLIP\RL(n,i\ge p,m_i)$, is a quantum equivalence class representative of liftings or an elliptic quantum deformation\/} associated with the exact sequence:
\[ 1\To \ellip\RL(n,[j])\To \ellip\RL(n,[p+j],m_{p+j}) \To \ellip\RL(n,[p],m_p] \To 1\;.\]
Let then
\[ \DD^{[p+j+k]\to[p+j]}\Rl(1): \quad
\ellip\RL(1,[p+j+k],m_{p+j+k})
\To\ellip\RL(1,[p+j],m_{p+j})\]
denote the inverse elliptic quantum deformation of a one-dimensional global elliptic subbisemimodule of class $[p+j+k]$ towards a class $[p+j]$.

{\bbf This inverse elliptic quantum deformation corresponds to the endomorphism\/}:
\begin{multline*} \End^{[p+j+k]\to[p+j]}\Rl(1): \qquad
\ellip\RL(1,[p+j+k],m_{p+j+k})\\
\To\ellip\RL(1,[p+j],m_{p+j})+
\ellip\RL(1,[k],m_{k})\end{multline*}
where $\ellip\RL(1,[k],m_k]$ denotes the product, right by left, of a right semicircle at $k$ quanta, i.e. characterized by a global residue degree $f_{v_k}=k$ according to section 2.1, and localized in the lower half space by its left equivalent localized in the upper half space.

If we consider now the inverse quantum deformation of a closed curve  isomorphic to the left (or right) undoubled semicircle of class 
$[p'+j'+k']$, we get the following relation, associated with its endomorphism:
\[ f_{v_{p'+j'+k'}}=f_{v_{p'+j'}}+f_{v_{k'}}\qquad \text{or} \qquad
p'+j'+k' = (p'+j')+k'\]
for the resulting global residue degrees and corresponding to a splitting of the closed curve
$c^1_{2_L}[p'+j'+k']$ of class $[p'+j'+k']$ into two complementary curves of classes $[p'+j']$ and $[k']$.

{\bbf Taking into account that the global residue degree $f_{v_{p'+j'+k'}}=p'+j'+k'$ of a closed curve must be an even integer $G_{\rm even}$, {\bbf the relation $G_{\rm even}=f_{v_{p'+j'}}+f_{v_{k'}}$ directly leads to the Goldbach conjecture\/} on the basis of the developments of chapter 5 and of \cite{Pie2} dealing with the results of the author on the Riemann hypothesis\/}.

\section{New concepts of non abelian global class field theory}

\subsection*{2.1 Global class field concepts}

Let $k$ be a global number field of characteristic $0$ and let $\wt L$ denote a finite extension of $k$.

$\wt L=\wt L_R\cup \wt L_L$ is assumed to be a {\bbf symmetric splitting field\/} composed of a right and a left algebraic extension semifields $\wt L_R$ and $\wt L_L$ in one-to-one correspondence.  $\wt L_L$ and $\wt L_R$ are respectively the sets of positive and symmetric negative simple roots of a polynomial ring over $k$.

$\wt L_L$ and $\wt L_R$ are commutative division semirings, i.e. semifields, because they lack for opposite elements with respect to the addition.

Let
$\wt L_{v_1}\subset \dots \subset\wt L_{v_i}\subset \dots \subset\wt L_{v_t}$
\resp{$\wt L_{\o v_1}\subset \dots \subset\wt L_{\o v_i}\subset \dots \subset\wt L_{\o v_t}$}
denote the set of increasing subsemifields of 
$\wt L_L$ \resp{$\wt L_R$}.

{\bbf The completion $ L_{v_i}$ \resp{$ L_{\o v_i}$} associated with 
 $\wt L_{v_i}$ \resp{$\wt L_{\o v_i}$}\/}
 is an isomorphism of compactification
 $c_{v_i}:\wt L_{v_i} \to L_{v_i}$
 \resp{$c_{v_i}:\wt L_{\o v_i} \to L_{\o v_i}$}
of $\wt L_{v_i}$ 
\resp{$\wt L_{\o v_i}$} onto the subsemifield 
$L_{v_i}$
\resp{$L_{\o v_i}$} which is a closed compact subset of $\rit_+$
\resp{$\rit_-$} \cite{Kna}, \cite{Ser3}.

The equivalence classes of completions of 
$\wt L_L$
\resp{$\wt L_R$}, characterized by they number of elements, are the \lr {\bbf infinite real places\/} of $\wt L_L$
\resp{$\wt L_R$}.

They are noted
$v=\{v_1,\dots,v_i,\dots,v_t\}$
\resp{$\o v=\{\o v_1,\dots,\o v_i,\dots,\o v_t\}$}.

Let $ L_{v_i}$ \resp{$ L_{\o v_i}$} denote the {\bbf infinite pseudo-ramified completion\/} proceeding from the subsemifield
$\wt L_{v_i}$ \resp{$\wt L_{\o v_i}$}:
\Bi
\item it is characterized by an achimedean absolute value in its topology;
\item it is generated from an irreducible central $k$-semimodule
$L_{v^1_i}$
\resp{$L_{\o v^1_i}$} of rank\linebreak 
$[L_{v^1_i}:k]=N$
\resp{$[L_{\o v^1_i}:k]=N$};

\item it is defined by its rank, i.e. its Galois extension degree
$[L_{v_i}:k]=[L_{\o v_i}:k]=\star+i\dot N$ which is an integer modulo $N$, $N\in\nit$, where $\star$ denotes an integer inferior to $N$.
\Ei

The corresponding pseudo-unramified completion
$L^{nr}_{v_i}$
\resp{$L^{nr}_{\o v_i}$}
is defined by:
\[ L^{nr}_{v_i}=L_{v_i}\big/ L_{v^1_i}\qquad
\rresp{L^{nr}_{\o v_i}=L_{\o v_i}\big/ L_{\o v^1_i}}\]
and is characterized by its global residue degree (in analogy with the local $p$-adic treatment)
$f_{v_i}=[L^{nr}_{v_i}:k]=i$
\resp{$f_{\o v_i}=[L^{nr}_{\o v_i}:k]=i$}.

The {\bbf infinite pseudo-ramified real place $v_i$
\resp{$\o v_i$}} is composed of the basic completion
$L_{v_i}$
\resp{$L_{\o v_i}$} and of the set
$\{L_{v_{i,m_i}}\}^{\sup(m_i)}_{m_i=1}$
\resp{$\{L_{\o v_{i,m_i}}\}^{\sup(m_i)}_{m_i=1}$}
of equivalent completions characterized by the same rank, where $\sup (m_i)$ denotes the multiplicity of $v_i$ \resp{$\o v_i$}.

The set \resp{the direct sum} of the real infinite pseudo-ramified completions is given by\linebreak 
$L_v=\{L_{v_{i,m_i}}\}_{i,m_i}$ or $L_{\o v}=\{L_{\o v_{i,m_i}}\}_{i,m_i}$
\resp{$L_{v_\oplus}=\bigoplus^t_{i=1} \bigoplus_{m_i} L_{v_{i,m_i}}$
 or $L_{\o v_\oplus}=\bigoplus^t_{i=1} \bigoplus_{m_i} L_{\o v_{i,m_i}}$}
 and the product of the primary real infinite pseudo-ramified completions gives rise to the adele semiring:
 \[ \Aa^\infty_{L_v}=\prod_{j_p,m_{j_p}} L_{v_{j_p,m_{j_p}}} \qquad
 \rresp{\Aa^\infty_{L_{\o v}}=\prod_{j_p,m_{j_p}} L_{\o v_{j_p,m_{j_p}}}}\]
 where $j_p$ denotes the $j$-th primary completion.
 \vskip 11pt
 
 
 \subsection*{2.2 Galois, inertia and Weil groups}

Let $\widetilde L_{v_i}$ (resp. $\widetilde L_{\o v_i}$) and $\widetilde L_{v_i,m_i}$ (resp. $\widetilde L_{\o v_i,m_i}$) 
denote respectively the the basic and the equivalent Galois extensions corresponding to the basic completion $L_{v_i}$ 
(resp. $L_{\o v_i}$) and to the equivalent completion $L_{v_i,m_i}$ (resp. $L_{\o v_i,m_i}$) at the $v_i$-th archimedean place.

Let $\Gal^D(\widetilde L_{v_i}\big/ k)$ \resp{$\Gal^D(\widetilde L_{\o v_i}\big/ k)$} denote the Galois subgroup of 
$\widetilde L_{v_i}$ \resp{$\widetilde L_{\o v_i}$}  and let
$\Gal(\widetilde L_{v_i,m_i}\big/ k)$ \resp{$\Gal(\widetilde L_{\o v_i,m_i}\big/ k)$} denote the Galois subgroup of the equivalent Galois extension $\widetilde L_{v_i,m_i}$ \resp{$\widetilde L_{\o v_i,m_i}$}.
\vskip 11pt 

So, the Galois subgroup associated with the $v_i$-th \resp{$\o v_i$-th} infinite pseudo-ramified real place will be given by:
\begin{align*}
\Gal(\widetilde L_{v_i}\big/ k)
&= \Gal^D(\widetilde L_{v_i}\big/ k)\bigoplus_{m_i} \Gal(\widetilde L_{v_i,m_i}\big/ k)\\[11pt]
\text{(resp.}\quad
\Gal(\widetilde L_{\o v_i}\big/ k)
&= \Gal^D(\widetilde L_{\o v_i}\big/ k)\bigoplus_{m_i} \Gal(\widetilde L_{\o v_i,m_i}\big/ k)\ ).\end{align*}
For the corresponding pseudo-unramified Galois extensions, we should have:
\begin{align*}
\Gal(\widetilde L ^{nr} _{v_i}\big/ k)
&= \Gal^D(\widetilde L^{nr} _{v_i}\big/ k)\bigoplus_{m_i} \Gal(\widetilde L^{nr} _{v_i,m_i}\big/ k)\\[11pt]
\text{(resp.}\quad
\Gal(\widetilde L^{nr} _{\o v_i}\big/ k)
&= \Gal^D(\widetilde L^{nr} _{\o v_i}\big/ k)\bigoplus_{m_i} \Gal(\widetilde L^{nr} _{\o v_i,m_i}\big/ k)\ ).\end{align*}
On the other hand, the Galois subgroup of the irreducible central extension $\widetilde L_{v^1_i}$ \resp{$\widetilde L_{\o v^1_i}$}, corresponding to the irreducible completion $L_{v^1_i}$ \resp{$L_{\o v^1_i}$} having a rank $N$, is obviously the {\bbf global inertia subgroup\/} $I_{L_{v_i}}$ \resp{$I_{L_{\o v_i}}$} 
which can be defined by:
\begin{align*}
I_{L_{v_i^1}}=\Gal(\widetilde L_{v_i^1}\big/ k) \Big / \Gal(\widetilde L^{nr}_{v_i^1}\big/ k)\\
\text{(resp.}\quad
I_{L_{\o v_i^1}}=\Gal(\widetilde L_{\o v_i^1}\big/ k) \Big / \Gal(\widetilde L^{nr}_{\o v_i^1}\big/ k)\ )\end{align*}
or by the equivalent exact sequence:
\[
\bM
& 1 & \To & I_{L_{v_i^1}} & \To & \Gal(\widetilde L_{v_i}\big/ k) &\To &\Gal(\widetilde L^{nr}_{v_i}\big/ k)
&\To & 1 \quad \\
\text{(resp.} & 1 & \To & I_{L_{\o v_i^1}} & \To & \Gal(\widetilde L_{\o v_i}\big/ k) &\To &\Gal(\widetilde L^{nr}_{\o v_i}\big/ k)
&\To & 1 \ ). \eM\]
As the global inertia subgroups are of Galois type, they are all isomorphic:
\[
\bM
I_{L_{v^1_1}} & \simeq \;\cdots \; \simeq & 
I_{L_{v_p^1}} & \simeq \; \cdots \; \simeq & 
I_{L_{v_i^1}} & \simeq \; \cdots \; \simeq & 
I_{L_{v_t^1}}\;.
\eM\]
Let $\wt L_L$ \resp{$\wt L_R$} denote the union of all finite abelian extensions of $k$.  Then, we have that:
\Bi
\item \quad $\ds \Gal(\wt L_L\big/k) = \bigoplus^t_{i=1}\Gal(\widetilde L_{v_i}\big/k)$,
\item \quad $\ds \Gal(\wt L_R\big/k) = \bigoplus^t_{i=1}\Gal(\widetilde L_{\o v_i}\big/k)$,
\item \quad $\ds \Gal(\wt L^{nr}_L\big/k) = \bigoplus^t_{i=1}\Gal(\widetilde L^{nr}_{v_i}\big/k)$,
\item \quad $\ds \Gal(\wt L^{nr}_R\big/k) = \bigoplus^t_{i=1}\Gal(\widetilde L^{nr}_{\o v_i}\big/k)$.
\Ei

In analogy with the $p$-adic case where the Weil group is the Galois subgroup of the elements inducing on the residue field an integer power of a Frobenius element, it was assumed \cite{Pie2} that the {\bbf Weil group in the archimedean case will be the Galois subgroup of the finite pseudo-ramified extensions characterized by extension degrees $d=0\mod N$\/}.

In this respect, if $\dot{\wt L}_{v_i}$ 
(and $\dot{\wt L}_{\o v_i}$) denotes a pseudo-ramified real Galois extension of degree\linebreak 
$[\dot{\wt L}_{v_i}:k]=i\cdot N$, the Weil group
$W_{\wt L_{v_\oplus}}$
\resp{$W_{\wt L_{\o v_\oplus}}$},
corresponding to the Galois group
$\Gal(\wt L_L\big /k)$
\resp{$\Gal(\wt L_R\big / k)$}, will be given by:
\[ W_{\wt L_{v_\oplus}} = \bigoplus^t_{i=1} \Gal (\dot{\wt L}_{v_i}\big /k) \qquad
\rresp{W_{\wt L_{\o v_\oplus}} = \bigoplus^t_{i=1} \Gal (\dot{\wt L}_{\o v_i}\big /k)}.\]
\vskip 11pt 

\subsection*{2.3 Non abelian global class field concepts}  

The set of \lr real pseudo-ramified completions is, in fact, isomorphic to a one-dimensional \lr affine scheme $\SS^1_L$ \resp{$\SS^1_R$}.  So, the challenge consists in introducing the $n$-dimensional analog of $\SS^1_L$ \resp{$\SS^1_R$} which is a 
$n$-dimensional linear algebraic group.  
But, as 
the endomorphisms $\End_k(A)$ of a $k$-algebra $A$ can be handled throughout its enveloping algebra $A^e=A\otimes_k A^{\rm op}$, where $A^{\rm op}$ denotes the opposite algebra of $A$, because $A^e\simeq \End_k(A)$, and as fundamental algebras, as the algebra of
modular forms, are intrinsically defined in the upper half space, {\bbf bilinearity instead of linearity will be envisaged \cite{Pie2}, \cite{Pie4}\/}.
\vskip 11pt 

Then, the $n$-dimensional equivalent of the product $\SS^1_R\times \SS^1_L$ of the one-dimensional affine schemes  $\SS^1_R$ and $\SS^1_L$ is a $n^{(2)}$-dimensional bilinear algebraic semigroup 
$G^{(n)}( {L_{\o v}}\times {L_{v}})$   isomorphic to the product 
$\GL_n ( {L_{\o v}}\times {L_v}) \equiv T^t_n( {L_{\o v}}) \times T_n({L_v}) $ of the group 
$T^t_n ( {L_{\o v}} )$ of lower triangular matrices with entries in ${L_{\o v}} $ by the group 
$T_n (  {L_v}) $ of upper triangular matrices with entries in
${L_v}$ where $L_v=\{L_{v_{i,m_i}}\}^t_{i=1}$ and
 $L_{\o v}=\{L_{v_{i,m_i}}\}^t_{i=1}$.
\vskip 11pt 

As the algebraic bilinear semigroup $G^{(n)} ( {L_{\o v}}\times {L_v}) $ is constructed over
$ {L_{\o v}}\times {L_v}$, it can be decomposed into $t$ conjugacy classes, $1\le i\le t$, having multiplicities $m^{(t)}=\sup (m_t)$, $m_t\in \NN$, in such a way that $m^{(t)}$ denotes the number of equivalent representatives in the $t$-th conjugacy class.
\vskip 11pt 

{\bbf The algebraic representation of the algebraic bilinear semigroup $\GL_n (  {L_{\o v}}\times {L_v}) $ in the
$G^{(n)} ( {L_{\o v}}\times {L_v}) $-bisemimodule $M_R\otimes M_L$ results from the morphism from 
$\GL_n (  L_{\o v}\times {L_v}) $ into $\GL(M_R\otimes M_L)$ where $\GL(M_R\otimes M_L)$ is the group of automorphisms of $M_R\otimes M_L$.

\noindent So, $\GL(M_R\otimes M_L)$ becomes the $n$-dimensional equivalent of the product 
$W^{ab}_{\wt L_{\o v_\oplus}}\times 
W^{ab}_{\wt L_{v_\oplus}}$ of the global Weil groups and the $n$-dimensional bilinear algebraic semigroup 
 $G^{(n)} (  {L_{\o v}}\times  {L_v}) $ is the $n$-dimensional representation space of 
$W^{ab}_{\wt L_{\o v_\oplus}}\times 
W^{ab}_{\wt L_{v_\oplus}}$\/}.
\vskip 11pt 

Referring to the algebraic bilinear semigroup of matrices $\GL_n (  {L_{\o v}}\times  {L_v}) \equiv T^t_n (  {L_{\o v}})\times T_n(  {L_v}) $, we see that it is submitted to the following {\bbf Gauss bilinear decomposition\/}:
\[ \GL_n (  {L_{\o v}}\times  {L_v}) 
= [D_n (  {L_{v}})\times D_n(  {L_{\o v}}) ]\times [UT_n (  {L_{v}})\times UT^t_n( {L_{\o v}}) ]
\]
where:
\Bi
\item $D_n(\cdot)$ is the group of diagonal matrices of order $n$, also called the $n$-split Cartan subgroup;
\item $UT_n(\cdot)$ is the group of upper unitriangular matrices.
\Ei
In fact, the diagonal bilinear algebraic semigroup of matrices $D_n(  {L_{\o v}}\times  {L_v}) 
\equiv D_n(  {L_{\o v}})\times D_n( {L_v}) $ is more exactly $D_n ( {L_{\o v_D}}\times {L_{v_D}}) $ where 
$L_{\o v_D}$ and $L_{v_D}$ are given respectively by
$\wt L_{\o v_D}=\{L_{\o v_i}\}^t_{i=1}$ and
$\wt L_{ v_D}=\{L_{ v_i}\}^t_{i=1}$ with $m^{(i)}=1$, $1\le i\le t$.

Due to the Gauss bilinear decomposition of the algebraic bilinear semigroup $\GL_n (  {L_{\o v}}\times  {L_v}) $, its conjugacy classes can be partitioned into:
\Bi
\item diagonal conjugacy classes whose representatives $g^{(n)}\RL[i,m^{(i)}=1]$ refer to the representation space
 $\Repsp (D_n (  {L_{\o v}}\times  {L_v}) )$ of the diagonal bilinear algebraic semigroup 
$D_n (  {L_{\o v}}\times  {L_v}) $;
\item off-diagonal conjugacy classes whose representatives $g^{(n)}\RL[i,m^{(i)}>1]$ are generated from the nilpotent biaction of 
$UT^t_n ( {L_{\o v}})\times UT_n( {L_v}) $ on the diagonal conjugacy classes\linebreak $g^{(n)}\RL[i, m^{(i)}=1]$.
\Ei
So, {\bbf the conjugacy class representatives $g^{(n)}\RL[i,m_i]$ of the representation space\linebreak 
$\Repsp(\GL_n (  {L_{\o v}}\times  {L_v}) )$ are the $G^{(n)}( {L_{\o v_i,m_i}}\times L_{v_i,m_i}) $-subbisemimodules $M_{\o v_i,m_i}\otimes M_{v_i,m_i}$ of the algebraic representation space of $\GL_n (  {L_{\o v}}\times  {L_v}) $ given by the 
$G^{(n)} (  {L_{\o v}}\times  {L_v}) $-bisemimodule $M_R\otimes M_L$\/}.  So, we have that:
\[\Repsp (\GL_n ( {L_{\o v_\oplus}}\times  {L_{v_\oplus}}))
= \bigoplus^t_{i=1} \bigoplus_{m_i} g^{(n)}\RL[i,m_i]\]
such that:
\Bi
\item \quad $\ds \Repsp (\GL_n (  {L_{\o v}}\times  {L_v}) )\equiv M_R\otimes M_L$;
\item \quad $\ds M_{R_\oplus}\otimes M_{L_\oplus} =\bigoplus^t_{i=1}\bigoplus_{m_i} (M_{\o v_i,m_i}\otimes M_{v_i,m_i})$,
\item \quad $\ds g^{(n)}\RL[i,m_i]\equiv M_{\o v_i,m_i}\otimes M_{v_i,m_i}$, $1\le i\le t$.
\Ei
\vskip 11pt 

Then, we can state the following propositions:
\vskip 11pt 

\subsection*{2.4 Proposition}

{\em Let $G^{(n)} ( {L_{\o v}}\times  {L_v}) $ denote a $n^{(2)}$-dimensional bilinear algebraic semigroup  isomorphic to the bilinear algebraic semigroup of matrices $\GL_n (  {L_{\o v_\oplus}}\times  {L_{v_\oplus}}) $
.

\noindent Then, $\GL_n (  {L_{\o v}}\times {L_v}) $, having the bilinear Gauss decomposition
\[\GL_n (  {L_{\o v}}\times {L_v}) 
=[D_n (  {L_{\o v_D}}\times  {L_{v_D}}) ] \times [UT_n^t ( {L_{\o v}})\times UT_n( {L_v})] \;, \]
is such that its algebraic representation space $\Repsp (\GL_n (  {L_{\o v_\oplus}}\times  {L_{v_\oplus}}) )$, which is a\linebreak $G^{(n)} (  {L_{\o v}}\times  {L_v}) $-bisemimodule $M_{R_\oplus}\otimes M_{L_\oplus}$, decomposes into diagonal and off-diagonal conjugacy class representatives $g^{(n)}\RL[i,m_i]$ which are $G^{(n)}(L_{\o v_i,m_i}\times L_{v_i,m_i})$-subbisemimodules $M_{\o v_i,m_i}\otimes M_{_i,m_i}$:
\[ M_{R_\oplus}\otimes M_{L_\oplus} = \bigoplus_i\bigoplus_{m_i} (M_{\o v_i,m_i}\otimes M_{v_i,m_i})\]
where:\[ M_R\otimes M_L\equiv \Repsp(\GL_n (  {L_{\o v}}\times  {L_v}) )\;.\]
}\vskip 11pt 

\begin{proof} results from section 2.3.\end{proof}

\subsection*{2.5 Proposition}

{\em Let $I _{L_{v_i^1}}$ and $I_{L_{v_i^1,m_i}}$ be two global inertia subgroups as introduced in section 2.2 and leading to:
\[  I_{L_v} =\bigoplus_iI _{L_{v_i^1}}\bigoplus_{m_i} I_{L_{v_i^1,m_i}} \qquad 
\text{(resp.} \quad
I_{L_{\o v}} =\bigoplus_iI_{L_{\o v_i^1}}\bigoplus_{m_i} I_{L_{\o v_i^1,m_i}} \ ).\]
Let
\[
\Gal(\wt L _L\big/ k)  = \bigoplus^t_{i=1}\Gal (\widetilde L_{v_i}\big/ k)\qquad
\resp{\Gal(\wt L _R\big/ k)  = \bigoplus^t_{i=1}\Gal (\widetilde L_{\o v_i}\big/ k)}\]
be the Galois groups of all finite abelian extensions of $k$.

\noindent Then, we get the {\bbf explicit $n$-dimensional representation spaces:
\Bi
\item \quad $\ds I_{L_{\o v}}\times I_{L_v} \To P^{(n)} ( {L^1_{\o v}}\times  {L^1_v})$,
\item \quad $\ds \Gal (\wt L_R\big/ k) \times \Gal (\wt L_L\big/ k)
\To G^{(n)} ( {L_{\o v}}\times  {L_v})$,
\Ei
where 
\Bi
\item $P^{(n)} ( {L^1_{\o v}}\times  {L^1_v})$ is the bilinear parabolic subgroup,
\item $L_{\o v^1}=\{L_{\o v^1_{i,m_i}}\}^t_{i=1}$,
\Ei
}
with $L_{\o v^1_i,m_i}$ an irreducible central $k$-subsemimodule of rank $N$.
}
\vskip 11pt 

\begin{proof}
\Be
\item The representation
\[  \Gal (\wt L _R\big/ k) \times \Gal (\wt L _L\big/ k) 
\To \GL_n ( {L_{\o v}}\times {L_v})\]
 results from non abelian class field concepts introduced in section 2.3 and leads to  the morphism from $\GL_n ( {L_{\o v}}\times  {L_v})$ into $\GL(M_R\otimes M_L)$.

\item $P^{(n)} ( {L^1_{\o v}}\times {L^1_v})$, being classically defined as the connected component of the identity in $G^{(n)}( {L_{\o v}}\times  {L_v})$,  constitutes a $n$-dimensional representation of $I_{L_{\o v}}\times I_{L_v}$.\qedhere
\Ee
\end{proof}
\vskip 11pt 

\subsection*{2.6 Proposition}
{\em Each $G^{(n)} ({L_{\o v_i,m_i}}\times {L_{v_i,m_i}})$-subbisemimodule $M_{\o v_i,m_i}\otimes M_{v_i,m_i}$ is characterized by a rank
\[ r^{(n)}_{\o v_i\times v_i} =i^{n^2}\cdot N^{n^2}\;.\]
}
\vskip 11pt 

\begin{proof}
 As $(M_{\o v_i,m_i}\otimes M_{v_i,m_i})$ is the $n$-dimensional analog of the one-dimensional bilinear algebraic subsemigroup $L_{\o v_i,m_i}\times L_{v_i,m_i}$ having a rank given by $r^{(1)}_{\o v_i\times v_i}=i^2\cdot N^2$ according to section 2.1, it is clear, according to the non abelian class field concepts developed in \cite{Pie4}, that the rank of $M_{\o v_i,m_i}\otimes M_{v_i,m_i}$ is given by:
\[
r^{(n)}_{\o v_i\times v_i}= (f_{\o v_i})^n\cdot N^n\cdot (f_{v_i})^n\cdot N^n=i^{n^2}\cdot N^{n^2}\;.\qedhere
\]
\end{proof}
\vskip 11pt 

\subsection*{2.7 Lattices and bilattices}  

Let $B_v$ and $B_{\o v}$ be two division semialgebras of dimension $n$ respectively over $L_v$ and $L_{\o v}$ such that $B_{\o v}$ be the opposite division semialgebra of $B_v$.

\noindent If we fix the isomorphisms:
\[ B_v\approx T_n( {L_v}) \quad \text{and} \quad B_{\o v}\approx T^t_n( {L_{\o v}})\;,\]
we have the following isomorphism:
\[ B_{\o v}\otimes B_v \approx T^t_n( {L_{\o v}})\times T_n( {L_v})\]
for the division bisemialgebra $B_{\o v}\otimes B_v$.

\noindent On the other hand, fix the maximal orders $\Os_{L,v}$ of $L_v$ and $\Os_{L,\o v}$ of $L_{\o v}$.  

\noindent Then, the maximal orders $\Lambda _v$ and $\Lambda _{\o v}$ respectively in the division semialgebras $B_v$ and $B_{\o v}$ are pseudo-ramified $\ZZ\big/ N\ \ZZ$-lattices in the $B_v$-semimodule $M_L$ and in the
$B_{\o v}$-semimodule $M_R$.
\vskip 11pt 

\subsection*{2.8 Proposition}

{\em Let $\Lambda _v$ and $\Lambda _{\o v}$ be pseudo-ramified $\ZZ\big/N\ \ZZ$-lattices respectively in the division semialgebras 
$B_v$ and $B_{\o v}$.

\noindent Then, the pseudo-ramified bisemilattice $\Lambda _{\o v}\otimes\Lambda _{v}$  in the 
$ B_{\o v}\otimes B_{v}$-bisemimodule
$M_R\otimes M_L$:
\Bi
\item verifies: 
\[   \Lambda _{\o v}\otimes\Lambda _{v}\simeq \GL_n(\Os_{L,\o v}\times \Os_{L,v})\;;\]
\item has the decomposition:
\[ \Lambda _{\o v}\otimes\Lambda _{v}=\bigoplus_i \bigoplus_{m_i} (\Lambda _{\o v_i,m_i}\otimes\Lambda _{v_i,m_i})\;.\]
\Ei
}
\vskip 11pt 

\begin{proof}
\Be
\item As $\Os_{L,v}$ \resp{$\Os_{L,\o v}$} is a maximal order in $L_v$ \resp{$L_{\o v}$}, and as 
$\Lambda _v$ \resp{$\Lambda _{\o v}$} is a maximal order in the division semialgebra $B_v$ \resp{$B_{\o v}$}, we have that:
\[ \Lambda _v\approx T_n(\Os_{L,v})\qquad \text{(resp.} \quad
 \Lambda _{\o v}\approx T^t_n(\Os_{L,\o v})\ )\]
leading to:
\[ \Lambda _{\o v}\otimes  \Lambda _v \simeq  T^t_n(\Os_{L,\o v})\times T_n(\Os_{L,v})\;.\]

\item As $ L_{\o v}\otimes   L_v $ can be decomposed into a sum of products of pseudo-ramified completions according to:
\[    L_v\otimes   L_v =\bigoplus_i \L( L_{\o v_i}\otimes L_{v_i}\R)\bigoplus_{m_i} \L(L_{\o v_i,m_i}\otimes L_{v_i,m_i}\R)\]
(see section 2.1) and, as $\Lambda _v$ and $\Lambda _{\o v}$ are supposed to be pseudo-ramified $\ZZ\big/N\ \ZZ$-lattices respectively in $B_v$ and $B_{\o v}$, we can conclude that $\Lambda _{\o v}\otimes \Lambda _v$ has the decomposition
\[ \Lambda _{\o v}\otimes \Lambda _v=\bigoplus_i\bigoplus_{m_i}(\Lambda _{\o v_i,m_i}\otimes \Lambda _{v_i,m_i})\]
where
\[ \Lambda _{v_i,m_i}\qquad \text{(resp.} \quad
 \Lambda _{\o v_i,m_i}\ )\]
is a pseudo-ramified sublattice in the $B_{v_i,m_i}$-subsemi\-module $M_{v_i,m_i}$ \resp{$B_{\o v_i,m_i}$-sub\-semimodule $M_{\o v_i,m_i}$}.\qedhere
\Ee
\end{proof}
\vskip 11pt 

\subsection*{2.9 Proposition}

{\em Let $\GL_n( (\ZZ\big/N\ \ZZ)^2 )$ be the general bilinear semigroup of matrices of order $n$ with entries in $(\ZZ\big/ N\ \ZZ)^2 $.

\noindent Then, the Hecke bisemialgebra of dimension $n^2$,  $\Hs\RL(n^2)$, is generated by all the  Hecke bioperators $T_R(n;t)\otimes T_L(n;t)$ having a representation in the subgroup of matrices\linebreak $\GL_n( (\ZZ\big/ N\ \ZZ)^2 )$.
}
\vskip 11pt 

\begin{proof}
\Be
\item A \lr Hecke operator $T_L(n;t)$ \resp{$T_R(n;t)$} is a \lr correspondence which associates to the \lr lattice $\Lambda _v$ \resp{$\Lambda _{\o v}$} the sum of its \lr sublattices $\Lambda _{v_{i,m_i}}$ \resp{$\Lambda _{\o v_{i,m_i}}$} of index $t$ and multiplicities $m^{(i)}=\sup(m_i)$:
\begin{align*}
T_L(n;t)\ \Lambda _v &= \bigoplus_{i=1}^t\bigoplus_{m_i} \Lambda _{v_i,m_i}	\\
\text{(resp.}\quad
T_R(n;t)\ \Lambda _{\o v} &= \bigoplus_{i=1}^t\bigoplus_{m_i} \Lambda _{\o v_i,m_i}	\ ).\end{align*}
So, the Hecke bioperator $T_R(n;t)\otimes T_L(n;t)$ is defined by the bicorrespondence:
\[ (T_R(n;t)\otimes T_L(n;t))(\Lambda _{\o v}\otimes \Lambda _v)
=\bigoplus_i\bigoplus_{m_i} (\Lambda _{\o v_i,m_i}\otimes \Lambda _{v_i,m_i})\;.\]

\item As 
\begin{align*}
&\Lambda _{\o v_i,m_i}\otimes \Lambda _{v_i,m_i}\simeq g_n(\Os_{L_{\o v_i,m_i}}\times \Os_{L_{v_i,m_i}})\\
&\qquad \in \GL_n(\Os_{L,\o v}\times \Os_{L,v})\quad \subset \quad \GL_n( (\ZZ\big/ N\ \ZZ)^2 )\;,\end{align*}
$g_n(\Os_{L_{\o v_i,m_i}}\times \Os_{L_{v_i,m_i}})$ can be chosen as a coset representative of the tensor product 
$T_R(n;t)\linebreak \otimes T_L(n;t)$ of Hecke operators \cite{Pie1}.\qedhere
\Ee
\end{proof}
\vskip 11pt 

\subsection*{2.10 Corollary}

{\em There exists an injective morphism:
\[ m_{\Lambda \RL\to M\RL}: \quad \Lambda _{\o v}\otimes \Lambda _v \To M_R\otimes M_L\]
from the pseudo-ramified bisemilattice $\Lambda _{\o v}\otimes \Lambda _v $ to the corresponding $G^{(n)}( {L_{\o v}}\times  {L_v})$-bisemimodule $M_R\otimes M_L$.
}
\vskip 11pt 

\begin{proof}
 Indeed, $\Lambda _{\o v}\otimes \Lambda _v$ is a pseudo-ramified bisemilattice into the $B_{\o v}\otimes B_v$-bisemimodule $M_R\otimes M_L$ according to proposition 2.8.  And thus, the decomposition of $M_R\otimes M_L$ into $G^{(n)}(L_{\o v_i,m_i}\times L_{v_i,m_i})$-subbisemimodules $M_{\o v_i,m_i}\otimes M_{v_i,m_i}$, being in one-to-one correspondence with the conjugacy class representatives of $\GL_n( {L_{\o v}}\times  {L_v})$, results from the similar decomposition of 
$\Lambda _{\o v}\otimes \Lambda _v $ into
\[
\Lambda _{\o v}\otimes \Lambda _v =\bigoplus_i\bigoplus_{m_i}(\Lambda _{\o v_i,m_i}\otimes \Lambda _{v_i,m_i} )\;.\qedhere
\]
\end{proof}
\vskip 11pt 

\subsection*{2.11 Toroidal compactification of lattice bisemispaces}

The space $X=\GL_n(\RR)\big/ \GL_n(\ZZ)$ corresponds to the set of lattices of $\RR^n$.  In this perspective, a \lr pseudo-ramified lattice semispace
\[
 X_L =T_n({L_L})\Big/ T_n(\ZZ\big/ N\ \ZZ)\qquad
\rresp{X_R =T^t_n({L_{R}})\Big/ T^t_n(\ZZ\big/ N\ \ZZ)}\]
 where $L_R$ and $L_L$ are compactified commutative division semirings corresponding to the semifields $\wt L_R$ and $\wt L_L$ introduced in section 2.1,
is introduced in such a way that the cosets of $X_L $ \resp{$X_R$} correspond to the conjugacy classes of $T_n( {L_v})$ \resp{$T_n^t( {L_{\o v}})$}.
\vskip 11pt

A toroidal compactification 
$\gamma_{X^T_L}$ 
\resp{$\gamma_{X^T_R}$}
is envisaged on $X_L$ \resp{$X_R$} in such a way that it corresponds to a projective mapping which can be decomposed into a two step sequence \cite{Pie3}:
\Be
\item the points $P_{a_L}[i,m_i]\in g^{(n)}_L[i,m_i]$ \resp{$P_{a_R}[i,m_i]\in g^{(n)}_R[i,m_i]$} of the conjugacy class representative $g^{(n)}_L[i,m_i]$ \resp{$g^{(n)}_R[i,m_i]$} of 
$T^{(n)}({L_v})\subset G^{(n)}(L_{\o v}\times L_v)$ 
\resp{$T^{(n)}({L_{\o v}})$}
 are mapped onto the origin of 
$ {L_{v}}$ \resp{$ {L_{\o v}}$}.
\item these points $P_{a_L}[i,m_i]$ \resp{$P_{a_R}[i,m_i] $} are then projected symmetrically from the origin of $ {L_v}$ \resp{$ {L_{\o v}}$} into a  connected compact  semivariety which is a $n$-dimensional real semitorus $T^n_L[i,m_i]$ \resp{$T^n_R[i,m_i]$} in ${L^T_v}$ \resp{${L^T_{\o v}}$} where ${L^T_v}$ \resp{${L^T_{\o v}}$} is given by:
\[
{L^T_v}=\{L^T_{v_{i,m_i}} \}^t_{i=1} \qquad
\rresp{L^T_{\o v}=\{L^T_{\o v_{i,m_i}} \}^t_{i=1}}
\]
with $L^T_{v_i,m_i}$ \resp{$L^T_{\o v_i,m_i}$} being a \lr toroidal  completion.
\Ee
{\bbf The toroidal compactification $\gamma _{X^T_L}$ \resp{$\gamma _{X^T_R}$} of the lattice semispace $X_L$ \resp{$X_R$} is thus the projective mapping\/}:
\begin{align*}
\gamma _{X^T_L} : \qquad &X_L\To X^T_L=T_n({L^T_L})\Big/ T_n(\ZZ\big/ N\ \ZZ)\\
\text{(resp.} \quad
\gamma _{X^T_R} : \qquad &X_R\To X^T_R=T^t_n( {L^T_R})\Big/ T^t_n(\ZZ\big/ N\ \ZZ)\ )\end{align*}
sending $X_L$ \resp{$X_R$} into the corresponding toroidal lattice semispace $X_L^T$ \resp{$X_R^T$} such that 
its cosets correspond to the conjugacy class representatives $g^{(n)}_{T_L}[i,m_i]$ \resp{$g^{(n)}_{T_R}[i,m_i]$} of $T^{(n)}( {L^T_v})$ \resp{$T^{(n)}( {L^T_{\o v}})$} which are $n$-dimensional real semitori
where $L^T_R$ 
\resp{$L^T_L$} is the toroidal equivalent of $L_R$ \resp{$L_L$}.
\vskip 11pt

A pseudo-ramified lattice bisemispace
\[ X\RL=X_R\otimes X_L= \GL_n( {L_R}
\times  {L_L} )\Big/\GL_n( (\ZZ\big/N \ \ZZ)^2)\]
can naturally be generated from the pseudo-ramified lattice semispaces $X_R$ and $X_L$.

The toroidal compactification $\gamma _{X^T_R}\times \gamma _{X^T_L}$ then maps the pseudo-ramified lattice bisemispace $X\RL$ into its toroidal equivalent $X^T\RL$ according to:
\[ \gamma _{X^T_R}\times \gamma _{X^T_L}: \qquad X\RL \To X^T\RL\]
where $X^T\RL$ is given by:
\[X^T\RL = \GL_n(  {L^T_R}
\times {L^T_L} )\Big/\GL_n( (\ZZ\big/N \ \ZZ)^2)\;.\]


\subsection*{2.12 Proposition}

{\em Let $X^T\RL=\GL_n(L^T_R\times L^T_L)\big/ \GL_n((\ZZ\big/N\ \ZZ)^2)$ be the toroidal pseudo-ramified lattice bisemispace.

As $X^T\RL$ corresponds to the representation space of
$\GL_n(L^T_{\o v}\times L^T_v)$ given by the algebraic bilinear semigroup $G^{(n)}(L^T_{\o v}\times L^T_v)$, it decomposes according to:
\[
X^T\RL \simeq G^{(n)}(L^T_{\o v_\oplus}\times L^T_{v_\oplus})
=\bigoplus^t_{i=1} \bigoplus_{m_i} g^{(n)}_{T\RL}[i,m_i]\]
where $g^{(n)}_{T\RL}$ is a conjugacy class representative given by the product, right by left, of $n$-dimensional semitori
$g^{(n)}_{T_R}[i,m_i]$ and
$g^{(n)}_{T_L}[i,m_i]$.}
\vskip 11pt

\begin{proof}
Indeed, the algebraic bilinear semigroup
$G^{(n)}(L^T_{\o v}\times L^T_v)$ results from
$G^{(n)}(L_{\o v}\times L_v)$ by the toroidal compactification:
\[
\gamma _{X^T_R}\times \gamma _{X^T_L}: \qquad
G^{(n)}(L_{\o v}\times L_v) \To
G^{(n)}(L^T_{\o v}\times L^T_v) \;.\]
As $G^{(n)}(L_{\o v}\times L_v)$ decomposes into conjugacy class representatives according to proposition 2.4, it is also the case for the algebraic bilinear semigroup
$G^{(n)}(L^T_{\o v}\times L^T_v)$ whose conjugacy class representatives are products, right by left, of $n$-dimensional semitori due to the projective morphism
$\gamma _{X^T_R}\times \gamma _{X^T_L}$.
\end{proof}
\section{Bilinear cuspidal representations}

\subsection{Bilinear global cuspidal representations on infinite real places}

\subsubsection{Bisemisheaf of rings on the bilinear algebraic semigroup $G^{(n)}(L^T_{\o v}\times L^T_v)$}

The set of differentiable smooth functions
$\phi (g^{(n)}_{T_L}[i,m_i])$
\resp{$\phi (g^{(n)}_{T_R}[i,m_i])$}, $1\le i\le t\le \infty $, on the conjugacy class representatives 
$g^{(n)}_{T_L}[i,m_i]$
\resp{$g^{(n)}_{T_R}[i,m_i]$} of the algebraic semigroup
$G^{(n)}(L^T_v)$
\resp{$G^{(n)}(L^T_{\o v})$} are the sections of a semisheaf of rings
$\Phi (G^{(n)}(L^T_v))$
\resp{$\Phi (G^{(n)}(L^T_{\o v}))$}.

According to proposition 2.12, each section
$\phi (g^{(n)}_{T_L}[i,m_i])\subset \Gamma (\Phi (G^{(n)}(L^T_v)))$
\resp{$\phi (g^{(n)}_{T_R}[i,m_i])\linebreak \subset \Gamma (\Phi (G^{(n)}(L^T_{\o v})))$}
is a differentiable function on a $n$-dimensional real semitorus as it will be developed in proposition 3.1.2.

Similarly, the set of differentiable smooth bifunctions
$\phi\RL (g^{(n)}_{T\RL}[i,m_i])=\phi (g^{(n)}_{T_R}[i,m_i])
\otimes\phi (g^{(n)}_{T_L}[i,m_i])$ on the conjugacy class representatives
$g^{(n)}_{T\RL}[i,m_i]$ of the algebraic bilinear semigroup
$G^{(n)}(L^T_{\o v}\times L^T_v)$ is the set
$\Gamma (\Phi (G^{(n)}(L^T_{\o v}\times L^T_v)))$ of bisections of the bisemisheaf of rings
$\Phi (G^{(n)}(L^T_{\o v}\times L^T_v))$ \cite{Pie4}.

On this set 
$\Gamma (\Phi (G^{(n)}(L^T_{\o v}\times L^T_v)))$ of bisections,
{\bbf a global elliptic
$\Gamma (\Phi (G^{(n)}(L^T_{\o v}\times L^T_v)))$- bisemimodule
$\ELLIP_R(n,i,m_i)\otimes_D\ELLIP_L(n,i,m_i)$} will be explicitly constructed in the following proposition.
\vskip 11pt

\subsubsection{Proposition}  

{\em  The functional representation space of the bilinear algebraic  semigroup $G^{(n)}( {L^T_{\o v_\oplus}}\times  {L^T_{v_\oplus}})$ can be given by the product, right by left, of $n$-dimensional global elliptic semimodules $\Ellip_R(n,i,m_i)$ and $\Ellip_L(n,i,m_i)$ according to \cite{Pie1}:
\[\FRepsp[G^{(n)}({L^T_{\o v_\oplus}}\times {L^T_{v_\oplus}})] =\Ellip_R(n,i,m_i)\otimes \Ellip_L(n,i,m_i)\]
where:
\Bi
\item $\Ellip_L(n,i,m_i)=\bigoplus^t_{i=1}\bigoplus_{m_i}\lambda ^{\half}(n,i,m_i)\ e^{2\pi i (i)x}$~, $1\le i\le t\le \infty $~,

$\Ellip_R(n,i,m_i)=\bigoplus^t_{i=1}\bigoplus_{m_i}\lambda ^{\half}(n,i,m_i)\ e^{-2\pi i (i)x}$

is a (truncated) Fourier series with
\Bi
\item[\textbullet] $x=\sum\limits^n_{\beta=1}x_\beta\ \vec e\ ^\beta$ a vector of $\RR^n$~;
\item[\textbullet]  $\lambda(n,i,m_i)=\prod^n_{c=1}\lambda_c(n,i,m_i)$ a product, right by left, of Hecke characters since\linebreak $\lambda_c(n,i,m_i)$ is an eigenbivalue of $g_n(\Os_{L^T_{\o v_i,m_i}}\times \Os_{L^T_{v_i,m_i}})$~;
\Ei
\item each term 
$\lambda ^{\half}(n,i,m_i)\ e^{2\pi i(i)x})$
\resp{$\lambda ^{\half}(n,i,m_i)\ e^{-2\pi i(i)x})$}
of
$\ELLIP_L(n,i,m_i)$
\resp{$\ELLIP_R(n,i,m_i)$} is a section of the semisheaf 
$\Phi  (G^{(n)}(L^T_v))$
\resp{$\Phi  (G^{(n)}(L^T_{\o v}))$}, i.e. a smooth differentiable function on a $n$-dimensional semitorus
$g^{(n)}_{T_L}[i,m_i]$
\resp{$g^{(n)}_{T_R}[i,m_i]$}.
\Ei
} \vskip 11pt

\begin{proof}
\Be
\item According to proposition  2.12, the representation space $\Repsp (
\GL_n ({L^{T}_{\o v_\oplus}} \times {L^{T}_{v_\oplus}} ))$ decomposes according to  its conjugacy classes representatives:
\[
\Repsp (\GL_n ({L^{T}_{\o v_\oplus}} \times {L^{T}_{v_\oplus}} ))
=G^{(n)}(L^T_{\o v_\oplus}\times L^T_{v_\oplus})
= \bigoplus^t_{i=1}\bigoplus_{m_i} g^{(n)}_{T\RL}[i,m_i]\]
where
\[ g^{(n)}_{T\RL}[i,m_i] =
g^{(n)}_{T_R}[i,m_i]\times g^{(n)}_{T_L}[i,m_i] \;.\]
According to section 2.11, $g^{(n)}_{T_L}[i,m_i] $ (resp. $g^{(n)}_{T_R}[i,m_i]$~) is a $n$-dimensional real semitorus localized in the upper (resp. lower) half space.
\item The decomposition of  $\Repsp (\GL_n ({L^{T}_{\o v}} \times {L^{T}_{v}} ))
$ into conjugacy class representatives\linebreak $g^{(n)}_{T\RL}[i,m_i]$ results from an endomorphism of
$\Repsp (\GL_n ( {L^{T}_{\o v}} \times  {L^{T}_{v}} ))
$ into itself generated by the action of  Hecke bioperators $T_R(n;t)\otimes T_L(n;t)$ \cite{Pie1}.
\item Each smooth continuous function on the \lr conjugacy class representative $g^{(n)}_{T_L}[i,m_i] $ (resp. $g^{(n)}_{T_R}[i,m_i]$~) is (a function on) a $n$-dimensional real semitorus $T^n_L[i,m_i]$ (resp. $T^n_R[i,m_i]$~) which has the following analytic representation:
\begin{align*}
T^n_L[i,m_i]
&\simeq \lambda ^{\half}(n,i,m_i)\ e^{2\pi i(i)x}\\
\text{(resp.} \quad 
T^n_R[i,m_i]
&\simeq \lambda ^{\half}(n,i,m_i)\ e^{-2\pi i(i)x}\ ).\end{align*}
Indeed, as $g^{(n)}_{T_L}[i,m_i] $ (resp. $g^{(n)}_{T_R}[i,m_i]$~) is the non abelian equivalent of the toroidal completion $L^T_{v_i,m_i}$ (resp. $L^T_{\o v_i,m_i}$~) according to section 2.3, we have to consider the global Frobenius substitution at the \lr place $v_i$ (resp. $\o v_i$~)  given by the mapping:
\begin{alignat*}{3}
e^{2\pi ix} &\To e^{2\pi i(i)x} \qquad \qquad &&x\in \RR^n\;, \\
\text{(resp.} \quad
e^{-2\pi ix} &\To e^{-2\pi i(i)x}\ ) \qquad \qquad && i=\sqrt{-1}\;,  \end{alignat*}
$(i)\in\nit$ being the global residue degree of this infinite place $v_i$ \resp{$\o v_i$}.

On the other hand, as ${T^n\RL}[i,m_i] $ results from an endomorphism of
$\Repsp (\GL_n ({L^{T}_{\o v}} \times {L^{T}_{v}} ))
$ into itself, the scalar $\lambda (n,i,m_i)$ will correspond to the eigenvalues of the associated coset representative of the product of Hecke operators.

 This coset representative of $T_R(n;t)\otimes T_L(n;t)$ is then given by
$g_n(\Os_{L^T_{\o v_i,m_i}}\times \Os_{L^T_{v_i,m_i}})$ according to proposition 2.9.

 If $\{\lambda _c(n,i,m_i)\}_{c=1}^n$ denotes the set of eigenvalues of
$g_n(\Os_{L^T_{\o v_i,m_i}}\times \Os_{L^T_{v_i,m_i}})$~, then
\[ \lambda (n,i,m_i)=\prod^n_{c=1}\lambda _c(n,i,m_i)\]
can be considered as a product, right by left, of Hecke characters and its square root $\lambda ^{\half}(n,i,m_i)$ can be chosen as the coefficient of $T^n_L[i,m_i]$ (resp. $T^n_R[i,m_i]$~).\qedhere
\Ee
\end{proof}
\vskip 11pt

\subsubsection{Analytic representation of semitori with respect to Hecke characters}
Let
\begin{align*}
 \FRepsp [G^{(n)}(L^T_{\o v_p,m_p}\times L^T_{v_p,m_p})]
&= \ellip\RL(n,[p],m_p)\\
&=
\lambda^{\half}(n,p,m_p)\ e^{-2\pi ipx}\otimes
\lambda^{\half}(n,p,m_p)\ e^{2\pi ipx}\end{align*}
be the analytic representation of the algebraic bilinear subsemigroup $G^{(n)}(L^T_{\o v_p,m_p}\times L^T_{v_p,m_p})$ with respect to the $(p,m_p)$-th completion of $L_R\times L_L$~.\\
 $\FRepsp[G^{(n)}(L^T_{\o v_p}\times L^T_{v_p})]$ has the analytic development $\ellip\RL(n,[p],m_p)$ which is in bijection with the product of a right $n$-dimensional semitorus  by its left equivalent. \vskip 11pt

Indeed, a \lr semitorus $T^n_L[p,m_p]$ (resp. $T^n_R[p,m_p]$~) is diffeomorphic to \[
T^1_{1_L}[p,m_p] \times \cdots \times T^1_{c_L}[p,m_p]\times \cdots \times T^1_{n_L}[p,m_p]\]
\[ \text{(resp.}\quad 
T^1_{1_R}[p,m_p] \times \cdots \times T^1_{c_R}[p,m_p]\times \cdots \times T^1_{n_R}[p,m_p]\ )\]
in such a way that the one-dimensional semitorus $T^1_{c_L}[p,m_p]$ (resp. $T^1_{c_R}[p,m_p]$~) has the representation given by the  following analytic development:
\begin{alignat*}{3}
T^1_{c_L}[p,m_p] &= r_c(p,m_p)\ e^{2\pi ipx_c}\;, \qquad && x_c\in \RR^1\; ,\\
\text{(resp.} \quad 
T^1_{c_R}[p,m_p] &= r_c(p,m_p)\ e^{-2\pi ipx_c}\;, \qquad && x_c\in \RR^1\; ),\end{alignat*}
whose radius $r_c(p,m_p)$ can be expressed with respect to $\lambda_c(p,m_p)$ according to \cite{Pie2}.  So, we have that:
\begin{align*}
\ellip_L(n,[p],m_p) &\approx \prod^n_{c=1} r_c(p,m_p)\ e^{2\pi ipx_c}\\
\text{(resp.} \quad
\ellip_R(n,[p],m_p) &\approx \prod^n_{c=1} r_c(p,m_p)\ e^{-2\pi ipx_c}\ ).\end{align*} \vskip 11pt

\subsubsection{Proposition}  

{\em Let $\ellip\RL(n,[p],m_p)\in \Ellip\RL(n,i,m_i)$ be a global elliptic
$(L^T_{\o v_p,m_p}\times L^T_{v_p,m_p})$-subbisemi\-module.\newline
Then,
\begin{align*}
\ellip\RL(n,[p],m_p) &= \lambda^{\half}(n,p,m_p)\ e^{-2\pi ipx}\otimes \lambda^{\half}(n,p,m_p)\ e^{2\pi ipx}\;, && x\in\rit^n\;,\\
&\approx\prod^n_{c=1}T^1_{c_R}[p,m_p]\otimes T^1_{c_L}[p,m_p]\\
&= \prod^n_{c=1} r_c(p,m_p)\ e^{-2\pi ipx_c}\otimes r_c(p,m_p)\ e^{2\pi ipx_c}\end{align*}
is the analytic representation of $G^{(n)}(L^T_{\o v_p,m_p}\times L^T_{v_p,m_p})$~.
} \vskip 11pt


\begin{proof}[Sketch of proof]

In fact,
\[ \ellip\RL(n,[p],m_p) =
\lambda ^\half(n,p,m_p)\ e^{-2\pi ipx}\otimes
\lambda ^\half(n,p,m_p)\ e^{2\pi ipx}\]
is a deformation of the product, right by left, of $n$-dimensional semitori:
\[ T^n_R[p,m_p]\otimes T^n_L[p,m_p]
\simeq \prod^n_{c=1}\L(
T^1_{c_R}[p,m]\times T^1_{c_L}[p,m_p]\R)\]
resulting from the isomorphism:
\[ I_{\rm {EL}\to T}: \quad
\ellip\RL(n,[p],m^p)
\overset{\sim}{\To} T^n_R[p,m_p]\otimes T^n_L[p,m_p]\]
sending $\lambda ^\half(n,p,m_p)$ into
$\ds\prod^n_{c=1}
r_c(p,m_p)$.\end{proof}
\vskip 11pt

\subsubsection{Proposition (Langlands global correspondence)}

{\em The functional representation space of the toroidal compactification of
$G^{(n)}(L_{\o v_\oplus}\otimes L_{v_\oplus})$:
\[\gamma _{G^{(n)}_{T\RL}}: \qquad
G^{(n)}(L_{\o v_\oplus}\times L_{v_\oplus})\To
G^{(n)}(L^T_{\o v_\oplus}\times L^T_{v_\oplus})\;, \]
where $\gamma _{G^{(n)}_{T\RL}}$ is an isomorphism, leads to the morphism
\[ G^{(n)}(L_{\o v_\oplus}\times L_{v_\oplus})
\To \ELLIP\RL(n,i,m_i)\]
which is equivalent to the following Langlands correspondence:
\[ \Irr\Rep^{(n)}(W_{\wt L_{\o v_\oplus}}\times W_{\wt L_{v_\oplus}})
\overset{\sim}{\To}
\Irr\cusp(\GL_n(L_{\o v}\times L_v))\]
where
\Bi
\item $\Irr\Rep^{(n)}(W_{\wt L_{\o v_\oplus}}\times W_{\wt L_{v_\oplus}})$ is the sum of products, right by left, of the equivalence classes of the irreducible $n^2$-dimensional representation of the bilinear global Weil group\linebreak 
$(W_{\wt L_{\o v }}\times W_{\wt L_{v }})$;

\item $\Irr\cusp(\GL_n(L_{\o v}\times L_v))$ is the sum of the products, right by left, of the conjugacy classes of the irreducible cuspidal representation space of 
$\GL_n(L_{\o v}\times L_v)$.
\Ei}
\vskip 11pt

\begin{proof}
\Be
\item Consider the sequence of morphisms:
\[
G^{(n)}(L_{\o v_\oplus}\times L_{v_\oplus}) \overset{\gamma _{G^{(n)}_{T\RL}}}{\To}
G^{(n)}(L^T_{\o v_\oplus}\times L^T_{v_\oplus}) \overset{\FRepsp ( {G^{(n)}_{T\RL}})}{\To}
\ELLIP\RL(n,i,m_i)\]
where $\FRepsp({G^{(n)}_{T\RL}})$, being the functional representation of the algebraic bilinear semigroup
$G^{(n)}(L^T_{\o v_\oplus}\times L^T_{v_\oplus})$, is the global elliptic bisemimodule $\ELLIP\RL(n,i,m_i)=\linebreak \ELLIP_R(n,i,m_i)
\otimes \ELLIP_L(n,i,m_i)$ according to proposition 3.1.2.

$G^{(n)}(L_{\o v_\oplus}\times L_{v_\oplus})$ is then clearly in bijection with   $\ELLIP\RL(n,i,m_i)$;

\item $G^{(n)}(L_{\o v_\oplus}\times L_{v_\oplus})$ is the $n$-dimensional representation space of
$(W^{ab}_{L_{\o v_\oplus}}\times W^{ab}_{L_{v_\oplus}})$ according to section 2.3.

On the other hand, $\ELLIP\RL(n,i,m_i)$ constitutes a cuspidal representation space of
$\GL_n(L_{\o v }\times L_{v })$ as developed in \cite{Pie1}.

\item the searched bijection
$\Irr\Rep^{(n)}(W_{\wt L_{\o v_\oplus}}\times W_{\wt L_{v_\oplus}})
\overset{\sim}{\To}
\Irr\cusp(\GL_n(L_{\o v}\times L_v))$ thus follows.  \qedhere
\Ee
\end{proof}
\vskip 11pt

\subsection{Bilinear global cuspidal representations (on infinite real places) covered by $p^\ell$-th roots}

\subsubsection{Proposition (Covering of one-dimensional semitori (i.e. semicircles) by $p^\ell$-th roots)}

{\em Let $p$ be a prime number, $v_p$ \resp{$\o v_p$} the $p$-th primary real infinite place and 
$v_{p+h}$
\resp{$\o v_{p+h}$} the $h$-th real infinite place above $v_p$
\resp{$\o v_p$}.

Then, {\bbf the \lr one-dimensional semitorus
$T^1_{c_L}[p+h,m_{p+h}]$
\resp{$T^1_{c_R}[p+h,m_{p+h}]$} can be covered by the $p^\ell$-th complex roots if we introduce the (etale) covering map\/}:
\begin{align*}
 R^1_{L_{p^\ell\to p+h}}:\qquad
&T^1_{c_L}[p+h,m_{p+h},\o p^\ell] = r_c(p+h,m_{p+h},\o p^\ell)^{p^\ell}\ e^{2\pi i(p^\ell)x'_c} \\
& \quad \To T^1_{c_L}[p+h,m_{p+h}] = r_c(p+h,m_{p+h})\ e^{2\pi i(p+h)x_c} && x_c\in\rit\;, \\
\text{(resp.} \quad R^1_{R_{p^\ell\to p+h}}:\qquad
&T^1_{c_R}[p+h,m_{p+h},\o p^\ell] = r_c(p+h,m_{p+h},\o p^\ell)^{p^\ell}\ e^{-2\pi i(p^\ell)x'_c} \\
& \quad \To T^1_{c_R}[p+h,m_{p+h}] = r_c(p+h,m_{p+h})\ e^{-2\pi i(p+h)x_c} \ )
\end{align*}
in such a way that:
\begin{align*}
T^1_{c_L}[p+h,m_{p+h}] &= T^1_{c_L}[p+h,m_{p+h},\o p^\ell] \\
\rresp{T^1_{c_R}[p+h,m_{p+h}] &= T^1_{c_R}[p+h,m_{p+h},\o p^\ell]}\ , \end{align*}
i.e. if:
\Bean
\item $r_c(p+h,m_{p+h})=
r_c(p+h,m_{p+h},\o p^\ell)^{p^\ell}$;

\item $p^\ell x'_c =x_c+2k\pi $, $0\le k\le (p^\ell-1)/2$;

\item $\ds\frac{p^\ell-1}2 \ge (p+h)$.
\Ee
}
\vskip 11pt

\begin{proof}
\Be
\item The etale covering map $R^1_{L_{p^\ell\to p+h}}$
\resp{$R^1_{R_{p^\ell\to p+h}}$} is equivalent to finding complex numbers
$z'_L$ \resp{$z'_R$} (i.e. $p^\ell$-th complex roots) of which
$p^\ell$-th power is equal to $z_L\in T^1_{c_L}[p+h,m_{p+h}]$
\resp{$z_R\in T^1_{c_R}[p+h,m_{p+h}]$}, i.e.
\[ (z'_L)^{p^\ell}=z_L \qquad \rresp{(z'_R)^{p^\ell}=z_R }\]
or
\begin{align*}
r_c(p+h,m_{p+h},\o p^\ell)^{p^\ell}\ e^{2\pi ip^\ell x'_c}&=r_c(p+h,m_{p+h})\ e^{2\pi i(p+h)x_c}\\
\rresp{r_c(p+h,m_{p+h},\o p^\ell)^{p^\ell}\ e^{-2\pi ip^\ell x'_c}&=r_c(p+h,m_{p+h})\ e^{-2\pi i(p+h)x_c}}
\end{align*}
which leads to the conditions a) and b) of this proposition.

\item Condition c) results from the fact that $x_c$ is a point of order $\#\Nu\times N$, and, thus, that to each point $x_c$ correspond $p^\ell/2$ roots $p^\ell$-th.

So, the number of points on $T^1_{c_L}[p+h,m_{p+h}]$
(or $T^1_{c_R}[p+h,m_{p+h}]$) is equal to:
\[ n_{T^1[p+h]}=\#\Nu\times N\times (p+h)\;,\]
where $(p+h)$ is the global residue degree $f_{v_{p+h}}$, $\#\Nu$ is the number of nonunits and $N$ is the degree of an irreducible completion according to section 2.1, while the number of points on
$T^1_{c_L}[p+h,m_{p+h},\o p^\ell]$
(or $T^1_{c_R}[p+h,m_{p+h},\o p^\ell]$) is equal to:
\[ n_{T^1[p+h,\o p^\ell]}=\#\Nu\times N\times (p^\ell-1)/2\;.\]
The covering map
$R^1_{L_{p^\ell\to p+h}}$
\resp{$R^1_{R_{p^\ell\to p+h}}$} is an isomorphism if $(p^\ell-1)/2=p+h$
and  an epimorphism if $(p^\ell-1)/2>p+h$.
\qedhere
\Ee
\end{proof}
\vskip 11pt

\subsubsection{Corollary}

{\em To a semitorus $T^1_{c_L}[p,m_{p}]$ 
\resp{$T^1_{c_R}[p,m_{p}]$} at the $v_p$-th
\resp{$\o v_p$-th} real archimedean place corresponds an etale covering map 
$R^1_{L_{p\to p}}$
\resp{$R^1_{R_{p\to p}}$} by $p$-th complex roots.
}
\vskip 11pt

\begin{proof}
Referring to proposition 3.2.1, we see that this case corresponds to the condition c) with $h=0$ and thus $\ell=1$.  The covering map
$R^1_{L_{p^\ell\to p+h}}$
(and $R^1_{R_{p^\ell\to p+h}}$) of this proposition occurs then at the conditions:
\Bean
\item $r_c(p,m_{p})=
r_c(p,m_{p},\o p)^{p}$;

\item $p x'_c =x_c+2k\pi $, $0\le k\le (p-1)/2$;

\item $({p^\ell}-1)/2=p+h$ with $\ell=1$ and $h=0$.\qedhere
\Ee
\end{proof}
\vskip 11pt

\subsubsection{Corollary}

{\em The semitori $T^1_{c_L}[j,m_j]$ 
\resp{$T^1_{c_R}[j,m_j]$} at the $v_j$-th
\resp{$\o v_j$-th} real archimedean places 
below the $v_p$-th
\resp{$\o v_p$-th} place, i.e. for $j<p$, cannot be covered by $p$-th complex roots.
}
\vskip 11pt

\begin{proof}
Indeed, if $v_j<v_p$
\resp{$\o v_j<\o v_p$}, then there are points of 
$T^1_{c_L}[j,m_j]$ 
\resp{$T^1_{c_R}[j,m_j]$} which are not covered by $p$-th complex roots.
\end{proof}
\vskip 11pt

\subsubsection{Proposition}

{\em {\bbf The $n$-dimensional semitorus $T^n_L[p+h,m_{p+h}]$
\resp{$T^n_R[p+h,m_{p+h}]$}, at the $v_{p+h}$-th
\resp{$\o v_{p+h}$-th} archimedean place, will be covered by $p^\ell$-th complex roots according to\/}:
\begin{align*}
 R^n_{L_{p^\ell\to p+h}}:\qquad
&
T^n_{L}[p+h,m_{p+h},\o p^\ell] \simeq 
\prod^n_{c=1} r_c(p+h,m_{p+h},\o p^\ell)^{p^\ell}\ e^{2\pi i(p^\ell)x'_c}
 \\
& \quad \To T^n_{L}[p+h,m_{p+h}] \simeq 
\prod^n_{c=1} r_c(p+h,m_{p+h})\ e^{2\pi i(p+h)x_c}
\\
\text{(resp.} \quad R^n_{R_{p^\ell\to p+h}}:\qquad
&T^n_{R}[p+h,m_{p+h},\o p^\ell] \simeq
\prod^n_{c=1} r_c(p+h,m_{p+h},\o p^\ell)^{p^\ell}\ e^{-2\pi i(p^\ell)x'_c} \\
& \quad \To  T^n_{R}[p+h,m_{p+h}]\simeq 
\prod^n_{c=1}r_c(p+h,m_{p+h})\ e^{-2\pi i(p+h)x_c} \ )
\end{align*}
leading to the same conditions as these considered in proposition 3.2.1.
}
\vskip 11pt

\begin{proof}
This is evident since this proposition is the $n$-dimensional generalization of proposition 3.2.1.
\end{proof}
\vskip 11pt

\subsubsection{Proposition (Base change under covering by $p^\ell$-th complex roots)}

{\em Let the covering map $R^n_{L_{p^\ell\to p+h}}$
\resp{$R^n_{R_{p^\ell\to p+h}}$} of the semitorus 
$T^n_{L}[p+h,m_{p+h}]$
\resp{$T^n_{R}[p+h,m_{p+h}]$} on the $v_{p+h}$-th
\resp{$\o v_{p+h}$-th} place by the semitorus 
$T^n_{L}[p+h,m_{p+h},\o p^\ell]$
\resp{$T^n_{R}[p+h,m_{p+h},\o p^\ell]$} be an isomorphism.

Then, this covering map by $p^\ell$-th complex roots corresponds to an equivariant base change from a base of dimension $(p+h)^n$ to a covering base of dimension $((p^\ell-1)/2)^n$.
}
\vskip 11pt

\begin{proof}
This results directly from proposition 3.2.1, condition c), and proposition 3.2.4.
\end{proof}
\vskip 11pt

\subsubsection{Corollary}

{\em {\bbf Only the $n$-dimensional global elliptic semimodules\/}
\begin{align*}
\ELLIP_L(n,p\le i,m_i)&= \bigoplus^t_{i=p} \bigoplus_{m_i} \lambda ^\half(n,p\le i,m_i)\ e^{2\pi i(i)x}\;, \\
\ELLIP_R(n,p\le i,m_i)&= \bigoplus^t_{i=p} \bigoplus_{m_i} \lambda ^\half(n,p\le i,m_i)\ e^{-2\pi i(i)x}\;, 
\end{align*}
$x\in\rit^n$, $i\equiv p+h$, $h$ running from $0$ to $\infty $, $p\le i\le t\le\infty $,\\
{\bbf with terms $i\ge p$ can be covered by $p^\ell$-th complex roots, $\ell$ varying.}}
\vskip 11pt

\begin{proof}
Indeed, according to proposition 3.2.1 and corollary 3.2.3, only $n$-dimensional semitori 
$T^n_L[p+h,m_{p+h}]$
\resp{$T^n_R[p+h,m_{p+h}]$} having a rank $(p+h)^n$ can be covered by $p^\ell$-th complex roots.  They thus correspond to terms $i=p+h$, $0\le h\le\infty $.
\end{proof}
\vskip 11pt

\subsubsection{Semisheaves associated with $\ELLIP_L(n,p\le i,m_i)$
and $\ELLIP_R(n,p\le i,m_i)$}

Let $\Phi (G^{(n)}(L^T_v))$
\resp{$\Phi (G^{(n)}(L^T_{\o v}))$} be the semisheaf on the algebraic semigroup\linebreak
$G^{(n)}(L^T_v)$
\resp{$G^{(n)}(L^T_{\o v})$} and let 
$\ELLIP_L(n,i,m_i)$
\resp{$\ELLIP_R(n, i,m_i)$} be the associated global elliptic
$\Gamma (\Phi (G^{(n)}(L^T_v)))$
\resp{$\Gamma (\Phi (G^{(n)}(L^T_{\o v})))$} semimodule in such a way that to each section
$\phi (g^{(n)}_{T_L}[i,m_i])$
\resp{$\phi (g^{(n)}_{T_R}[i,m_i])$} of 
$\Phi (G^{(n)}(L^T_v))$
\resp{$\Phi (G^{(n)}(L^T_{\o v}))$} corresponds one term
$\ellip_L(n,[i],m_i)=\lambda ^\half(n,i,m_i)\ e^{2\pi i(i)x}$
\resp{$\ellip_R(n,[i],m_i)=\lambda ^\half(n,i,m_i)\ e^{-2\pi i(i)x}$} of 
$\ELLIP_L(n,i,m_i)$
\resp{$\ELLIP_R(n, i,m_i)$}.

Let then 
$\Phi (G^{(n)}(L^T_{[v_p]}))$
\resp{$\Phi (G^{(n)}(L^T_{[\o v_p]}))$}
denote the semisheaf on the algebraic semigroup over the set of completions 
$[v_p]=\{v_p,v_{p+1},\dots,v_t\}$
\resp{$[\o v_p]=\{\o v_p,\o v_{p+1},\dots,\o v_t\}$}
superior and equal to $p$ and let
$\ELLIP_L(n,p\le i,m_i)$
\resp{$\ELLIP_R(n,p\le i,m_i)$} be the associated global elliptic semimodule as introduced in corollary 3.2.6.

Then, there exists {\bbf the epimorpism\/}
\begin{align*}
\Em^{(n)}_{v\to[v_p]}: \qquad \Phi (G^{(n)}(L^T_v))&\To\Phi (G^{(n)}(L^T_{[v_p]}))\\
\rresp{\Em^{(n)}_{\o v\to[\o v_p]}: \qquad \Phi (G^{(n)}(L^T_{\o v}))&\To\Phi (G^{(n)}(L^T_{[\o v_p]}))}
\end{align*}
{\bbf of which kernel}
$\ker(\Em^{(n)}_{v\to[v_p]})=\Phi (G^{(n)}(L^T_{v-[v_p]}))$
\resp{$\ker(\Em^{(n)}_{\o v\to[\o v_p]})=\Phi (G^{(n)}(L^T_{\o v-[\o v_p]}))$}
{\bbf is the complementary semisheaf\/}
$\Phi (G^{(n)}(L^T_{v-[v_p]}))$
\resp{$\Phi (G^{(n)}(L^T_{\o v-[\o v_p]}))$}
on the algebraic semigroup
$G^{(n)}(L^T_{ v-[ v_p]})$
\resp{$G^{(n)}(L^T_{\o v-[\o v_p]})$} over the set of completions
$v-[v_p]=\{v_1,\dots,v_{p-1}\}$
\resp{$\o v-[\o v_p]=\{\o v_1,\dots,\o v_{p-1}\}$}.
\vskip 11pt

\subsubsection{Proposition}

{\em {\bbf Each semitorus $T^1_{c_L}[p+h,m_{p+h},\o p^\ell]$
\resp{$T^1_{c_R}[p+h,m_{p+h},\o p^\ell]$} covering the global elliptic subsemimodule\/}
$\ellip_L(1,[p+h],m_{p+h})=\lambda ^\half(1,p+h,m_{p+h})\ e^{2\pi i(p+h)x_c}$
\resp{\linebreak $\ellip_R(1,[p+h],m_{p+h})=\lambda ^\half(1,p+h,m_{p+h})\ e^{-2\pi i(p+h)x_c}$}
{\bbf is a discrete valuation semiring\/} of which:
\Bean
 \item the uniformizing element is
 $r_c(p+h,m_{p+h},\o p^\ell)^p$;
 
 \item the units are the invertible elements 
 $e^{2\pi i(p^\ell)x'_c}$
\resp{$e^{-2\pi i(p^\ell)x'_c}$}
each $x'_c$ verifying\linebreak $p^\ell x'_c=x_c+2k\pi $, $x'_c\in\rit$.
\Ee}
\vskip 11pt

\begin{proof}
\Be
\item Referring to proposition 3.2.1, we see that the semitorus
$T^1_{c_L}[p+h,m_{p+h},\o p^\ell]=r_c(p+h,m_{p+h},\o p^\ell)^{p^\ell}\ e^{2\pi i(p^\ell)x'_c}$
\resp{$T^1_{c_R}[p+h,m_{p+h},\o p^\ell]=r_c(p+h,m_{p+h},\o p^\ell)^{p^\ell}\ e^{-2\pi i(p^\ell)x'_c}$}
covering the global elliptic subsemimodule
$\ellip_L(1,[p+h],m_{p+h})$
\resp{$\ellip_R(1,[p+h],m_{p+h})$} must be a discrete valuation semiring because it is a principal ideal domain having a unique nonzero prime ideal given by:
\begin{align*}
T^1_{c_L}[p+h,m_{p+h},\o p]&=r_c(p+h,m_{p+h},\o p)^{p}\ e^{2\pi i(p)x'_c}\\
\rresp{T^1_{c_R}[p+h,m_{p+h},\o p]&=r_c(p+h,m_{p+h},\o p)^{p}\ e^{-2\pi i(p)x'_c}}.\end{align*}
Each element of this discrete valuation semiring
$T^1_{c_L}[p+h,m_{p+h},\o p^\ell]$
\resp{$T^1_{c_R}[p+h,m_{p+h},\o p^\ell]$} then writes as the product of the $\ell$ power of the uniformizing element
$r_c(p+h,m_{p+h},\o p^\ell)^{p}$
by a unit
$e^{2\pi i(p^\ell)x'_c}$
\resp{$e^{-2\pi i(p^\ell)x'_c}$}.

The valuation of this element is the integer $\ell$.\qedhere
\Ee
\end{proof}
\vskip 11pt

\subsubsection{Corollary}

{\em The etale covering of the global elliptic subsemimodules
$\ellip_L(1,[p+h],m_{p+h})$
\resp{$\ellip_R(1,[p+h],m_{p+h})$}
by the semitori
$T^1_{c_L}[p+h,m_{p+h},\o p^\ell]$
\resp{$T^1_{c_R}[p+h,m_{p+h},\o p^\ell]$} {\rm \cite{Pie2}}
can give rise to other discrete valuation semirings of which
\Bean
\item the uniformizing element is $r_c(p+h,m_{p+h})$;

\item the units are the invertible elements $e^{2\pi ik}$, $0\le k\le (p^\ell-1)$, in one-to-one correspondence with the $p^\ell$-th complex roots $e^{2\pi ik/p^\ell}$ of unity;

\item the unique valuation is the integer 1.
\Ee}
\vskip 11pt

\begin{proof}
Indeed, the semitorus
$T^1_{c_L}[p+h,m_{p+h},\o p^\ell]=r_c(p+h,m_{p+h},\o p^\ell)^{p^\ell}\ e^{2\pi i(p^\ell)x'_c}$
\resp{\linebreak $T^1_{c_R}[p+h,m_{p+h},\o p^\ell]=r_c(p+h,m_{p+h},\o p^\ell)^{p^\ell}\ e^{-2\pi i(p^\ell)x'_c}$}
can also we written according to:
\begin{align*}
T^1_{c_L}[p+h,m_{p+h}]&=r_c(p+h,m_{p+h})\ e^{2\pi i(p+h)x_c}\cdot e^{2\pi ik}\\
\rresp{T^1_{c_R}[p+h,m_{p+h}]&=r_c(p+h,m_{p+h})\ e^{-2\pi i(p+h)x_c}\cdot e^{-2\pi ik}}
\end{align*}
since $p^\ell x'_c=x_c+2k\pi $, $i\le k\le (p^\ell-1)/2$ according to proposition 3.2.1.

Consequently, it is a principal ideal domain of which the unique prime ideal is
$T^1_{c_L}[p+h,m_{p+h}]=r_c(p+h,m_{p+h})\ e^{2\pi i(p+h)x_c}$
\resp{$T^1_{c_R}[p+h,m_{p+h}]=r_c(p+h,m_{p+h})\ e^{-2\pi i(p+h)x_c}$}.

The uniformizing element is $r_c(p+h,m_{p+h})$ or
 $r_c(p+h,m_{p+h})\ e^{2\pi i(p+h)x_c}$ with $x_c=0$ and the units are $e^{2\pi ik}$.
 
 Consequently, the $(p^\ell-1)/2\times \Nu\times N$ points of
 $T^1_{c_L}[p+h,m_{p+h}]$
\resp{$T^1_{c_R}[p+h,m_{p+h}]$} can be expressed from the product of the uniformizing element
$r_c(p+h,m_{p+h})\ e^{2\pi i(p+h)x_c}$ 
\resp{$r_c(p+h,m_{p+h})\ e^{-2\pi i(p+h)x_c}$} at $x_c=0$, by the $p^\ell$ units 
$e^{2\pi ik}$
\resp{$e^{-2\pi ik}$}, $0\le k\le (p^\ell-1)/2$.
\end{proof}
\vskip 11pt

\subsubsection{Corollary}

{\em {\bbf The Kronecker-Weber theorem follows directly from the existence of discrete valuation semirings\/}
$T^1_{c_L}[p+h,m_{p+h}]$
\resp{$T^1_{c_R}[p+h,m_{p+h}]$}, $0\le h\le\infty $, $h\in\NN$.
}
\vskip 11pt

\begin{proof}
Corollary 3.2.8 shows that the points of the valuation semirings
$T^1_{c_L}[p+h,m_{p+h}]$
\resp{$T^1_{c_R}[p+h,m_{p+h}]$} can be expressed by multiplying the uniformizing element
$r_c(p+h,m_{p+h})\linebreak e^{2\pi i(p+h)x_c}$ 
\resp{$r_c(p+h,m_{p+h})\ e^{-2\pi i(p+h)x_c}$} at $x=0$ by the units $e^{2\pi ik}$, $0\le k\le (p^\ell-1)/2$ which are in one-to-one correspondence with the $p^\ell$-th roots of unity 
$e^{2\pi ik/p^\ell}$, solution of the equation $x^{p^\ell}-1=0$.  (These roots form a cyclic group having as generator
$e^{2\pi i/p^\ell}$).

Note that there is an inflation map:
\begin{align*}
\INF_L: \quad e^{2\pi ik}&\To r_c(p+h,m_{p+h})\ e^{2\pi ix_c}\cdot e^{2\pi ik}\\
\rresp{\INF_R: \quad e^{-2\pi ik}&\To r_c(p+h,m_{p+h})\ e^{-2\pi ix_c}\cdot e^{-2\pi ik}}
\end{align*}
from the units, points of a circle having a radius equal to $1$, to the complex numbers 
$r_c(p+h,m_{p+h})\ e^{2\pi ix_c}\ \cdot e^{2\pi ik}$ 
\resp{$r_c(p+h,m_{p+h})\ e^{-2\pi ix_c}\ \cdot e^{-2\pi ik}$},
points of a circle having a radius equal to 
$r_c(p+h,m_{p+h})$.

When the covering map 
$R^1_{L_{p^\ell\to p+h}}$
\resp{$R^1_{R_{p^\ell\to p+h}}$}
is an epimorphism (see proposition 3.2.1), i.e. the case where $(p^\ell-1)/2>p+h$, the finite abelian extension of $k=\QQ$, characterized by a global residue degree $f_{v_{p+h}}=p+h$, is thus related to the cyclotomic field of $p^\ell$-th roots of unity which corresponds to a Galois cyclotomic extension of order $(p^\ell-1)/2$.

The Kronecker-Weber theorem, expressing that each finite abelian extension of $\QQ$ is contained in a cyclotomic extension of $\QQ$, is thus reached here.
\end{proof}
\vskip 11pt

\subsubsection{Proposition}

{\em {\bbf Each $n$-dimensional semitorus 
$T^n_L[p+h,m_{p+h},\o p^\ell]$
\resp{$T^n_R[p+h,m_{p+h},\o p^\ell]$}}
covering the global elliptic subsemimodule
\begin{align*}
\ellip_L(n,[p+h],m_{p+h})& =\lambda ^\half(n,p+h,m_{p+h})\ e^{2\pi i(p+h)x} && x\in\rit^n\;,\\
\rresp{\ellip_R(n,[p+h],m_{p+h})& =\lambda ^\half(n,p+h,m_{p+h})\ e^{-2\pi i(p+h)x}}\end{align*}
{\bbf is a discrete valuation semiring of which\/}:
\Bean
\item the uniformizing element is 
$r(p+h,m_{p+h},\o p^\ell)^p=
\prod^n_{c=1}r_c(p+h,m_{p+h},\o p^\ell)^p$;

\item the units are the invertible elements 
$e^{2\pi i(p^\ell)x'}=\prod^n_{c=1}\ e^{2\pi i(p^\ell)x'_c}$
\resp{$e^{-2\pi i(p^\ell)x'_c}$}, $x'\in\rit^n$.
\Ee
}
\vskip 11pt

\begin{proof} This proposition is the $n$-dimensional generalization of proposition 3.2.5, taking into account corollary 3.2.6.
\end{proof}
\vskip 11pt

\subsubsection{Global elliptic semimodules covered by $p^\ell$- roots}

Let $\ELLIP_L(n,i\ge p,m_i)$
\resp{$\ELLIP_R(n,i\ge p,m_i)$} be the 
$\Gamma (\Phi (G^{(n)}(L^T_{[v_p]})))$-semi\-module
\resp{$\Gamma (\Phi (G^{(n)}(L^T_{[\o v_p]})))$-semimodule} with terms $i\ge p$.

The $n$-dimensional global elliptic semimodule
\begin{align*}
\ELLIP_L(n,i\ge p,m_i,\o p^{(\ell)}) &= \bigoplus^w_{h=0}\oplus_{m_{p+h}} \ellip_L ( n,p+h,m_{p+h},\o p^{(\ell ) } ) \\ 
\noalign{\mbox{}\hfill $0\le h\le w\le\infty\;,$}
\rresp{\ELLIP_R(n,i\ge p,m_i,\o p^{(\ell)}) &= \bigoplus^w_{h=0}\oplus_{m_{p+h}} \ellip_R(n,p+h,m_{p+h},\o p^{(\ell)})}
\end{align*}
covers 
$\ELLIP_L(n,i\ge p,m_i)$
\resp{$\ELLIP_R(n,i\ge p,m_i)$} in the sense that each term
\begin{align*}
\ellip_L(n,p+h,m_{p+h},\o p^\ell) & \simeq \prod^n_{c=1} r_c(p+h,m_{p+h},\o p^\ell)^{p^\ell}\ e^{2\pi i(p^\ell)x'_c}\\
\rresp{\ellip_R(n,p+h,m_{p+h},\o p^\ell) & \simeq \prod^n_{c=1} r_c(p+h,m_{p+h},\o p^\ell)^{p^\ell}\ e^{-2\pi i(p^\ell)x'_c}}
\end{align*}
of $\ELLIP_L(n,i\ge p,m_i,\o p^{(\ell)})$
\resp{$\ELLIP_R(n,i\ge p,m_i,\o p^{(\ell)})$}
is the covering by $p^{(\ell)}$ roots ($\ell$ varying from one term to another) of each term
$\ellip_L(n,[p+h],m_{p+h})$
\resp{$\ellip_R(n,[p+h],m_{p+h})$}
of $\ELLIP_L(n,i\ge p,m_i)$
\resp{$\ELLIP_R(n,i\ge p,m_i)$}.
\vskip 11pt

A semisheaf $\Phi (G^{(n)}_{\o p^\ell}(L^T_{[v_p]}))$
\resp{$\Phi (G^{(n)}_{\o p^\ell}(L^T_{[\o v_p]}))$}, of which sections are the functions\linebreak
$\ellip_L(n,p+h,m_{p+h},\o p^{(\ell)})$
\resp{$\ellip_R(n,p+h,m_{p+h},\o p^{(\ell)})$}, is associated with the $n$-dimensional global elliptic semimodule
$\ELLIP_L(n,i\ge p,m_i,\o p^{(\ell)})$
\resp{$\ELLIP_R(n,i\ge p,m_i,\o p^{(\ell)})$} as introduced precedently.
\vskip 11pt

\subsection{Bilinear local cuspidal representations associated with their global correspondents covered by $p^\ell$-th roots}

\subsubsection{Non-archimedean local fields}

Let $K^+_p/L^+_p$
\resp{$K^-_p/L^-_p$} denote a finite Galois extension of a non-archimedean $p$-adic \lr semifield $L^+_p$ \resp{$L^-_p$} which is finite extension of
$\QQ^+_p=\ZZ_p\otimes\QQ_+$
\resp{$\QQ^-_p=\ZZ_p\otimes\QQ_-$}
where $\ZZ_p=\ds\lim_{\overleftarrow r}\ZZ/(p^r)$
and where $\QQ^+_p$
\resp{$\QQ^-_p$} is the completion of $\QQ_+$ (the positive rational numbers)
\resp{$\QQ_-$ (the negative rational numbers)} in the $p$-adic metric.
\vskip 11pt

Let $[K^+_p:L^+_p]=q$ denote the degree of this extension, $v_p$ a discrete valuation of $L^+_p$ with semiring $A^+_p$ and $\omega _r$ the different prolongations of $v_p$ to $K^+_p$ \cite{Ser3}.

Let $B^+_p$ be the integral closure of $A^+_p$ into $K^+_p$ in such a way that the $A^+_p$-semimodule $B^+_p$ be finitely generated and that its field of fractions be $K^+_p$.
\vskip 11pt

Let $\beta ^+_{p_1}\subset\beta ^+_{p_2}\subset\cdots
\subset\beta ^+_{p_r}$ be a chain of distinct prime ideals of $B^+_p$ and let $m(A^+_p)=\beta ^+_{p_r}\cap A^+_p$ define the division of $\beta ^+_{p_r}$ by the maximal ideal $m(A^+_{p_r})$, noted $\beta ^+_{p_r}\mid m(A^+_p)$ ($\beta ^+_{p_r}$ contains the ideal $m(A^+_p)B^+_p$ generated by $m(A^+_p)$).  If the set $\{\beta ^+_{p_r}\}_r$ of prime ideals of $B^+_p$ verifies $\{\beta ^+_{p_r}\}_r\cap A^+_p=\Os_{K^+_p}\cap A^+_p$, then $\{\beta ^+_{p_r}\}_r=\Os_{K^+_p}$ where $\Os_{K^+_p}$ is the semiring of integers of $K^+_p$.

For each prime ideal $\beta ^+_{p_r}$ above $m(A^+_p)$, $B^+_p/\beta ^+_{p_r}$ is an extension of $A^+_p/m(A^+_p)$ of which extension degree is the ``local'' residue degree
$f_{\beta ^+_{p_r}}=[B^+_p/\beta ^+_{p_r}:A^+_p/m(A^+_p)]$ of $\beta ^+_{p_r}$ in the extension $K^+_p/L^+_p$.

The exponent $e_{\beta ^+_{p_r}}$ of $\beta ^+_{p_r}$ in the decomposition of $m(A^+_p)B^+_p$ into prime ideal is the ramification index of $\beta ^+_{p_r}$ in the extension $K^+_p/L^+_p$.

The semiring $B^+_p/m(A^+_p)B^+_p$ is an $A^+_p/m(A^+_p)$-semialgebra of degree  $q=\sum\limits_{\beta _{p_r}\mid p_p} f_{\beta ^+_{p_r}}\ e_{\beta ^+_{p_r}}=[K^+_p:L^+_p]$ and is isomorphic to the product $\prod\limits_{\beta^+ _{p_r}\mid p^+_p} B^+_p/\beta ^{+^{e_{\beta ^+_{p_r}}}}_{p_r}$.

Remark that the ramification indices $e_{\beta ^+_{p_r}}$, referring to
$m(A^+_p)$ are all equal to $e_{\beta ^+_{p}}$ in Galois extensions \cite{Ser3}.
\vskip 11pt

To each prime ideal $\beta ^+_{p_r}$ corresponds a residue semifield
$k_{K_{\beta ^+_{p_r}}}=\Os_{K^+_p\mid \beta ^+_{p_r}}/m(A^+_p)=\beta ^+_{p_r}/m(A^+_p)$ where $\Os_{K^+_p\mid \beta ^+_{p_r}}$ is the semiring of integers of $K^+_p$ restricted to $\beta ^+_{p_r}$.

According to J.P. Serre, the semiring $B^+_p$ is a discrete valuation semiring \cite{Ser3}.  Taking into account that the different valuations $w_r$, corresponding to the prime ideal $\beta ^+_{p_r}$, define each one a norm on $K^+_p$ making $K^+_p$ a Hausdorff topological vector semispace over $L^+_p$.  As $L^+_p$ is assumed to be complete and as the topology $\Fs_r$ defined by $w_r$ is a product topology on $K^+_p$ not depending on the index $r$, there is only one $w_r$ which is relevant.
\vskip 11pt

Let $\tilde\omega _{K^+_p}$ denote the uniformizer (i.e. a prime element) in $\Os_{K^+_p}$.  $\Os_{K^+_p}$ is then the inverse limit of $\Os_{K^+_p}/(\tilde \omega _{K^+_p})^r$.

The number of elements in $K^+_p$ is thus $p^q=p^{\sum\limits_r f_r\cdot e_r}\equiv p^{\sum\limits_rq_r}$ \cite{Kna}, where $f_{\bp}$ and $e_{\bp}$ have been written in condensed form respectively as $f_r$ and $e_r$.
\vskip 11pt

Remark that the right case can be handled similarly with the evident following notations: semifields $L^-_p$ and $K^-_p$, semiring $A^-_p$, $A^-_p$-semimodule $B^-_p$, prime ideals $\beta ^-_{p_r}$, semiring of integers $\Os_{K^-_p}$, residue semifield $k_{K_{\beta ^-_{p_r}}}$ and so on.
\vskip 11pt

\subsubsection{Proposition (Global $\leftrightarrow $ local correspondences between extension (semi)fields)}

{\em {\bbf The set of intermediate subsemifields\/}
$\{\wt L_{v_{p+h,m_{p+h}}}\}^\infty _{h=0}$
\resp{$\{\wt L_{\o v_{p+h,m_{p+h}}}\}^\infty _{h=0}$}
of $\wt L_L$
\resp{$\wt L_R$}, extensions of a numberfield $k$ of char $0$ as introduced in section 2.1, or equivalently the set of corresponding archimedean completions
$\{L_{v_{p+h,m_{p+h}}}\}^\infty _{h=0}$
\resp{$\{L_{\o v_{p+h,m_{p+h}}}\}^\infty _{h=0}$},
{\bbf can be covered in an etale way by a (set of) $p$-adic finite extension semifield(s) leading to a global $\leftrightarrow $ local isomorphism if\/}

\[ p^{\sum\limits_rf_r\cdot e_r}=\#\Nu\times N\times \sum_h f_{v_{p+h,m_{p+h}}}\;,\]
i.e. if {\bbf their numbers of elements correspond\/},
where:
\Bi
\item $\#\Nu$ is the number of global nonunits;
\item $N\in\NN$ is the Galois extension degree of a quantum, i.e. an irreducible subsemifield $\wt L_{v_i^1}$ or $\wt L_{\o v_i^1}$;
\item $f_{v_{p+h,m_{p+h}}}=p+h\in\NN$ is the global residue degree of $\wt L_{v_{p+h,m_{p+h}}}$;
\item $p^{\sum\limits_rf_r\cdot e_r}$ is the number of elements in the finite $p$-adic Galois extension(s) as developed in section 3.3.1.
\Ei
}
\vskip 11pt

\begin{proof}
\Be
\item According to section 2.1, the Galois extension degree of the pseudo-ramified completion $L_{v_{p+h,m_{p+h}}}$ (or
$L_{\o v_{p+h,m_{p+h}}}$) at the infinite real place 
$v_{p+h}$
(or $\o v_{p+h}$) is
\[ [L_{v_{p+h,m_{p+h}}}:k]\equiv [L_{\o v_{p+h,m_{p+h}}}:k]=(p+h)\ N\]
(in the residue class zero of the integers modulo $N$).

Consequently, the number of elements in
$L_{v_{p+h,m_{p+h}}}$ or in the subsemifield
$\wt L_{v_{p+h,m_{p+h}}}$ is equal to:
\[ \left|\wt L_{v_{p+h,m_{p+h}}}\R|=\#\Nu\times N\times f_{v_{p+h}}\;.\]

\item Referring to section 3.3.1, if the number of elements
$\left|\L\{\wt L_{v_{p+h,m_{p+h}}}\R\}_h\R|=\#\Nu\times N\times \sum\limits_hf_{v_{p+h,m_{p+h}}}$ of the set of subsemifields
$\wt L_{v_{p+h,m_{p+h}}}$ is a power of $p$, i.e. if
$\left|\L\{\wt L_{v_{p+h,m_{p+h}}}\R\}\R|=p^q$, where $q=\sum\limits_rf_r\cdot e_r$, then
$\left\{\wt L_{v_{p+h,m_{p+h}}}\R\}_h$ is isomorphic to the $p$-adic extension semifield $K^+_p$ of dimension $q$.

Similarly, if $\left|\L\{\wt L_{v_{p+h,m_{p+h}}}\R\}_h\R|<p^q$,
$\left\{\wt L_{v_{p+h,m_{p+h}}}\R\}_h$ is monomorphic to $K^+_p$, and if\linebreak
$\left|\L\{\wt L_{v_{p+h,m_{p+h}}}\R\}_h\R|>p^q$,
$\left\{\wt L_{v_{p+h,m_{p+h}}}\R\}_h$ is epimorphic to $K^+_p$.
\qedhere
\Ee
\end{proof}
\vskip 11pt

\subsubsection{Corollary}

{\em The set of completions
$\left\{\wt L_{v_{p+h,m_{p+h}}}\R\}_h$
\resp{$\left\{\wt L_{\o v_{p+h,m_{p+h}}}\R\}_h$}
is a $p$-adic semifield if:
\Bean
\item its number of elements is a power of $p$, i.e. if
\[ \left| \L\{   L_{v_{p+h,m_{p+h}}} \R\}_h\R|\equiv
 \left| \L\{   L_{\o v_{p+h,m_{p+h}}} \R\}_h\R|=p^q\;;\]

\item they are defined as completion(s) of $K^+_p$ \resp{$K^-_p$} for the topology defined by its $p$-adic absolute value.
\Ee}
\vskip 11pt

\begin{proof}
\Bean
\item According to proposition 3.3.2, if the number of elements of
$\left\{  L_{v_{p+h,m_{p+h}}}\R\}_h$ is a power of $q$, this set of infinite completions is isomorphic to a finite $p$-adic extension semifield.

\item As each infinite completion
$ L_{v_{p+h,m_{p+h}}}$ results from an isomorphism of compactification into a closed compact subset of
$\rit_+$, $\left\{ L_{v_{p+h,m_{p+h}}}\R\}$ will define a completion of a subset of $K^+_p$ if the considered topology refers to an ultrametric $p$-adic absolute value.

Similarly, $\left\{ L_{\o v_{p+h,m_{p+h}}}\R\}_h$ can generate a completion of a subset of $K^-_p$.\qedhere
\Ee
\end{proof}
\vskip 11pt

\subsubsection{Proposition (Local elliptic semimodule)}

{ \em Let $\Phi (G^{(2)}(L^T_{[v_p]}))$
\resp{$\Phi (G^{(2)}(L^T_{[\o v_p]}))$}
be the two-dimensional global semisheaf to which is associated the global elliptic semimodule
$\ELLIP_L(2,i\ge p,m_i)$
\resp{$\ELLIP_R(2,i\ge p,m_i)$} and let
$\Phi (G^{(2)}_{\o p^{(\ell)}}(L^T_{[v_p]}))$
\resp{$\Phi (G^{(2)}_{\o p^{(\ell)}}(L^T_{[\o v_p]}))$}
be its etale covering global semisheaf by $p^{(\ell)}$ roots to which is associated the global elliptic semimodule
$\ELLIP_L(2,i\ge p,m_i,\o p^{(\ell)})$
\resp{$\ELLIP_R(2,i\ge p,m_i,\o p^{(\ell)})$}.

Then, $\Phi (G^{(2)}_{\o p^{(\ell)}}(L^T_{[v_p]}))$
\resp{$\Phi (G^{(2)}_{\o p^{(\ell)}}(L^T_{[\o v_p]}))$} is covered by a $p$-adic local semigroup
$G^{(2)}(K^+_p)\equiv T_2(K^+_p)$
\resp{$G^{(2)}(K^-_p)\equiv T^t_2(K^-_p)$}
on a semischeme $S^+$ \resp{$S^-$} {\rm \cite{Mes}} which is a flat $\Os_{K^+_p}$-semimodule
\resp{$\Os_{K^-_p}$-semimodule} if:
\Bean
\item $p^m=\#\Nu\times N$, i.e. if the number of elements $\#\Nu\times N$ in a global quantum is a power of $p$;

\item $p^{2q}=\L(\sum\limits_\ell((p^\ell-1)/2)\times m_{(p^\ell-1)/2}\times p^m\R)^2$ (case $n=2$, two-dimensional case), where $m_{(p^\ell-1)/2}$ denotes the multiplicity of the covering sections by $p^\ell$ roots, or 
\[p^{2q}=\L(\sum\limits_h(p+h)\times m_{p+h}\times p^m\R)^2\;;\]

\item there are Frobenius endomorphisms with generator 
$x\to x^{p^{f_r}}$
\resp{$-x\to -x^{p^{f_r}}$} resulting form the cyclicity of the Galois semigroup
$\Gal(k_{K^+_{\bp}}/k_{L^+_p})$
\resp{$\Gal(k_{K^-_{\bmp}}/k_{L^-_p})$}
where $k_{K^+_{\bp}}$
\resp{$k_{K^-_{\bmp}}$} and $k_{L^+_p}$ \resp{$k_{L^-_p}$}
are respectively residue semifields of 
$K^+_p$
\resp{$K^-_p$} and of $L^+_p$
\resp{$L^-_p$}.

On the other hand, the image of a member of
$\Gal(K^+_p/L^+_p)$
\resp{$\Gal(K^-_p/L^-_p)$} is of the form 
$x\to x^{p^{f_r\cdot e_r}}=\mu ^{q_r}$
\resp{$-x\to -x^{p^{f_r\cdot e_r}}=-\mu ^{q_r}$};

\item there are embeddings
\begin{align*}
e_{L_{[v_p]}\to K^+_p}: \qquad \lambda (2,p+h,m_{p+h})&\To\lambda _p(2,r,m_r)\\
\rresp{e_{L_{[\o v_p]}\to K^-_p}: \qquad \lambda (2,p+h,m_{p+h})&\To\lambda _p(2,r,m_r)}
\end{align*}
of the product, right by left, of Hecke characters
$\lambda (2,p+h,m_{p+h})$ over $L_{[v_p]}$
\resp{$L_{[v_p]}$} into their equivalents
$\lambda _p(2,r,m_r)$ over $K^+_p$
\resp{$K^-_p$},
\Ee
in such a way that the {\bbf global elliptic semimodule
$\ELLIP_L(2,p\le i,m_i)$
\resp{\linebreak $\ELLIP_R(2,p\le i,m_i)$} and its covering by $p^{(\ell)}$ roots
$\ELLIP_L(2,p\le i,m_i,\o p^{(\ell)})$
\resp{$\ELLIP_R(2,p\le i,m_i,\o p^{(\ell)})$} be covered by the local elliptic
$\End(G^{(2)}(K^+_p))$-semimodule
\resp{$\End(G^{(2)}(K^-_p))$-semimodule}\/} referring to Drinfeld {\rm \cite{Drin}} and Anderson {\rm \cite{And}}:
\begin{align*}
\ELLIP(2,x,K^+_p) &= \bigoplus_r\L(\lambda ^\half_p(2,r,m_r)(x))\ f(\mu ^{q_r\cdot r} \R)\\
\rresp{\ELLIP(2,-x,K^-_p) &= \bigoplus_r\L(\lambda ^\half_p(2,r,m_r)(-x))\ f(-\mu ^{q_r\cdot r} \R)}
\end{align*}
for every closed point $x$ \resp{$-x$} of $K^+_p$ \resp{$K^-_p$} where
$f(\mu ^{q_r\cdot r} )$
\resp{$f(-\mu ^{q_r\cdot r} )$} is a function of the Frobenius endomorphism.
}
\vskip 11pt

\begin{proof}
It is thus asserted that there exists an isomorphism
\begin{align*}
i_{G^{(2)}_{\o p^{(\ell)}}\to \Os(G^{(2)}(K^+_p))} : \qquad \Phi (G^{(2)}_{\o p^{(\ell)}}(L^T_{[v_p]})&\To G^{(2)}(K^+_p)\\
\rresp{i_{G^{(2)}_{\o p^{(\ell)}}\to \Os(G^{(2)}(K^-_p))} : \qquad \Phi (G^{(2)}_{\o p^{(\ell)}}(L^T_{[\o v_p]})&\To G^{(2)}(K^-_p)}
\end{align*}
from the global semisheaf
$\Phi (G^{(2)}_{\o p^{(\ell)}}(L^T_{[v_p]}))$
\resp{$\Phi (G^{(2)}_{\o p^{(\ell)}}(L^T_{[\o v_p]}))$}
to the local semigroup
$ G^{(2)}(K^+_p)$\linebreak 
\resp{$ G^{(2)}(K^-_p)$}
in such a way that the diagram
\[ \begin{CD}
\Phi (G^{(2)}_{\o p^{(\ell)}}(L^T_{[v_p]}) @>\sim>>
 G^{(2)}(K^+_p)\\
@VVV @VVV\\
\ELLIP_L(2,i\ge p,m_i,\o p^{(\ell)}) @>\sim>>
\ELLIP(2,x,K^+_p) \end{CD}\]
be commutative.
\vskip 11pt

This can be achieved if the number of points of the local semigroup is equal to the number of points of the global semisheaf.

The number of points of the local semigroup 
$G^{(2)}(K^+_p)$
\resp{$G^{(2)}(K^-_p)$} is $p^{2q}$ according to proposition 3.3.2, the factor ``$2$'' resulting from the dimension $n=2$.

The number of points of the global semisheaf
$\Phi (G^{(2)}_{\o p^{(\ell)}}(L^T_{[v_p]})$, covered by $G^{(2)}(K^+_p)$, is\linebreak  $\L(\sum\limits_\ell ((p^\ell-1)/2)\times m_{p^{\ell/2}}\times p^m\R)^2$ since its sections of type 
$T^1_{c_L}[p+h,m_{p+h},\o p^\ell]$, according to proposition 3.2.1, have a number of points
\[ n_{T^1[p+h,\o p^\ell]}=\#\Nu\times N\times (p^\ell-1)/2\;.\]
From proposition 3.2.5, it results that the covering map
$R^2_{L_{p^\ell\to p+h}}$ is an isomorphism if
$((p^\ell-1)/2)^2=(p+h)^2$, which explains that
\[ p^{2q}=\L(\sum_h(p+h)\times m_{(p+h)}\times p^m\R)^2\;.\]
By this way, each global point of the algebraic semigroup
$G^{(2)}_{\o p^{(\ell)}}(L^T_{[v_p]})$ is in one-to-one correspondence with a local closed point of the local group
$G^{(2)}(K^+_p)$: this results from the conditions a) and b) of this proposition.
\vskip 11pt

In order that the local $p$-adic elliptic
$\End(G^{(2)}(K^+_p))$-semimodule $\ELLIP(2,x,K^+_p)$
\resp{\linebreak $\End(G^{(2)}(K^-_p))$-semimodule $\ELLIP(2,-x,K^-_p)$}
corresponds to a local cuspidal \lr form, the conditions c) and d) must be fulfilled in analogy with the global case considered in proposition 3.1.2, i.e.
\Be
\item a Frobenius substitution:
\[ \mu \To \mu ^{q_r} \qquad
\rresp{-\mu \To -\mu ^{q_r}}\]
on every local Frobenius endomorphism
\[ \mu : \quad x\To x^p \qquad
\rresp{-\mu : \quad -x\To -x^p}\;;\]

\item an embedding 
\begin{align*}
e_{L_v\to K^+_p}: \qquad \lambda (2,p+h,m_{p+h})&\To \lambda _p(2,r,m_r)\\
\rresp{e_{L_{\o v}\to K^-_p}: \qquad \lambda (2,p+h,m_{p+h})&\To \lambda _p(2,r,m_r)}\end{align*}
of the Hecke character 
$\lambda (2,p+h,m_{p+h})$ into $\lambda _p(2,r,m_r)$ which is the square of the coefficient for every local point $x$ \resp{$-x$}
of the local elliptic
$\End(G^{(2)}(K^+_p))$-semimodule $\ELLIP(2,x,K^+_p)$
\resp{$\End(G^{(2)}(K^-_p))$-semimodule $\ELLIP(2,-x,K^-_p)$}.
\Ee

This condition corresponds to the embedding
$i(a_\ell)=\tr(\Frob_\ell)$ into $\QQ_p$ of the ring of the integers of a finite extension $E_f$ (i.e. the ring of the coefficients of the cuspidal form $f$) of $\QQ$, $a_\ell$ being the coefficient of the cuspidal form, as introduced by Deligne \cite{Del} and Serre \cite{Ser4} and mentioned in \cite{Win}.
\end{proof}
\vskip 11pt

\subsubsection{Corollary}

{\em There is an epimorphism
\begin{align*}
e_{\Phi  ( G^{(2)}(K^+_p) )\to\Phi (G^{(2)}_{\o p^{(\ell )}}}:\qquad G^{(2)}(K^+_p)&\To\Phi  ( G^{(2)}_{\o p^{(\ell )}}(L^T_{[v_p]}))\\
\rresp{e_{\Phi (G^{(2)}(K^-_p))\to\Phi (G^{(2)}_{\o p^{(\ell )}}}:\qquad G^{(2)}(K^-_p)&\To\Phi  ( G^{(2)}_{\o p^{(\ell )}}(L^T_{[\o v_p]}))}
\end{align*}
from the local semigroup 
$G^{(2)}(K^+_p)$
\resp{$G^{(2)}(K^-_p)$} into the global semisheaf
$\Phi  ( G^{(2)}_{\o p^{(\ell )}}(L^T_{[v_p]}))$
\resp{$\Phi  ( G^{(2)}_{\o p^{(\ell )}}(L^T_{[\o v_p]}))$}
such that 
$\ELLIP(2,x,K^+_p)$
\resp{$\ELLIP(2,-x,K^-_p)$} be projected onto
$\ELLIP_L(2,i\ge p,m_i,\o p^{(\ell)})$
\resp{$\ELLIP_R(2,i\ge p,m_i,\o p^{(\ell)})$}
if the number of ``local'' points is superior to the number of ``global points'', i.e. if
\[ p^{2q}>\L(\sum_\ell (p^\ell/2)\times m_{p^\ell/2}\times p^m\R)^2\;.\]
}
\vskip 11pt

\begin{proof}
This directly results from proposition 3.3.4 and,  more particularly, from condition b), the other conditions a), b) and d) being unchanged.
\end{proof}
\vskip 11pt

\subsubsection{Proposition}

{\em {\bbf The local $p$-adic elliptic 
$\End(G^{(2)}(K^+_p))$-semimodule $\ELLIP(2,x,K^+_p)$
\resp{\linebreak $\End ( G^{(2)}(K^-_p))$-semimodule 
$\ELLIP ( 2,-x,K^-_p)$}, 
corresponding to a local $p$-adic \lr cuspidal form\/}, results from the \lr semisheaf 
$\Phi  ( G^{(2)}(L^T_v))$
\resp{$\Phi  ( G^{(2)}(L^T_{\o v}))$} on the global algebraic semigroup
$G^{(2)} ( L^T_{v} )$
\resp{$G^{(2)}(L^T_{\o v}))$} by the commutative diagram (the left case being only considered here):
\begin{footnotesize}
\[\begin{CD}
{\Phi (G^{(2)}(L^T_v)) \hspace{-2mm}} 
@>>{\Em^{(2)}_{v\to[v_p]}}>
{\hspace{-2mm} \Phi  ( G^{(2)} ( L^T_{[v_p]} ) ) \hspace{-2mm}}
 @>{\sim}>>
{\hspace{-2mm} \Phi (G^{(2)} _{\o p^{(\ell)}} (L^T_{[v_p]})) \hspace{-2mm}}
 @>>{i_{G^{(2)}_{\o p^{(\ell)}}\to\Phi ( G^{(2)}(K^+_p)}}>
{\hspace{-2mm} G^{(2)}(K^+_p)} \\
@VVV @VVV \mbox{} @VVV@VVV \\
{\hspace{-2mm}\ELLIP_L(2,i,m_i) \hspace{-2mm}}
 @>>>
{\hspace{-2mm}\ELLIP_L(2,p\le i,m_i)}
 @>>>
{\ELLIP_L(2,p\le i,m_i,\o p^{(\ell)})\hspace{-5mm}} @>>>
{\hspace{-5mm}\ELLIP_L(2,x,K^+_p)\hspace{-5mm}}\end{CD}\]
\end{footnotesize}
}
\vskip 11pt

\begin{proof}
This is a consequence of corollary 3.2.6, sections 3.2.7 and 3.2.12 as well as proposition 3.3.4.

That is to say that:
\Bean
\item the epimorphism $\Em^{(2)}_{v\to[v_p]}$ sends the semisheaf
$\Phi (G^{(2)}(L^T_v))$ on the algebraic semigroup $ G^{(2)}(L^T_v)$
over ``$t$'' sets of toroidal archimedean completions, $1\le i\le t\le\infty $, into the semisheaf $\Phi  ( G^{(2)} ( L^T_{[v_p]} ) )$ on the algebraic semigroup $G^{(2)} ( L^T_{[v_p]} ) $ restricted to toroidal completions above $v_p$.

\item this allows to find the covering semisheaf
$\Phi (G^{(2)} _{\o p^{(\ell)}} (L^T_{[v_p]}))$ by $p^{(\ell)}$ roots of the semisheaf
$\Phi (G^{(2)}  (L^T_{[v_p]}))$ and the $p$-adic local semigroup
$ G^{(2)}  (K^+_p)$.
\Ee
The local $p$-adic elliptic
$\End(G^{(2)}(K^+_p))$-semimodule $\ELLIP(2,x,K^+_p)$ then results from the global elliptic
$\Gamma (\Phi (G^{(2)}(L^T_v)))$-semimodule $\ELLIP_L(2,i,m_i)$.
\end{proof}
\vskip 11pt

\subsubsection{Proposition}

{\em {\bbf The Serre (Eichler, Deligne, Shimura) conjecture\/}, asserting that Galois representations $\rho :G_\QQ\to \GL_2(\o\QQ_p)$ can be associated to modular forms, directly
{\bbf results from proposition 3.3.6, and more particularly, from the epimorphism\/}:
\begin{align*}
\Em^{(2)}_{L^T_v\to K^+_p}: 
\qquad \GL_2(L^T_v)&\To \GL_2(K^+_p)\;, && \text{with\ } T^2(L^T_v)\equiv \GL_2(L^T_v)\;,\\
\rresp{\Em^{(2)}_{L^T_{\o v}\to K^-_p}: \qquad \GL_2(L^T_{\o v})&\To \GL_2(K^-_p)}\;,
\end{align*}
sending the algebraic semigroup
$\GL_2(L^T_v)$
\resp{$\GL_2(L^T_{\o v})$} over the set $L^T_v$
\resp{$L^T_{\o v}$} of archimedean completions into the algebraic semigroup 
$\GL_2((K^+_p)$
\resp{$\GL_2((K^-_p)$} at the following conditions:
\Be
\item $L^T_v$
\resp{$L^T_{\o v}$} is extended to $\o\QQ_+$
\resp{$\o\QQ_-$},
$K^+_p$ to $\QQ^+_p$ and 
$K^-_p$ to $\QQ^-_p$;

\item there is an epimorphsm
\begin{align*}
\Em^{(2)}_{L^T_v\to L^T_{[v_p]}}: \qquad \GL_2(L^T_v) &\To \GL_2(L^T_{[v_p]})\\
\rresp{\Em^{(2)}_{L^T_{\o v}\to L^T_{[\o v_p]}}: \qquad \GL_2(L^T_{\o v}) &\To \GL_2(L^T_{[\o v_p]})}
\end{align*}
from $\GL_2(L^T_{v})$
\resp{$\GL_2(L^T_{\o v})$} into the algebraic semigroup
$\GL_2(L^T_{[v_p]}$
\resp{$\GL_2(L^T_{[\o v_p]}$} over the set of completions
$L^T_{[v_p]}$
\resp{$L^T_{[\o v_p]}$} superior or equal to $v_p$ \resp{$\o v_p$};

\item $p^m=\#\Nu\times N$;

\item $p^{2q}=\sum\limits_h (p+h)\times m_{(p+h)}\times p^m)^2$;

\item there are Frobenius endomorphisms
$\mu :x\to x^{p^{f_r\cdot e_r}}=\mu ^{q_r}$
\resp{$-\mu :-x\to -x^{p^{f_r\cdot e_r}}=-\mu ^{q_r}$};

\item there are embeddings
$e_{L_v\to K^+_p}:\lambda (2,i,m_i)\to\lambda _p(2,r,m_r)$
\resp{$e_{L_{\o v}\to K^-_p}:\lambda (2,i,m_i)\to\lambda _p(2,r,m_r)$} of the product, right by left, of Hecke characters
$\lambda (2,i,m_i)$ over $L_v$\resp{$L_{\o v}$}
into their equivalents
$\lambda_p (2,r,m_r)$ over $K^+_p$\resp{$K^-_p$}.
\Ee
}
\vskip 11pt

\begin{proof}
(for the left case, the right case being handled similarly)
\Be
\item First, we have to consider the mapping of 
$G_\QQ=\Gal(\o\QQ/\QQ)$, or, more restrictively, of
$\Gal(\wt L_L/k)$ or of $W_{\wt L_v}$, defined in section 2.4, into the set (or the sum) of the equivalence classes of the representations space
$\Irr\Rep^{(2)}(W_{\wt L_v})$ of the global Weil group $W_{\wt L_v}$ in such a way that:
\[ \Irr\Rep^{(2)}(W_{\wt L_v})=G^{(2)}(L_v)\]
as developed in proposition 3.1.5.

Then, the algebraic semigroup  (of matrices)
$\GL_2(L^T_v)$, isomorphic to the algebraic semigroup
$G^{(2)}(L_v)$ according to the preceeding developments of this paper, is sent by the epimorphism
$\Em^{(2)}_{L^T_v\to L^T_{[v_p]}}$ into the algebraic semigroup
$\GL_2(L^T_{[v_p]})$ with coefficients on completions superior or equal to $v_p$.

Finally, $\GL_2(L^T_{[v_p]})$ is sent by the isomorphism
$\Im^{(2)}_{L^T_{[v_p]}\to K^+_p}$ into $GL_2(K^+_p)$ at the conditions a), b), c) and d) of proposition 3.3.4.

It then results that the Galois representation
\[\rho :\qquad G_\QQ\To \GL_2(\o\QQ_p)\]
is in fact an epimorphism corresponding mainly to the compositions of morphisms:
\[\Im^{(2)}_{L^T_{[v_p]}\to K^+_p}\circ
\Em^{(2)}_{L^T_{v}\to L^T_{[v_p]}}\circ
\Irr\Rep^{(2)}(W^{ab}_{\wt L_v)})\;;\]

\item The fact that the Galois representation
$\rho : G_\QQ\to \GL_2(\o\QQ_p)$ can be associated to a modular form results from the commutative diagram of proposition 3.3.6, taking into account that the global elliptic semimodule
$\ELLIP_L(2,i,m_i)$ constitutes a cuspidal representation of $G^{(2)}(L_v)$ according to proposition 3.1.5 and is in one-to-one correspondence with a cuspidal form as developed in \cite{Pie2}.\qedhere
\Ee
\end{proof}
\vskip 11pt

\subsubsection{Corollary}

{\em {\bbf The Shimura-Taniyama-Weil conjecture\/}, associated with the action of $G_\QQ$ on the elliptic curve $E[p]$ leading to a continuous representation
\[ \rho _{E,p}:\qquad G_\QQ\To \GL_2(\FF_p)\]
in such a way that:
\[\tr (\rho _{E,p}(\Frob_p))=p+1-\#E(\FF_p) \pmod p\;,\]
also directly results from propositions 3.3.4 and 3.3.6.
}
\vskip 11pt

\begin{proof}
In the new context proposed here, the Shimura-Taniyama-Weil conjecture is a special case of the Serre conjecture since, referring to proposition 3.3.6 and to the extension of the residue field $k_{K_p}$ to $\FF_p$, it results from a global elliptic semimodule 
$\ELLIP_L(2,p\le i,m_i)$
\resp{$\ELLIP_R(2,p\le i,m_i)$} restricted to the $i$-th$=p$-th terms, i.e. to the case where $h=0$ and $\ell=1$, according to corrolaries 3.2.2 and 3.2.6.

Remark that the Shimura-Taniyama-Weil-conjecture was specifically studied in \cite{Pie2} from this new point of view.
\end{proof}
\section{Deformations of Galois representations}

Two kinds of deformations of $n$-dimensional representations of Galois or Weil groups given by bilinear algebraic semigroups over complete global and local noetherian bisemirings, in reference with the work of B. Mazur \cite{Maz1}, are envisaged:
{\bbf
\Be
\item global and local bilinear deformations inducing the invariance of their respective global and local bilinear residue (semi)fields.

\item global and local bilinear ``quantum'' deformations leaving invariant the orders of the inertia subgroups.
\Ee
}

\subsection{Local and global coefficient semiring homomorphisms}

\subsubsection{Local coefficient semiring homomorphisms (left case)}


A coefficient semiring $\BB^+_p$ is, according to B. Mazur \cite{Maz2}, a complete noetherian local semiring with finite residue semifield $k_{K^+_p}$.  It is characterized by a profinite topology given by a base of prime ideals
$\wt\omega _{K^+_p}\BB^+_p$ in such a way that:
\[ \BB^+_p=\lim_{r\to\infty }\BB^+_p/\wt\omega ^r_{K^+_p}\BB^+_p\]
where $\wt\omega _{K^+_p}\BB^+_p$ is the maximal ideal of $\BB^+_p$.

The discrete valuation semiring $B^+_p$, introduced in section 3.3.1 as the integral closure of $A^+_p$ in the finite Galois extension $K^+_p/L^+_p$ of the $p$-adic semifield $L^+_p$, is a noetherian local semiring if the chain
$\beta ^+_{p_1}\subset 
\beta ^+_{p_2}\subset \dots \subset
\beta ^+_{p_r}$ of prime ideals of $B^+_p$ tends to $\infty $, i.e. if $r\to\infty $.

It will then be assumed in this chapter that $B^+_p$ is a noetherian local semiring on the completion of $K^+_p$.
\vskip 11pt

Let $B^{'+}_p$ be another coefficient semiring being the integral closure of $A^+_p$ in another finite Galois extension 
$B^{'+}_p/L^+_p$ and let $k_{K^{'+}_p}=\{\Os_{K^{'+}_p\mid \beta ^{'+}_{p_r}}\}_r$ be its residue semifield defined on the set of residue subsemifields
$\Os_{K^{'+}_p\mid \beta ^{'+}_{p_r}}/m(A^+_p)$ restricted to the prime ideals
$\beta ^{'+}_{p_r}$ and
$B ^{'+}_{p}$, referring to section 3.3.1.

The semiring $B ^{'+}_{p}/m(A^+_p)B ^{'+}_{p}$ is also an
$A^+_p/m(A^+_p)$-semialgebra of degree
$q'=\sum\limits_{\beta ^{'+}_{p_r}\mid p_p}
f_{\beta ^{'+}_{p_r}}e_{\beta ^{'+}_{p_r}}=[K ^{'+}_{p}:L^+_p]$
where $f_{\beta ^{'+}_{p_r}}$ is the residue degree of
${\beta ^{'+}_{p_r}}$ in the extension
$K^{'+}_p/L^+_p$ and $e_{\beta ^{'+}_{p_r}}$ is the corresponding ramification index.
\vskip 11pt

{\bbf A coefficient semiring homomorphism\/} \cite{Maz2}
\[ h_{B ^{'+}_{p}\to B ^{+}_{p}}: \qquad
B ^{'+}_{p}\To B^+_p\]
sending $ B ^{+}_{p}$ into 
 $B ^{'+}_{p}$ {\bbf is such that\/}:
 \Bean
 \item the inverse image of the maximal ideal
$\wt\omega _{K^+_p}\BB^+_p$ of $B^+_p$
is the maximal ideal
$\wt\omega _{K^+_p}B^{'+}_p$ of 
$B^{'+}_p$;

\item {\bbf the induced homomorphism\/}
\[ h_{k_{K^{'+}_p}\to k_{K^+_p}}: \qquad  k_{K^{'+}_p}
\overset{\sim}{\To} k_{K^+_p}\]
{\bbf on the residue semifields is an isomorphism\/} leading to the evident condition
\[\sum_rf_{\beta ^{'+}_{p_r}}=\sum_rf_{\beta ^+_{p_r}}\]
on the residue degrees.
\Ee
\vskip 11pt

\subsubsection{Proposition}

{\em The kernel of the coefficient semiring homomorphism
$ h_{B ^{'+}_{p}\to B ^{+}_{p}}:
B ^{'+}_{p}\to B^+_p$ is characterized by a degree of extension:
\begin{align*}
[K^{'+}_p:L^+_p]-[K^+_p:L^+_p]
&= \sum_rf_{\beta ^{'+}_{p_r}} e_{\beta ^{'+}_{p_r}}-
\sum_rf_{\beta ^{+}_{p_r}} e_{\beta ^{+}_{p_r}}\\
&=(e_{\beta ^{'+}_{p_r}}-e_{\beta ^{+}_{p_r}})
\L(\sum_rf_{\beta ^{+}_{p_r}}\R)\;.
\end{align*}
}
\vskip 11pt

\begin{proof} This is a consequence of section 4.1.1 leading generally to the inequality
\[\sum_rf_{\beta ^{'+}_{p_r}}  e_{\beta ^{'+}_{p_r}}
>\sum_rf_{\beta ^{+}_{p_r}}  e_{\beta ^{+}_{p_r}}\;.\]
On the other hand, as
$ \sum_rf_{\beta ^{'+}_{p_r}}=\sum_rf_{\beta ^{+}_{p_r}}$ and as the ramification indices 
$e_{\beta ^{'+}_{p_r}}$ are equal to $e_{\beta ^{'+}_{p}}$ and the
$e_{\beta ^{+}_{p_r}}$ are equal to $e_{\beta ^{+}_{p}}$ according to section  3.3.1, the preceding inequality results from the fact that $e_{\beta ^{'+}_{p}}>e_{\beta ^{+}_{p}}$.\end{proof}
\vskip 11pt

\subsubsection{Corollary}

{\em The coefficient semiring homomorphism
$ h_{B ^{'+}_{p}\to B ^{+}_{p}}:
B ^{'+}_{p}\to B^+_p$ corresponds to a base change from $K^+_p$ to $K^{'+}_p$.
}
\vskip 11pt

\begin{proof} The degree of this change of basis is thus:
\[ [K^{'+}_p:L^+_p] -[K^+_p:L^+_p]=q'-q\]
if $q=e_{\beta ^+_p}\sum\limits_rf_{\beta ^+_{p_r}}$ and
 $q'=e_{\beta ^{'+}_p}\sum\limits_rf_{\beta ^{'+}_{p_r}}$.
 \end{proof}
 \vskip 11pt
 
 \subsubsection{Global coefficient semiring homomorphisms}
 
 A global coefficient semiring
 $\wt L_{L_p}$
 \resp{$\wt L_{R_p}$} is, according to section 2.1, a complete noetherian global semiring characterized by a set of embedded subsemifields above ``$p$'':
 \[
 \wt L_{v_p} \subset \dots \subset \wt L_{v_{p+h}} \subset \dots\qquad 
 \rresp{\wt L_{\o v_p} \subset \dots \subset \wt L_{\o v_{p+h}} \subset \dots}\]
 which, being compactified, give rise to the corresponding infinite pseudo-ramified completions:
 \begin{align*}
L_{v_p} \subset \dots \subset  L_{v_{p+h}} \subset \dots \in L_{v_p}\in L_{L_p} & \equiv L_{[v_p]}\;, && 1\le h\le \infty \;, \\
\rresp{L_{\o v_p} \subset \dots \subset  L_{\o v_{p+h}} \subset \dots \in L_{\o v_p}\in L_{R_p} & \equiv L_{[\o v_p]}}.
\end{align*}
Referring to proposition 3.3.2, the number of elements of
$L_{L_p}$
\resp{$L_{R_p}$}, a more manageable notation than
$L_{[v_p]}$
\resp{$L_{[\o v_p]}$}, is:
\[ |L_{L_p}|=|L_{R_p}|=\#\Nu\times N\times \sum_h\sum_{m_{p+h}}
f_{v_{p+h,m_{p+h}}}\]
while the number of elements of the corresponding global unramified compactified coefficient semiring
$L^{nr}_{L_p}$
\resp{$L^{nr}_{R_p}$} is:
\[ |L^{nr}_{L_p}|=|L^{nr}_{R_p}|=\#\Nu\times \sum_h\sum_{m_{p+h}}
f_{v_{p+h,m_{p+h}}}\]
where $N$ is the order of the global inertia subgroup(s).
\vskip 11pt

Let 
$L'_{L_p}$
\resp{$L'_{R_p}$} denote another global compactified coefficient semiring characterized by the same set of embedded infinite pseudo-ramified completions
 \[
 L'_{v_p} \subset \dots \subset  L'_{v_{p+h}} \subset \dots\in L'_{L_p} \qquad 
\rresp{L'_{\o v_p} \subset \dots \subset  L'_{\o v_{p+h}} \subset \dots\in L'_{R_p} }\]
of which number of elements is
\[ |L'_{L_p}|=|L'_{R_p}|=\#(\Nu)'\times N'\times \sum_h\sum_{m_{p+h}}
f'_{v_{p+h,m_{p+h}}}\;.\]
$L'_{L_p}$ \resp{$L'_{R_p}$} then differs from 
$L_{L_p}$ \resp{$L_{R_p}$} by the number of non units $\#(\Nu)'$ and by the order $N'$ of the inertia subgroup, the (unramified) maximal orders being by hypothesis the same in
$L_{L_p}$ \resp{$L_{R_p}$} and in
$L'_{L_p}$ \resp{$L'_{R_p}$} and characterized essentially by the global residue degrees
$f_{v_{p+h,m_{p+h}}}$ as developed subsequently.
\vskip 11pt

{\bbf A coefficient semiring homomorphism (isomorphism)}
\[ h_{L'_{L_p}\to L_{L_p}}: \qquad L'_{L_p}\To L_{L_p}\]
{\bbf induces a homomorphism \resp{an isomorphism}\/}
\begin{align*}
 h_{L^{'nr}_{L_p}\to L^{nr}_{L_p}}: \qquad L^{'nr}_{L_p}&\To L^{nr}_{L_p}\\
 \rresp{i_{L^{'nr}_{L_p}\to L^{nr}_{L_p}}: \qquad L^{'nr}_{L_p}&\overset{\sim}{\To} L^{nr}_{L_p}}
 \end{align*}
 {\bbf on the global unramified compactified left coefficient semirings 
$L^{'nr}_{L_p}$ and $L^{nr}_{L_p}$\/} (which corresponds in characteristic $0$ to (global) residue semifields) at the condition that:
\begin{align*}
\sum_h\sum_{m_{p+h}} f'_{v_{p+h,m_{p+h}}}
&= \sum_h\sum_{m_{p+h}}f_{v_{p+h,m_{p+h}}}\\
\rresp{
\#(\Nu)'\times \sum_h\sum_{m_{p+h}} f'_{v_{p+h,m_{p+h}}}
&=\#\Nu\times \sum_h\sum_{m_{p+h}} f_{v_{p+h,m_{p+h}}}}\;.\end{align*}
\vskip 11pt

\subsubsection{Proposition (left case)}

{\em The kernel of the coefficient semiring homomorphism:
\[ h_{L'_{L_p}\to L_{L_p}}: \qquad L'_{L_p}\To L_{L_p}\;,\]
inducing the homomorphism $h_{L^{'nr}_{L_p}\to L^{nr}_{L_p}}$
\resp{the isomorphism $i_{L^{'nr}_{L_p}\to L^{nr}_{L_p}}$}
on their global residue semifields, is characterized by an extension degree:
\[
[L'_{L_p}:k]-[L_{L_p}:k]=(N'-N)\times \L(\sum_h\sum_{m_{p+h}} f_{v_{p+h,m_{p+h}}}\R)\]
and a number of elements:
\begin{align*}
 |L'_{L_p}|-|L_{L_p}|
&= [ (\#(\Nu)'\times N')-(\#\Nu\times N)]\times \sum_h\sum_{m_{p+h}}f_{v_{p+h,m_{p+h}}}\\
\rresp{ |L'_{L_p}|-|L_{L_p}|
&= (N'-N)\times\#\Nu
\times \sum_h\sum_{m_{p+h}}f_{v_{p+h,m_{p+h}}}}\;.
\end{align*}
}
\vskip 11pt

\begin{proof}
This is a consequence of section 4.1.4 restricted to the left case, the right case being handled similarly.\end{proof}
\vskip 11pt

\subsubsection{Corollary}

{\em The coefficient semiring homomorphism
\[ h_{L'_{L_p}\to L_{L_p}} : \qquad L'_{L_p}\To L_{L_p}\]
corresponds to a base change from $L_{L_p}$ into $L'_{L_p}$ of which degree is:
\[
[L'_{L_p}:k]-[L_{L_p}:k]=(N'-N)\times \sum_h\sum_{m_{p+h}} f_{v_{p+h,m_{p+h}}}\;.\]
}
\vskip 11pt

\subsubsection{Proposition (left case)}

{\em {\bbf The inverse image of the homomorphism
$h_{L'_{L_p}\to L_{L_p}}$ between global coefficient semirings is isomorphic to the inverse image of the homomorphism
$h_{B^{'+}_p\to B^+_p}:B^{'+}_p\to B^+_p$ between local coefficient semirings if:\/}
\Be
\item  {\bbf the number of elements of the kernel
$K(h_{L'_{L_p}\to L_{L_p}})$
 is equal to the number of elements of the kernel $K(h_{B^{'+}_p\to B^+_p})$\/},
 i.e. if
\[ |K(h_{L'_{L_p}\to L_{L_p}})|=|K(h_{B^{'+}_p\to B^+_p})|\]
given by:
\[(N'-N)\times \#\Nu\times \sum_h\sum_{m_{p+h}} f_{v_{p+h,m_{p+h}}}=p^{q'-q}\]
where $q=e_{\beta^+_p}\sum_r f_{\beta ^+_{p_r}}$ and
 $q'=e_{\beta^{'+}_p}\sum_r f_{\beta ^{'+}_{p_r}}$;
 
 \item $B^+_p$ covers isomorpically $L_{L_p}$.
 \Ee
 }
 \vskip 11pt
 
 \begin{proof}
 In order that the inverse image of 
 $h_{L'_{L_p}\to L_{L_p}}$ be isomorphic to the inverse image of 
$h_{B^{'+}_p\to B^+_p}$, it is necessary that the induced homomorphism on the respective residue semifields be an isomorphism as it was seen in sections 4.1.1 and 4.1.4 where
$h_{L^{'nr}_{L_p}\to L^{nr}_{L_p}}$ must be the isomorphism
$i_{L^{'nr}_{L_p}\overset{\sim}{\to} L^{nr}_{L_p}}$.
\vskip 11pt

On the other hand, the inverse image of the homomorphism
$h_{L'_{L_p}\to L_{L_p}}$ is given by:
\[ L'_{L_p}=L_{L_p}+K\L(h_{L'_{L_p}\to L_{L_p}}\R)\]
with respect to its kernel $K\L(h_{L'_{L_p}\to L_{L_p}}\R)$ and the inverse image of the homomorphism 
$h_{B^{'+}_p\to B^+_p}$ is similarly given by:
\[ B^{'+}_p=B^+_p+K\L(h_{B^{'+}_p\to B^+_p}\R)\;.\qedhere\]
\end{proof}
\vskip 11pt

\subsubsection{Proposition (left case)}

{\em
\Be
\item To each global coefficient semiring $L_{L_p}$ corresponds
{\bbf the category $\Cs(L'_{L_{p_c}})$, $1\le c\le\infty $, associated with the set $\{h_{L'_{L_{p_c}}\to L_{L_p}}\}_{N'_c}$ of coefficient semiring homomorphisms\/} in such a way that:
\Be
\item the extension degree $[L'_{L_{p_c}}:k]$  of $L'_{L_{p_c}}$ differ from the extension degree of $L_{L_p}$ by the orders ``$N'_c$'' of the inertia subgroups, $N'_c\neq N$.

\item the set $\{f_{v_{p+h,m_{p+h}}}\}_h$ of global residue degrees is an invariant in $L_p$ and in the set $\{L'_{L_{p_c}}\}_{N'_c}$.
\Ee

\item Similarly, to each local coefficient semiring $B^+_p$ corresponds {\bbf the category $\Cs(B^{'+}_{p_d})$, $1\le d\le \infty $, associated with the set
$\{h_{B^{'+}_{p_d}\to B^+_p}\}_{e_{\beta ^{'+}_d}}$ of coefficient semiring homomorphisms\/} in such a way that:
\Be
\item the extension degrees $[K^{'+}_{p_d}:L_p]$ of $K^{'+}_{p_d}$ differ from the extension degree of $K^+_p$ by the ramification indices $e_{\beta ^{'+}_{p_d}}$ different from $e_{\beta ^+_p}$.

\item the set $\{f_{\beta ^{'+}_{p_r}}\}_r$ of residue degrees is an invariant in $K^+_p$ and in the set $\{K^{'+}_{p_d}\}_d$.
\Ee\Ee
}
\vskip 11pt

\begin{proof}
\Bean
\item The existence of the category $\Cs(L'_{L_{p_c}})$ whose objects are global coefficient semirings $L'_{L_{p_c}}$ and whose morphisms are the homomorphisms $h_{L'_{L_{p_c}}\to L_{L_p}}$ results from the orders $N'_c$ of the inertia subgroups of $L'_{L_{p_c}}$, $c\in \NN$, $1\le c\le \infty $.

\item In the same way, the existence of the category $\Cs(B^{'+}_{p_d})$, $1\le d\le \infty $, whose objects are local coefficient semirings $B^{'+}_{p_d}$ and whose morphisms are homomorphisms
$h_{B^{'+}_{p_d}\to B^+_p}$ results from the ramification indices
$e_{\beta ^{'+}_{p_d}}$ of $B^{'+}_{p_d}$ different from the ramification index $e_{\beta ^+_p}$ of $B^+_p$.\qedhere
\Ee
\end{proof}
\vskip 11pt

\subsection{Deformations of Galois representations over local and global noetherian bisemirings}

\subsubsection{Definition (Global bilinear deformation representative)}

A global bilinear deformation representative resulting from a global bilinear coefficient semiring homomorphism
\[ h_{L'_{R_p}\times L'_{L_p}\to L_{R_p}\times L_{L_p}} : \qquad
L'_{R_p}\times L'_{L_p}\To L_{R_p}\times L_{L_p}\;, \]
inducing the isomorphism
$i_{L^{'nr}_{R_p}\times L^{'nr}_{L_p}}\to L^{nr}_{R_p}\times L^{nr}_{L_p}$ 
on their global bilinear residue semifields is {\bbf an equivalence class representative $\rho'_L$ of lifting}
\[\begin{CD}
\Gal(\dot{\wt L}'_{R_p}/k)\times\Gal(\dot{\wt L}'_{L_p}/k)
@>{h_{L'\to L}}>> \Gal(\dot{\wt L}_{R_p}/k)
\times \Gal(\dot{\wt L}_{L_p}/k)\\
@VV{\rho _{L'}}V  @VV{\rho _{L}}V \\
\GL_n(L'_{R_p}\times L'_{L_p}) @>{h_{G'\to G}}>> \GL_n(L_{R_p}\times L_{L_p})
\end{CD}\]
with the evident bilinear notations introduced in sections 4.1.4 and 4.1.5 for the right and left global semiring homomorphisms, in section 2.4 for the Galois and Weil groups and in section 2.5 for the algebraic bilinear semigroups 
$\GL_n(L'_{R_p}\times L'_{L_p})$ and
$\GL_n(L_{R_p}\times L_{L_p})$.
\vskip 11pt

\subsubsection{Proposition (Global bilinear deformation)}

{\em Let $K(h_{L'_{R_{p_c}}\times L'_{L_{p_c}}\to
L_{R_{p_c}}\times L_{L_{p_c}}})$ 
be the kernel of the bihomomorphism $h_{L'_{R_{p_c}}\times L'_{L_{p_c}}\to
L_{R_{p_c}}\times L_{L_{p_c}}}$
where $L'_{R_{p_c}}\times 
L'_{L_{p_c}}$ belongs to the bicategory
$\Cs(L'_{R_{p_c}}\times L'_{L_{p_c}})$, $1\le c\le \infty $, defined similarly as for the left case in proposition 4.1.8.

Let 
\[
\Gal(\delta \dot{\wt L}'_{R_{p_c}}/k)
\times
\Gal(\delta \dot{\wt L}'_{L_{p_c}}/k)
=
[\Gal( \dot{\wt L}'_{R_{p_c}}/k)-
\Gal( \dot{\wt L}'_{R_{p}}/k)]
\times
[\Gal( \dot{\wt L}'_{L_{p_c}}/k)-
\Gal( \dot{\wt L}'_{L_{p}}/k)]\]
be the Weil subgroup associated with this kernel 
$K(h_{L'_{R_{p_c}}\times
L'_{L_{p_c}}}\to
L_{R_{p_c}}\times
L_{L_{p_c}})$.

Then, {\bbf a $n$-dimensional global bilinear deformation of $\rho _L$ is an equivalence class of liftings $\{\rho _{L'_c}\}$, $1\le c\le\infty $, described by the following commutative exact sequence\/}
\[
\begin{CD}
{1 \quad \to \begin{array}[t]{r}
\Gal(\delta \dot{\wt L}'_{R_{p_c}}/k)\quad \\[-8pt]
\times
\Gal(\delta \dot{\wt L}'_{L_{p_c}}/k)
\end{array}}
@>>>
\begin{array}[t]{r}
\Gal( \dot{\wt L}_{R_{p_c}}/k)\quad\\[-8pt]
 \times
\Gal(  \dot{\wt L}_{L_{p_c}}/k)
\end{array}
@>{h_{L'_c\to L}}>>
\begin{array}[t]{r}
\Gal(  \dot{\wt L}_{R_{p }}/k)\quad\\[-8pt]
 \times
\Gal( \dot{\wt L}_{L_{p }}/k)\end{array}\to 1\\
@VV{\delta \rho _{L'_c}}V
@VV{\rho _{L'_c}}V
@VV{\rho _{L}}V\\
{1\to \GL_n(\delta { L}'_{R_{p_c}}\times\delta {L}'_{L_{p_c}})}
@>>>
{\GL_n( { L}'_{R_{p_c}}\times {L}'_{L_{p_c}})} 
@>>{h_{G'_c\to G}}>
{\GL_n( { L}_{R_{p}}\times {L}_{L_{p}})\to 1} 
\end{CD}\]
of which ``Weil kernel'' is 
$\Gal(\delta \dot{\wt L}'_{R_{p_c}}/k)\times
\Gal(\delta \dot{\wt L}'_{L_{p_c}}/k)$ and
``$\GL_n(\pt\times\pt)$ kernel'' is\linebreak
$\GL_n(\delta {L}'_{R_{p_c}}\times\delta { L}'_{L_{p_c}})$
where $(\delta {L}'_{R_{p_c}}\times\delta { L}'_{L_{p_c}})$
is given by
\[
(\delta {L}'_{R_{p_c}}\times\delta { L}'_{L_{p_c}})=[(
 {L}'_{R_{p_c}}-{ L}_{R_{p}})\times
 ({L}'_{L_{p_c}}-{ L}_{L_{p}} )]\;.\]
 }
 \vskip 11pt

\begin{proof}
A $n$-dimensional global bilinear deformation of $\rho _L$ is thus  an equivalence class of liftings
\[ \rho _{L'_c}=\rho _L+\delta \rho _{L'_c}\qquad\forall\ c\;, \quad 1\le c\le \infty \;,\]
where two liftings $\rho _{L'_{c_1}}$ and $\rho _{L'_{c_2}}$ are strictly equivalent if they can be transformed one into another by conjugation by bielements of 
$\GL_n( {L}'_{R_{p_c}}\times{ L}'_{L_{p_c}})$ in the kernel
$\GL_n(\delta  {L}'_{R_{p_c}}\times\delta { L}'_{L_{p_c}})$ of $h_{G'_c\to G}$: this is the worked out condition of (bilinear) deformations proposed by B. Mazur for example in \cite{Maz2}.

In other words, the lifting  $\rho _{L'_{c_1}}$ generates the bilinear algebraic semigroup
$\GL_n( {L}'_{R_{p_{c_1}}}\times{ L}'_{L_{p_{c_1}}})$ having a (presumed) rank
\[r_{G_n(L'_{R_{c_1}}\times L'_{L_{c_1}})}
=((N'_{c_1})^n\cdot f^n_{\o v})\times ((N'_{c_1})^n\cdot f^n_v)
\]
where:
\Bi
\item $N'_{c_1}$ is the order of the global inertia subgroup according to section 4.1.4;
\item $f_v$ ($\equiv f_{\o v}$) is a condensed notation for
$\sum\limits_h\sum\limits_{m_{p+h}}f_{v_{p+h,m_{p+h}}}$.
\Ei
And, the lifting $\rho _{L'_{c_2}}$ generates the algebraic bilinear semigroup
$\GL_n(L'_{R_{p_{c_2}}}\times L'_{L_{p_{c_2}}})$ of which rank
\[r_{G_n(L'_{R_{c_2}}\times L'_{L_{c_2}})}
=((N'_{c_2})^n\cdot f^n_{\o v})\times ((N'_{c_2})^n\cdot f^n_v)
\]
differs from the rank
$r_{G_n(L'_{R_{c_1}}\times L'_{L_{c_1}})}$ of
$\GL_n(L'_{R_{c_1}}\times L'_{L_{c_1}})$ by
\[\delta r_{G_{n_{(2-1)}}}
=[(N'_{c_2})^{n^2}-(N'_{c_1})^{n^2}](f_v)^{n^2}\]
taking into account that $f_v=f_{\o v}$ if we refer to proposition 4.1.5.

As a consequence, the deformed algebraic bilinear semigroups
$\GL_n(L'_{R_{p_{c_1}}}\times L'_{L_{p_{c_1}}})$ and\linebreak
$\GL_n(L'_{R_{p_{c_2}}}\times L'_{L_{p_{c_2}}})$ differ from the algebraic bilinear semigroup
$\GL_n(L_{R_{p}}\times L_{L_{p}})$ by their respective kernels
$\GL_n(\delta L'_{R_{p_{c_1}}}\times \delta L'_{L_{p_{c_1}}})$ and
$\GL_n(\delta L'_{R_{p_{c_2}}}\times \delta L'_{L_{p_{c_2}}})$ in such a way that
$\GL_n(L'_{R_{p_{c_1}}}\times L'_{L_{p_{c_1}}})$
can be transformed into
$\GL_n(L'_{R_{p_{c_2}}}\times L'_{L_{p_{c_2}}})$ by conjugation of bielements in the first kernel
$\GL_n(\delta L'_{R_{p_{c_1}}}\times \delta L'_{L_{p_{c_1}}})$ bringing it into the second kernel
$\GL_n(\delta L'_{R_{p_{c_2}}}\times\delta  L'_{L_{p_{c_2}}})$ and inversely.
\end{proof}
\vskip 11pt

\subsubsection{Corollary (The transformation of kernels)}

{\em \[\GL_n(\delta L'_{R_{p_{c_1}}}\times \delta L'_{L_{p_{c_1}}})
\To
\GL_n(\delta L'_{R_{p_{c_2}}}\times \delta L'_{L_{p_{c_2}}})\]
{\bbf corresponds to a base change\/} of
$\GL_n( L'_{R_{p_{c_1}}}\times  L'_{L_{p_{c_1}}})$
into
$\GL_n( L'_{R_{p_{c_2}}}\times  L'_{L_{p_{c_2}}})$
whose dimension is given by the difference of ranks
\[\delta r_{G_{n_{(2-1)}}}
=[(N'_{c_2})^{n^2}-(N'_{c_1})^{n^2}](f_v)^{n^2}\;.\]
}
\vskip 11pt

\begin{proof}
This results directly from proposition 4.2.2.
\end{proof}
\vskip 11pt

As the local semirings of category
$\Cs(B^{'+}_{p_d})$ differ between themselves by their ramification indices and as the global semirings of the category
$\Cs(L'_{p_c})$ differ between themselves by their orders of inertia subgroups, the local (bilinear) deformations can be described similarly as it was done for the global bilinear deformations.
\vskip 11pt

\subsubsection{Definition (Local bilinear deformation representative)}

A local bilinear deformation representative resulting from a local bilinear coefficient semiring homomorphism:
\[
h_{B^{'-}_p\times B^{'+}_p\to B^{-}_p\times B^{+}_p}: \qquad 
{B^{'-}_p\times B^{'+}_p\To B^{-}_p\times B^{+}_p}\;,\]
inducing an isomorphism on their residue (bi)semifields, is
{\bbf an equivalence class representative $\rho _{K'}$ of lifting\/}
\[\begin{CD}
\Gal(K^{'-}_p/L^-_p)\times\Gal(K^{'+}_p/L^+_p)
@>>{h_{K'\to K}}>
\Gal(K^{-}_p/L^-_p)\times\Gal(K^{+}_p/L^+_p)\\
@VV{\rho _{K'}}V @VV{\rho _{K}}V\\
\GL_n(K^{'-}_p\times K^{'+}_p)
@>>{h_{\GL'\to \GL}}>
\GL_n(K^{-}_p\times K^{+}_p)
\end{CD}\]
where:
\Bi
\item $K^+_p$ 
\resp{$K^-_p$} is a finite \lr $p$-adic Galois extension of the \lr $p$-adic semifield 
$L^+_p$
\resp{$L^-_p$} according to section 3.3.1 in such a way that 
$B^+_p$
\resp{$B^-_p$} be the (valuation) coefficient semiring in
$K^+_p$
\resp{$K^-_p$} referring to section 4.1.1.

\item $K^{'+}_p$
\resp{$K^{'-}_p$} is another finite \lr $p$-adic Galois extension of 
$L^+_p$
\resp{$L^-_p$} in such a way that the extension degree $q'$ of
$K^{'+}_p$
\resp{$K^{'-}_p$} is superior to the extension degree $q$ of
$K^+_p$
\resp{$K^-_p$} referring to section 4.1.1 and proposition 4.1.2.
\Ei
\vskip 11pt

\subsubsection{Proposition (Local bilinear deformation)}

{\em Let
$K(
h_{B^{'-}_{p_d}\times B^{'+}_{p_d}\to B^{-}_p\times B^{+}_p})$ be the kernel of the bihomomorphism
$h_{B^{'-}_{p_d}\times B^{'+}_{p_d}\to B^{-}_p\times B^{+}_p}$ where
$B^{'-}_{p_d}\times B^{'+}_{p_d}$ belongs to the bicategory
$\Cs(B^{'-}_{p_d}\times B^{'+}_{p_d})$, $1\le d\le\infty $, defined similarly as for the left case in proposition 4.1.8.

Let 
\[\Gal(\delta K^{'-}_{p_d}/L^-_p)\times \Gal(\delta K^{'+}_{p_d}/L^+_p)=
[\Gal( K^{'-}_{p_d}/L^-_p)- \Gal( K^{-}_{p}/L^-_p)]\times
[\Gal( K^{'+}_{p_d}/L^+_p)- \Gal( K^{+}_{p}/L^+_p)]\]
be the Weil subgroup corresponding to this kernel
$K(
h_{B^{'-}_{p_d}\times B^{'+}_{p_d}\to B^{-}_p\times B^{+}_p})$.

Then, {\bbf a $n$-dimensional local bilinear deformation of $\rho _K$ is an equivalence class of liftings
$\{\rho _{K'_d}\}$ described by the following commutative exact sequence\/}:
\[
\begin{CD}
{1 \quad \to \begin{array}[t]{r}
\Gal(\delta{ K}^{'-}_{p_d}/L^-_p)\quad \\[-8pt]
\times
\Gal(\delta{ K}^{'+}_{p_d}/L^+_p)
\end{array}}
@>>>
\begin{array}[t]{r}
\Gal({ K}^{'-}_{p_d}/L^-_p)\quad\\[-8pt]
 \times
\Gal({ K}^{'+}_{p_d}/L^+_p)\end{array}
@>{h_{K'_d\to K}}>>
\begin{array}[t]{r}
\Gal({ K}^{-}_{p}/L^-_p)\quad\\[-8pt]
 \times
\Gal({ K}^{+}_{p}/L^+_p)\end{array}\to 1\\
@VV{\delta \rho _{K'_d}}V
@VV{\rho _{K'_d}}V
@VV{\rho _{K}}V\\
{1\to \GL_n ( \delta { K}^{'-}_{p_d}\times\delta {K}^{'+}_{p_d} )}
@>>>
{\GL_n (  { K}^{'-}_{p_d}\times {K}^{'+}_{p_d} )} 
@>>{h_{\GL'_d\to \GL}}>
{\GL_n (  { K}^{-}_{p}\times {K}^{+}_{p} )\to 1} 
\end{CD}\]
of which ``Weil kernel'' is 
$\Gal({\delta K}^{'-}_{p_d}/L^-_p)
 \times
\Gal(\delta{ K}^{'+}_{p_d}/L^+_p)$
 and
``$\GL_n(\pt\times\pt)$ kernel'' is
$\GL_n ( \delta { K}^{'-}_{p_d}\times\delta {K}^{'+}_{p_d} )$
where 
\[
 ( \delta { K}^{'-}_{p_d}\times\delta {K}^{'+}_{p_d} )=
[ (  { K}^{'-}_{p_d}-K^-_p)\times(\delta {K}^{'+}_{p_d}-K^+_p)]
\;.\]
}
\vskip 11pt

\begin{proof}
A $n$-dimensional local bilinear deformation of $\rho _k$ is an equivalence class of liftings
\[
\rho _{K'_d}= \rho _{K}+\delta \rho _{K'_d}\;, \qquad \forall\ d\;, \quad 1\le d\le\infty \;, \]
where two local liftings $\rho _{K'_{d_1}}$ and
 $\rho _{K'_{d_2}}$ are strictly equivalent if they can be transformed one into another by conjugation of bielements of 
 $\GL_n(K^{'-}_{p_{d_1}}\times K^{'+}_{p_{d_1}})$ in the local kernel
 $\GL_n(\delta K^{'-}_{p_{d_1}}\times \delta K^{'+}_{p_{d_1}})$ of
 $h_{\GL '_d\to\GL}$ \cite{Maz2}.
 
 In other words, the lifting
 $\rho _{K'_{d_1}}$ generates the bilinear algebraic semigroup
 $\GL_n(K^{'-}_{p_{d_1}}\times K^{'+}_{p_{d_1}})$ having a rank
 \[
 r_{\GL_n (K^{'-}_{p_{d_1}}\times K^{+}_{p_{d_1}})}
 =((e_{\beta ^{'-}_{p_{d_1}}})^n\times (f_{\beta ^{'-}_{p_{d_1}}})^n)
 \times
((e_{\beta ^{'+}_{p_{d_1}}})^n\times (f_{\beta ^{'+}_{p_{d_1}}})^n)
=(e_{\beta '_{p_{d_1}}})^{n^2}\times (f_{\beta '_{p_{d_1}}})^{n^2}\]
where:
\Bi
\item $e_{\beta ^{'-}_{p_{d_1}}}=e_{\beta ^{'+}_{p_{d_1}}}$ is the ramification index of 
$K^{'-}_{p_{d_1}}$ or of
$K^{'+}_{p_{d_1}}$.

\item $f_{\beta ^{'-}_{p_{d_1}}}=f_{\beta ^{'+}_{p_{d_1}}}=
\sum\limits_r f_{\beta ^{'-}_{p_{r_{d_1}}}}$ is the residue degree of
$K ^{'-}_{p_{d_1}}$ or of $K ^{'+}_{p_{d_1}}$.
\Ei

And,  the lifting $\rho _{K ^{'}_{p_{d_2}}}$ generates the algebraic bilinear semigroup 
$\GL_n(K^{'-}_{p_{d_2}}\times K^{'+}_{p_{d_2}})$ of which rank 
\[
 r_{\GL_n (K^{'-}_{p_{d_2}}\times K^{'+}_{p_{d_2}})}
 =(e_{\beta '_{p_{d_2}}})^{n^2}\times (f_{\beta '_{p_{d_2}}})^{n^2}\]
 differs from the rank
 $
 r_{\GL_n (K^{'-}_{p_{d_1}}\times K^{+}_{p_{d_1}})}$
 of ${\GL_n (K^{'-}_{p_{d_1}}\times K^{+}_{p_{d_1}})}$ by
 \[\delta r_{\GL_n{{(2-1)}}}
=[(e_{\beta  '_{p_{d_2}}})^{n^2}-
(e_{\beta '_{p_{d_1}}})^{n^2}]
(f_{\beta  '_{p_{d}}})^{n^2})]
\]
because all the residue degrees of 
$f_{\beta '_{p_{d}}}$ are equal since the homomorphisms on the envisaged residue semifields are isomorphisms according to section 4.1.1.

Consequently, the deformed algebraic bilinear semigroups
$\GL_n (K^{'-}_{p_{d_1}}\times K^{'+}_{p_{d_1}})$ and
$\GL_n (K^{'-}_{p_{d_2}}\times K^{'+}_{p_{d_2}})$ differ from the algebraic bilinear semigroup
$\GL_n (K^{-}_{p}\times K^{+}_p)$ by their respective kernels
$\GL_n (\delta K^{'-}_{p_{d_1}}\times \delta K^{'+}_{p_{d_1}})$ and
$\GL_n (\delta K^{'-}_{p_{d_2}}\times \delta K^{'+}_{p_{d_2}})$ in such a way that
$\GL_n (K^{'-}_{p_{d_1}}\times K^{'+}_{p_{d_1}})$ can be transformed into
$\GL_n (K^{'-}_{p_{d_2}}\times K^{'+}_{p_{d_2}})$ by conjugation of bielements in the first kernel bringing it to the second kernel and inversely.
\end{proof}
\vskip 11pt

\subsubsection{Corollary (The transformation of kernels)}

{\em 
\[ \GL_n (\delta K^{'-}_{p_{d_1}}\times \delta K^{'+}_{p_{d_1}})
\To
\GL_n (\delta K^{'-}_{p_{d_2}}\times \delta K^{'+}_{p_{d_2}})\]
{\bbf corresponds to a base change\/} of 
$\GL_n ( K^{'-}_{p_{d_1}}\times  K^{'+}_{p_{d_1}})$
into
$\GL_n ( K^{'-}_{p_{d_2}}\times K^{'+}_{p_{d_2}})$
of which dimension is given by the difference of ranks
\[ \delta r_{\GL_n(2-1)}=
[(e_{\beta '_{p_{d_2}}})^{n^2}-
(e_{\beta '_{p_{d_1}}})^{n^2}]
(f_{\beta '_{p_{d}}})^{n^2})]\;.\]
}
\vskip 11pt

\begin{proof}
This results from proposition 4.2.5 similarly as for the global case.
\end{proof}
\vskip 11pt

\subsection{Global and local coefficient semiring quantum homomorphisms}

\subsubsection{Global coefficient semiring quantum homomorphisms (left case)}

Let $L_{L_p}$
\resp{$L_{R_p}$} denote a global compactified coefficient semiring characterized by a set of embedded infinite pseudo-ramified completions above ``$p$'':
\begin{align*}
& L_{v_p} \subset \dots \subset L_{v_{p+h}} \subset \dots 
\qquad \rresp{&L_{\o v_p} \subset \dots \subset L_{\o v_{p+h}} \subset \dots }&& 1\le h\le\infty \;.\end{align*}
\vskip 11pt

Let $L_{L_{p+j}}$
\resp{$L_{R_{p+j}}$} denote another global compactified coefficient semiring composed of the same number of corresponding embedded pseudo-ramified completions above ``$p+j$'':
\[ 
L_{v_{p+j}} \subset \dots \subset L_{v_{p+j+h}} \subset \dots \qquad
\rresp{L_{\o v_{p+j}} \subset \dots \subset L_{\o v_{p+j+h}} \subset \dots }\]
where:
\Bi
\item $q=p+j$ is assumed to be a prime number;
\item the global residue degree of 
$L_{v_{p+j+h}}$ (and
$L_{\o v_{p+j+h}}$) in
$L_{L_{p+j}}$
\resp{$L_{R_{p+j}}$} is given by:
\[ f_{v_{p+j+h}}=f_{v_{p+h}}+j = p+h+j\]
with respect to its correspondent $L_{v_{p+h}}$ in
$L_{L_p}$
\resp{$L_{R_p}$};

\item the number of non units and the degree of the inertia subgroup are the same in $L_{L_p}$
\resp{$L_{R_p}$} and in 
$L_{L_q}$
\resp{$L_{R_q}$}.
\Ei
\vskip 11pt

Let 
\[ Qh_{L_{L_{p+j}}\to L_{L_p}}: \qquad L_{L_{p+j}}\To L_{L_p}\]
be {\bbf a uniform quantum homomorphism between global compactified coefficient semirings inducing an isomorphism on their global inertia subgroups having the same degree $N$\/}
\vskip 11pt

\subsubsection{Proposition}

{\em
The kernel $K(Qh_{L_{L_{p+j}}\to L_{L_p}})$ of the uniform quantum homomorphism
\[ Qh_{L_{L_{p+j}}\to L_{L_p}}: \qquad {L_{L_{p+j}}\To L_{L_p}}\;,\]
inducing an isomorphism on their global inertia subgroups, is characterized by an extension degree:
\[ [L_{L{p+j}}:k]-[L_{L_p}:k]=N\times j\times \sum_hm_{p+h}\]
if $m_{p+h+j}=m_{p+h}$.
}
\vskip 11pt

\begin{proof}
The extension degree $[L_{L{p+j}}:k]$ is given by:
\[ [L_{L{p+j}}:k]=N\times \sum_h\sum_{m_{p+h+j}}f_{v_{p+h+j}}\;,\]
$j$ being fixed, where $f_{v_{p+h+j}}$ is the global residue degree of the completion $L_{v_{p+h+j}}$, and the extension degree
$[L_{L_p}:k]$ is given by:
\[ [L_{L_p}:k]=N\times\sum_h\sum_{m_{p+h}}
f_{v_{p+h}}\;, \qquad 0\le h\le \infty \;.\]
So, referring to section 4.3.1, we have that:
\begin{align*}
[L_{L{p+j}}:k]-[L_{L_p}:k]
&=\L[ N\times\L(\sum_hm_{p+h+j}\R) \times (p+h+j)\R]
-\L[ N\times\sum_h (m_{p+h}) \times (p+h)\R]\\
&= N\times j\times\sum_hm_{p+h}\;, \end{align*}
where $m_{p+h}$ is the multiplicity of the completion $L_{v_{p+h}}$, if $m_{p+h+j}=m_{p+h}$, condition resulting from the homomorphisms $Qh_{L_{p+j}\to L_{L_p}}$.
\end{proof}
\vskip 11pt

\subsubsection{Corollary}

{\em {\bbf The uniform quantum homomorphim:
\[Qh_{L_{L_{p+j}}\to L_{L_p}}: \qquad L_{L_{p+j}}\To L_{L_p}\]
corresponds to a base change from $L_{L_p}$ into $L_{L_{p+j}}$ of which extension degree
\[ [L_{L_{p+j}}:k]-[L_{L_p}:k]=N\times j\times\sum_hm_{p+h}\]
means an increment of $j$ quanta of degree $N$ on each completion\/} of the global compactified coefficient semiring $L_{L_p}$.
}
\vskip 11pt

\begin{proof}
\Be
\item From the preceding developments, it appears that the kernel
$K(Qh_{L_{L_{p+j}}\to L_{L_p}})$ of the uniform quantum homomorphism
$Qh_{L_{L_{p+j}}\to L_{L_p}}$ measures the extent of the base change from $L_{L_p}$ into ${L_{L_{p+j}}}$.

\item This base change corresponds to an increment of $j$ quanta on each completion of $L_{L_p}$, taking into account that a quantum was defined \cite{Pie3} as a (compact) closed-irreducible algebraic subset of degree $N$.
\qedhere
\Ee
\end{proof}
\vskip 11pt

\subsubsection{Remarks}

\Be
\item {\bbf The quantum homomorphism
$Qh_{L_{L_{p+\{j\}}}\to L_{L_p}}$ will be said non uniform\/} if the completions of $L_{L_p}$ are increased by different numbers of quanta, the integer ``$j$'' then varying from one completion of
$L_{p+\{j\}}$ to another.

\item The uniform quantum homomorphism
$Qh_{L_{L_{p+j}}\to L_{L_p}}$ is in fact induced by a quantum homomorphism
\[ Qh_{ L^{nr}_{p+j}\to L^{nr}_{L_p}}: \qquad L^{nr}_{L_{p+j}}\To L^{nr}_{L_p}\]
between the corresponding unramified global compactified coefficient semirings $L^{nr}_{L_p}$ and $L^{nr}_{L_{p+j}}$.
\Ee
\vskip 11pt

\subsubsection{``Quantum'' homomorphisms between local coefficient semirings (left case)}

Let $B^+_{p_r}$ denote a noetherian local left coefficient semiring on the completion of a finite Galois extension 
$K^+_{p_r}$ of the $p$-adic semifield $L^+_p$.

$B^+_{p_r}$ is a discrete valuation semiring, being the integral closure of $A^+_p$ in $K^+_{p_r}/L^+_p$.  It is characterised by a chain $\beta ^+_{p_1}\subset\dots\subset \beta ^+_{p_r}$ of $r$ prime ideals and has a finite residue semifield $k_{K^+_{p_r}}$ according to section 4.1.1.
\vskip 11pt

Let $B^+_{p_t}$ be another local left coefficient semiring, with $t=r+s$, $r\le t\le s+r$, in such a way that $B^+_{p_t}$ be the integral closure of $A^+_p$ in the finite Galois extension
$K^+_{p_t}/L^+_p$ and that
\[ k_{K^+_{p_t}}=\{\Os_{K^+_{p_t}\mid \beta ^+_{p_t}}\}^r_{t=s}\] 
be its residue semifield defined on the set of $r$ residue semifields $\Os_{K^+_{p_t}\mid \beta ^+_{p_t}}/m(A^+_p)$.

The semiring $B^+_{p_t}/m(A^+_p)B^+_{p_t}$ is also an
$A^+_p/m(A^+_p)$-semialgebra of degree
\[ d_t =[K^+_{p_t}:L^+_p]
=\sum_{\beta ^+_{p_t}\mid p_p}f_{\beta ^+_{p_t}}
e_{\beta ^+_p}\]
 where $f_{\beta ^+_{p_t}}$ is the residue degree of
 $\beta ^+_{p_t}$ in the extension $K^+_{p_t}/L^+_p$ and 
 $e_{\beta ^+_p}$ is the ramification index in the extensions
 $K^+_{p_r}/L^+_p$ and $K^+_{p_t}/L^+_p$.
 \vskip 11pt
 
 {\bbf A coefficient semiring quantum homomorphism\/}
 \[ Qh_{B^+_{p_t}\to B^+_{p_r}}: \qquad B^+_{p_t}\To B^+_{p_r}\]
 sending $B^+_{p_r}$ into $B^+_{p_t}$ is such that:
 \Bean
 \item the inverse image of the maximal ideal
 $\wt \omega _{K^+_{p_r}}$ of $B^+_{p_r}$ is the maximal ideal
$\wt \omega _{K^+_{p_t}}B^+_{p_t}$ of $B^+_{p_t}$.

\item the induced homomorphism
 \[Qh_{k_{K^+_{p_t}}\to k_{K^+_{p_r}}}:\qquad
 k_{K^+_{p_t}}\To k_{K^+_{p_r}}\]
 on the left residue semifields is an isomorphism.
 
 \item there is an isomorphism between the inertia subgroups $K^+_{p_r}$ and $K^+_{p_t}$ having the same ramification index $e_{\beta ^+_p}$.
 \Ee
 \vskip 11pt
 
 \subsubsection{Proposition}
 {\em {\bbf The kernel of the coefficient semiring quantum homomorphism 
  \[ Qh_{B^+_{p_t}\to B^+_{p_r}}: \qquad B^+_{p_t}\To B^+_{p_r}\] is characterized by a degree
  \begin{align*}
  [K^+_{p_t}:L^+_p]-
  [K^+_{p_r}:L^+_p]
  &= e_{\beta ^+_p}\times \L[\sum_{t=r}^{r+s}f_{\beta ^+_{p_t}}-\sum_rf_{\beta ^+_{p_r}}\R]\\
  &= e_{\beta ^+_p}\times \L[f_{k_{K^+_{p_t}}}-f_{k_{K^+_{p_r}}}\R]
  \end{align*}
  and measures the extent of the base change from $K^+_{p_r}$ to $K^+_{p_t}$\/}}.
  \vskip 11pt
  
  \begin{proof}
  This is evident from the preceding developments.
  \end{proof}
  \vskip 11pt
  
  \subsubsection{Categories associated with quantum homomorphisms (left case)}
  
  \Be
  \item {\bbf To each global coefficient semiring $L_{L_p}$ corresponds the category $\Cs(L_{L_{p+j}})$, $1\le j\le\infty $, associated with the set $\{Qh_{L_{p+j}\to L_{L_p}}\}^\infty _{j=1}$ of uniform quantum homomorphisms\/} in such a way that:
\Be
\item the extension degrees $[L_{L_{p+j}}:k]$ of $L_{L_{p+j}}$ differ from the extension degree of
$L_{L_{p}}$ by the numbers of $j$ quanta of degree $N$, $j>0$.

\item the order $N$ of the inertia subgroups is the same in
$L_{L_{p+j}}$ and in $L_{L_{p}}$: it is thus an invariant of
$\Cs(L_{L_{p+j}})$.
\Ee

\item {\bbf To each local coefficient semiring $B^+_{p_r}$ corresponds the category $\Cs(B^+_{p_t})$, $r\le t\le r+s$, associated with the set
$\{Qh_{B^+_{p_t}\to B^+_{p_r}}\}$ of uniform quantum homomorphisms\/} in such a way that:
\Be
\item the extension degrees $[K^+_{p_t}:L^+_p]$ of $K^+_{p_t}$ differ from the extension degree $[K^+_{p_r}:L^+_p]$
of $ K^+_{p_r}$ by the differences 
$[f_{k_{K^+_{p_t}}}- f_{k_{K^+_{p_r}}}]$ in their (local) residue degrees according to proposition 4.3.6.

\item the ramification index $e_{\beta ^+_p}$ is the same in $K^+_{p_t}$  and in $K^+_{p_r}$: it is thus an invariant of
$Qh_{B^+_{p_t}\to B^+_{p_r}}$.
\Ee
\Ee
\vskip 11pt

\subsection{Quantum deformations of Galois representations over local and global noetherian bisemirings}

\subsubsection{Definition (Global bilinear quantum deformation representative)}

A global bilinear quantum deformation representative, resulting from a global bilinear coefficient semiring quantum homomorphism:
\[ Qh_{L_{R_{p+j}}\times L_{L_{p+j}}\to
L_{R_{p}}\times L_{L_{p}}}: \qquad
L_{R_{p+j}}\times L_{L_{p+j}}\To
L_{R_{p}}\times L_{L_{p}}\;,\]
characterized by global inertia subgroups having the same degree $N$, is an equivalence class representative $\rho _{L_j}$ of lifting
\[\begin{CD}
\Gal(\dot{\wt L}_{R_{p+j}}/k)\times\Gal(\dot{\wt L}_{L_{p+j}}/k)
@>{Qh_{L_j\to L}}>> 
\Gal(\dot{\wt L}_{R_p}/k)\times
\Gal(\dot{\wt L}_{L_p}/k)\\
@VV{\rho _{L_j}}V  @VV{\rho _{L}}V \\
\GL_n(L_{R_{p+j}}\times L_{L_{p+j}}) @>{Qh_{G_j\to G}}>> \GL_n(L_{R_p}\times L_{L_p})
\end{CD}\]
with the evident bilinear notations.
\vskip 11pt

\subsubsection{Proposition (Global bilinear quantum deformation)}

{\em
Let $K(Qh_{L_{R_{p+j}}\times L_{L_{p+j}}\to
L_{R_{p}}\times L_{L_{p}}})$ be the kernel of the uniform quantum bihomomorphism\linebreak
$Qh_{L_{R_{p+j}}\times L_{L_{p+j}}\to
L_{R_{p}}\times L_{L_{p}}}$ where
${L_{R_{p+j}}\times L_{L_{p+j}}}$ belongs to the bicategory
$\Cs({L_{R_{p+j}}\times L_{L_{p+j}}})$, $1\le j\le \infty $.

Let
\[
\Gal(\delta \dot{\wt L}_{R_{p+j}}/k)\times\Gal(\delta \dot{\wt L}_{L_{p+j}}/k)=
[\Gal( \dot{\wt L}_{R_{p+j}}/k)-\Gal( \dot{\wt L}_{R_{p}}/k)]\times
[\Gal( \dot{\wt L}_{L_{p+j}}/k)-\Gal( \dot{\wt L}_{L_{p}}/k)]
\]
be the Weil subgroup associated with this kernel.

Then, {\bbf a $n$-dimensional global bilinear quantum deformation of $\rho _L$ is an equivalence class of liftings\/}
$\{\rho _{L_j}\}_j$, $1\le j\le\infty $, described by the following commutative exact sequence:
\[
\begin{CD}
{1 \quad \to \begin{array}[t]{r}
\Gal(\delta \dot{\wt L}_{R_{p+j}}/k)\quad \\[-8pt]
\times\Gal(\delta \dot{\wt L}_{L_{p+j}}/k)
\end{array}}
@>>>
\begin{array}[t]{r}
\Gal( \dot{\wt L}_{R_{p+j}}/k)\quad \\[-8pt]
\times\Gal( \dot{\wt L}_{L_{p+j}}/k)
\end{array}
@>{Qh_{L_j\to L}}>>
\begin{array}[t]{r}
\Gal( \dot{\wt L}_{R_{p}}/k)\quad \\[-8pt]
\times\Gal( \dot{\wt L}_{L_{p}}/k)
\end{array}\to 1\\
@VV{\delta \rho _{L_j}}V
@VV{\rho _{L_j}}V
@VV{\rho _{L}}V\\
{1\to 
\GL_n ( {\delta L}_{R_{p+j}}\times{\delta L}_{L_{p+j}})}
@>>>
{\GL_n ( { L}_{R_{p+j}}\times{ L}_{L_{p+j}})} 
@>{Qh_{G_j\to G}}>>
{\GL_n ( { L}_{R_{p}}\times{ L}_{L_{p}})\to 1} 
\end{CD}\]
of which ``Weil kernel'' is 
$\Gal(\delta \dot{\wt L}_{R_{p+j}}/k)
\times\Gal(\delta \dot{\wt L}_{L_{p+j}}/k)$ and `` $\GL_n(\pt\times \pt)$ kernel'' is
$\GL_n ( {\delta L}_{R_{p+j}}\times{\delta L}_{L_{p+j}})$.
}
\vskip 11pt

\begin{proof}
The proof is similar to that of proposition 4.2.2.

A $n$-dimensional global bilinear quantum deformation of $\rho _L$ is an equivalence class of liftings:
\[\rho _{L_j}=\rho _L+\delta \rho _{L_j}\]
where two liftings 
$\rho _{L_{j_1}}$ and
$\rho _{L_{j_2}}$ are strictly equivalent if they can be transformed one into another by conjugation by bielements of 
$\GL_n ( { L}_{R_{p+j}}\times{ L}_{L_{p+j}})$ in the kernel of
$ Qh_{G_j\to G}$.

The lifting $\rho _{L_{j_1}}$ generates the algebraic bilinear semigroup $\GL_n ( { L}_{R_{p+j_1}}\times{ L}_{L_{p+j_1}})$
having a rank
\[ r_{G_n ( { L}_{R_{j_1}}\times{ L}_{L_{j_1}})}
=((N)^n\cdot f^n_{\o v_{(+j_1)}})\times
((N)^n\cdot f^n_{v_{(+j_1)}})\]
where
\begin{align*}
f_{\o v_{(+j_1)}}\equiv
f_{v_{(+j_1)}}
&= \sum_h m_{p+h+j_1}(p+h+j_1)\\
&= \sum_h \sum_{m_{p+h+j_1}} f_{v_{p+h+j_1,m_{m+h+j_1}}}\;.\end{align*}
Proceeding similarly for the lifting $\rho _{L_{j_2}}$, we find that the difference of ranks between 
$\GL_n ( { L}_{R_{p+j_1}}\times{ L}_{L_{p+j_1}})$ and
$\GL_n ( { L}_{R_{p+j_2}}\times{ L}_{L_{p+j_2}})$ is
\[ \delta r_{G_n(j_2-j_1)}=N^{n^2}\ \L(
f^{n^2}_{v_{(+j_2)}}- f^{n^2}_{v_{(+j_1)}}\R)\;.\]
The deformed algebraic bilinear semigroup
$\GL_n ( { L}_{R_{p+j_1}}\times{ L}_{L_{p+j_1}})$
can be transformed into
$\GL_n ( { L}_{R_{p+j_2}}\times{ L}_{L_{p+j_2}})$ by conjugation of bielements in the first kernel
$\GL_n ( {\delta  L}_{R_{p+j_1}}\times{\delta  L}_{L_{p+j_1}})$
bringing it into the second kernel
$\GL_n ( {\delta  L}_{R_{p+j_2}}\times{\delta  L}_{L_{p+j_2}})$.
\end{proof}
\vskip 11pt

\subsubsection{Corollary}

{\em The transformation of kernels
\[\GL_n ( {\delta  L}_{R_{p+j_1}}\times{\delta  L}_{L_{p+j_1}})\To
\GL_n ( {\delta  L}_{R_{p+j_2}}\times{\delta  L}_{L_{p+j_2}})\]
corresponds to a base change of
$\GL_n ( {  L}_{R_{p+j_1}}\times{  L}_{L_{p+j_1}})$ into
$\GL_n ( {  L}_{R_{p+j_2}}\times{  L}_{L_{p+j_2}})$ of which dimension is given by the difference of ranks
\[ \delta r_{G_n(j_2-j_1)}=N^{n^2}\ \L(
f^{n^2}_{v_{(+j_2)}}- f^{n^2}_{v_{(+j_1)}}\R)\;.\]
}
\vskip 11pt

\subsubsection{Definition (Local bilinear ``quantum'' deformation representative)}

A local bilinear ``quantum'' deformation representative resulting from a local bilinear coefficient semiring ``quantum'' homomorphism:
\[ Qh_{B^-_{p_t}\times B^+_{p_t}\to B^-_{p_r}\times B^+_{p_r}}:\qquad
{B^-_{p_t}\times B^+_{p_t}\To B^-_{p_r}\times B^+_{p_r}}\;,\]
inducing an isomorphism on their inertia subgroups having the same ramification index $e_{\beta ^+_p}$, is an equivalence class representative $\rho _{K{p_t}}$ of lifting
\[\begin{CD}
\Gal(K^{-}_{p_t}/L^-_p)\times\Gal(K^{+}_{p_t}/L^+_p)
@>{Qh_{K_{p_t}\to K_{p_r}}}>>
\Gal(K^{-}_{p_r}/L^-_p)\times\Gal(K^{+}_{p_r}/L^+_p)\\
@VV{\rho _{K_{p_t}}}V @VV{\rho _{K_{p_r}}}V\\
\GL_n(K^{-}_{p_t}\times K^{+}_{p_t})
@>{Gh_{\GL_{(t)}\to \GL_{(r)}}}>>
\GL_n(K^{-}_{p_r}\times K^{+}_{p_r})
\end{CD}\]
of which notations refer to section 4.3.5.
\vskip 11pt

\subsubsection{Proposition (Local bilinear ``quantum'' deformation)}

{\em
Let $K(Qh_{B^-_{p_t}\times B^+_{p_t}\to B^-_{p_r}\times B^+_{p_r}})$
be the kernel of the bihomomorphism
$Qh_{B^-_{p_t}\times B^+_{p_t}\to B^-_{p_r}\times B^+_{p_r}}$ where
${B^-_{p_t}}\times B^+_{p_t}$ belongs to the bicategory
$\Cs({B^-_{p_t}\times B^+_{p_t}})$, $r\le t\le r+s$.

Let 
\begin{multline*}
\Gal(\delta K^{-}_{p_t}/L^-_p)\times\Gal(\delta K^{+}_{p_t}/L^+_p)\\
=
[\Gal(K^{-}_{p_t}/L^-_p)-\Gal(K^{-}_{p_r}/L^-_p)]\times
[\Gal(K^{+}_{p_t}/L^+_p)-\Gal(K^{+}_{p_r}/L^+_p)]
\end{multline*}
be the Weil subgroup corresponding to this kernel.

Then, {\bbf a $n$-dimensional local bilinear ``quantum'' deformation of $\rho _{K_{p_r}}$ is an equivalence class of liftings\/}
$\{\rho _{K_{p_t}}\}$ described by the following exact sequence:
\[
\begin{CD}
{1 \quad \to \begin{array}[t]{r}
\Gal(\delta{ K}^{-}_{p_t}/L^-_p)\quad \\[-8pt]
\times
\Gal(\delta{ K}^{+}_{p_t}/L^+_p)
\end{array}}
@>>>
\begin{array}[t]{r}
\Gal({ K}^{-}_{p_t}/L^-_p)\quad\\[-8pt]
 \times
\Gal({ K}^{+}_{p_t}/L^+_p)\end{array}
@>{Qh_{K_{p_t}\to K_{p_r}}}>>
\begin{array}[t]{r}
\Gal({ K}^{-}_{p_r}/L^-_p)\quad\\[-8pt]
 \times
\Gal({ K}^{+}_{p_r}/L^+_p)\end{array}\to 1\\
@VV{\delta \rho _{K_{p_t}}}V
@VV{\rho _{K_{p_t}}}V
@VV{\rho _{K_{p_r}}}V\\
{1\to \GL_n ( \delta { K}^{-}_{p_t}\times\delta {K}^{+}_{p_t} )}
@>>>
{\GL_n (  { K}^{-}_{p_t}\times {K}^{+}_{p_t} )} 
@>>{Qh_{\GL(t)\to\GL(r)}}>
{\GL_n (  { K}^{-}_{p_r}\times {K}^{+}_{p_r} )\to 1} 
\end{CD}\]
of which ``Weil kernel'' is 
$\Gal(\delta{ K}^{-}_{p_t}/L^-_p)\times\Gal(\delta{ K}^{+}_{p_t}/L^+_p)$ and ``$\GL_n(\pt\times\pt)$ kernel'' is\linebreak 
$\GL_n ( \delta  { K}^{-}_{p_t}\times \delta {K}^{+}_{p_t} )$.
}
\vskip 11pt

\begin{proof}
A $n$-dimensional local bilinear ``quantum'' deformation of $\rho  _{K_{p_r}}$ is an equivalence class of liftings
\[ \rho  _{K_{p_t}}=\rho  _{K_{p_r}}+\delta \rho _{K_{p_t}}\]
where two local liftings
$\rho  _{K_{p_{t_1}}}$ and
$\rho  _{K_{p_{t_2}}}$ are strictly equivalent if they can be transformed one into another by conjugation of bielements of
$\GL_n (  { K}^{-}_{p_{t_1}}\times {K}^{+}_{p_{t_1}} )$ in the kernel 
$\GL_n ( \delta  { K}^{-}_{p_{t_1}}\times\delta {K}^{+}_{p_{t_1}} )$ of
$Qh_{\GL(t_1)\to\GL(r)}$.

In other words, the difference of ranks of the two bilinear algebraic semigroups
$\GL_n (  { K}^{-}_{p_{t_2}}\times {K}^{+}_{p_{t_2}} )$
and
$\GL_n (  { K}^{-}_{p_{t_1}}\times {K}^{+}_{p_{t_1}} )$
is
\[\delta r_{\GL_n(t_2-t_1)}=e^{n^2}_{\beta _p}\ \L(
f^{n^2}_{k_{K^+_{p_{t_2}}}}-
f^{n^2}_{k_{K^+_{p_{t_1}}}}\R)\]
as it results from proposition 4.3.6.

Thus, the deformed algebraic bilinear semigroups
$\GL_n (  { K}^{-}_{p_{t_1}}\times {K}^{+}_{p_{t_1}} )$
and
$\GL_n (  { K}^{-}_{p_{t_2}}\times {K}^{+}_{p_{t_2}} )$
differ from the algebraic bilinear semigroup
$\GL_n (  { K}^{-}_{p_{r}}\times {K}^{+}_{p_{r}} )$ by their respective kernels in such a way that
$\GL_n (  { K}^{-}_{p_{t_1}}\times {K}^{+}_{p_{t_1}} )$
can be transformed into
$\GL_n (  { K}^{-}_{p_{t_2}}\times {K}^{+}_{p_{t_2}} )$
by conjugation of bielements in their first kernel bringing it into the second kernel and inversely.
\end{proof}

\subsubsection{Corollary}

{\em The transformation of kernels
\[\GL_n (  {\delta K}^{-}_{p_{t_1}}\times {\delta K}^{+}_{p_{t_1}} )\To
\GL_n (  {\delta K}^{-}_{p_{t_2}}\times {\delta K}^{+}_{p_{t_2}} )\]
corresponds to a base change of
$\GL_n (  { K}^{-}_{p_{t_1}}\times {K}^{+}_{p_{t_1}} )$
into
$\GL_n (  { K}^{-}_{p_{t_2}}\times {K}^{+}_{p_{t_2}} )$
of which dimension is given by the difference of ranks
$\delta r_{\GL_n(t_2-t_1)}$.
}
\section{Inverse quantum lifts and the Goldbach conjecture}

\subsection{Global bilinear elliptic quantum deformations}

\subsubsection{Definition ($n$-dimensional global elliptic bilinear quantum deformation)}

According to proposition 4.4.2, a $n$-dimensional global bilinear quantum deformation of
\[\rho _L: \qquad \Gal(\dot{\wt L}_{R_p}/k)\times
\Gal(\dot{\wt L}_{L_p}/k)\To \GL_n ( L_{R_p}\times L_{L_p})\]
is an equivalence class of liftings
\[\rho _{L_j}=\rho _L+\delta \rho _{L_j}\]
where $\rho _{L_j}$ is the morphism:
\[\rho _{L_j}: \qquad \Gal(\dot{\wt L}_{R_{p+j}}/k)\times
\Gal(\dot{\wt L}_{L_{p+j}}/k)\To \GL_n ( L_{R_{p+j}}\times L_{L_{p+j}})\;.\]
Taking into account the Langlands global correspondence recalled in proposition 3.1.5, we can define a $n$-dimensional global elliptic bilinear quantum deformation of
\[\rho _L^{\ELLIP}=\ELLIP\FREPsp(\GL_n(L_{R\times L_p})\circ\rho _L\]
where
\[ \ELLIP\FREPsp(\GL_n(L_{R\times L_p}):\qquad
\GL_n(L_{R_p}\times L_{L_p})\To
\ELLIP\RL(n,i\ge p,m_i)\]
is the epimorphism from the bilinear algebraic semigroup
$\GL_n(L_R\times L_p)$ into the $n$-dimensional global elliptic bisemimodule
\begin{multline*}
\ELLIP\RL(n,i\ge p,m_i)\\
=
\L(\sum^t_{i=p}\sum_{m_i}\lambda ^\half(n,i\ge p,m_i)\ e^{-2\pi i(i)x}\R)
\otimes_D\L(\sum^t_{i=p}\sum_{m_i}\lambda ^\half(n,i\ge p,m_i)\ e^{2\pi i(i)x}\R)\;,\\
x\in\rit^n\ , \; i=p+h\ , \; p\le i\le t\le\infty \ , \; 0\le h\le\infty \ ,\end{multline*}
introduced in corollary 3.2.6.
\vskip 11pt

{\bbf This $n$-dimensional global elliptic bilinear quantum deformation of
$\rho _L^{\ELLIP}$ is an equivalence class of liftings\/}:
\[
\rho^{\ELLIP} _{L_j}=\rho^{\ELLIP} _L+\delta \rho^{\ELLIP} _{L_j}\]
in such a way that $\rho^{\ELLIP} _{L_j}$ be given by:
\[
\rho^{\ELLIP} _{L_j}=
\ELLIP\FREPsp(\GL_n(L_{R\times L_{(p+j)}})\circ\rho _{L_j}\]
where
\[ \ELLIP\FREPsp(\GL_n(L_{R\times L_{(p+j)}}):\qquad
\GL_n(L_{R_{p+j}}\times L_{L_{p+j}})\To
\ELLIP\RL(n,i\ge p+j,m_i)\]
is the epimorphism from $\GL_n(L_{R_{p+j}}\times L_{L_{p+j}})$ into the $n$-dimensional global elliptic bisemimodule
\begin{multline*}
\ELLIP\RL(n,i\ge p+j,m_i)=
\L(\sum^{t+j}_{i=p+j}\sum_{m_i}\lambda ^\half(n,i\ge p+j,m_i)\ e^{-2\pi i(i)x}\R)\\
\qquad \qquad \qquad \otimes_D\L(\sum^{t+j}_{i=p+j}\sum_{m_i}\lambda ^\half(n,i\ge p+j,m_i)\ e^{2\pi i(i)x}\R)\;.\end{multline*}
\vskip 11pt

\subsubsection{Proposition}

{\em The $n$-dimensional elliptic global bilinear quantum deformation $\rho^{\ELLIP} _{L_j}$ of $\rho^{\ELLIP} _{L}$ results from the $n$-dimensional global bilinear quantum deformation
$\rho _{L_j}$ of $\rho  _{L}$ by the following commutative diagram:
\[\hspace*{1.5cm} \psmatrix[colsep=3cm]
\Gal(\dot{\wt L}_{R_{p+j}}/k)\times\Gal(\dot{\wt L}_{L_{p+j}}/k) & 
\Gal(\dot{\wt L}_{R_{p}}/k)\times\Gal(\dot{\wt L}_{L_{p}}/k) \\
\GL_n(L_{R_{p+j}}\times L_{L_{p+j}}) &
\GL_n(L_{R_{p}}\times L_{L_{p}})\\
\ELLIP\RL(n,i\ge p+j,m_i) & \ELLIP\RL(n,i\ge p,m_i) &
\psset{nodesep=10pt,arrows=->}
\everypsbox{\scriptstyle}
\ncline{1,1}{1,2}^{Qh_{L_j\to L}}
\ncline{2,1}{2,2}^{Qh_{G_j\to G}}
\ncline{3,1}{3,2}^{Qh_{\rm{EL}_j\to {\rm {EL}}}}
\ncline{1,1}{2,1}>{\rho _{L_j}}
\ncline{1,2}{2,2}>{\rho _{L}}
\ncline{2,1}{3,1}<{\ELLIP}>{\FRepsp(\GL_n(L_{R\times L_{(p+j)}}))}
\ncline{2,2}{3,2}<{\ELLIP}>{\FRepsp(\GL_n(L_{R\times L_{p}}))}
\ncangle[angleA=-180,angleB=180,armB=5mm]{1,1}{3,1}
\Aput{\begin{array}{c} \rho ^{\ELLIP}_{L_j}\\ \mbox{}\\ \mbox{}\end{array}}
\ncangle[angleA=360,angleB=-360,armB=1.5cm]{1,2}{3,2}
\Bput{\begin{array}{c} \rho ^{\ELLIP}_{L}\\ \mbox{}\\ \mbox{}\end{array}}
\endpsmatrix
\]
in such a way that the injective morphism
$\DD^{\{p\}\to\{p+j\}}\RL(n)=Qh^{-1}_{\rm{EL}_j\to\rm{EL}}$:
\[ \DD^{\{p\}\to\{p+j\}}\RL(n): \qquad
\ELLIP\RL(n,i\ge p,m_i)\To
\ELLIP\RL(n,i\ge p+j,m_i)\]
be a {\bbf quantum deformation of the $n$-dimensional global elliptic bisemimodule\linebreak
$\ELLIP\RL(n,i\ge p,m_i)$\/}.
}
\vskip 11pt

\begin{proof}
This results immediately from definition 5.1.1.
\end{proof}
\vskip 11pt

\subsubsection{Proposition}

{\em Let
\[ \DD^{[p]\to[p+j]}\RL(n): \qquad
\ellip\RL(n,[p],m_p)\To
\ellip\RL(n,[p+j],m_{p+j})\]
denote {\bbf a quantum equivalence class representative of liftings called in a more condensed form a quantum deformation of the $(p,m_p)$-th conjugacy class representative
$\ellip\RL(n,i\ge p,m_i)$ of $
\ELLIP\RL(n,i\ge p,m_i)$\/}, i.e. the $(p,m_p)$-th term of
$\ellip\RL(n,i\ge p,m_i)$.

This elliptic quantum deformation is then associated with the exact sequence:
\[
1\To
\ellip\RL(n,[j])\To
\ellip\RL(n,[p+j],m_{p+j})\To
\ellip\RL(n,[p],m_p)\To
1\;.\]
}
\vskip 11pt

\begin{proof}
Indeed, this exact sequence results from the commutative diagram:
\[
\begin{CD}
{1  \To \begin{array}[t]{r}
\Gal( \dot{\wt L}_{\o v_j}/k)\quad \\[-8pt]
\times
\Gal( \dot{\wt L}_{ v_j}/k)
\end{array}}
@>>>
\begin{array}[t]{r}
\Gal( \dot{\wt L}_{\o v_{p+j}}/k)\quad\\[-8pt]
 \times
\Gal( \dot{\wt L}_{v_{p+j}}/k)
\end{array}
@>>>
\begin{array}[t]{r}
\Gal( \dot{\wt L}_{\o v_{p}}/k)\quad\\[-8pt]
 \times
\Gal( \dot{\wt L}_{v_{p}}/k)\end{array}\To 1\\
@VVV
@VVV
@VVV\\
{1\To \gl_n({ L}_{\o v_j}\times {L}_{v_j})}
@>>>
{\gl_n({ L}_{\o v_{p+j}}\times {L}_{v_{p+j}})} 
@>>>
{\gl_n({ L}_{\o v_p}\times {L}_{v_p})\To 1} \\
@VVV
@VVV
@VVV\\
{1\To \ellip\RL(n,[j])}
@>>>
\ellip\RL(n,[p+j],m_{p+j})
@>>>
{\ellip\RL(n,[p],m_p)\To 1} 
\end{CD}\]
where:
\Bi
\item $\gl_n({ L}_{\o v_{p+j}}\times {L}_{v_{p+j}})\in
\GL_n(L_{R_{p+j}}\times L_{L_{p+j}})$ is the $(p+j)$-th conjugacy class of the bilinear algebraic semigroup $\GL_n(L_{R_{p+j}}\times L_{L_{p+j}})$;

\item $\gl_n({ L}_{\o v_j}\times {L}_{v_j})$ is the kernel at ``$j^{n^2}$'' biquanta of the exact sequence
\[
{1\To \gl_n({ L}_{\o v_j}\times {L}_{v_j})}
\To{\gl_n({ L}_{\o v_{p+j}}\times {L}_{v_{p+j}})} 
\To
{\gl_n({ L}_{\o v_p}\times {L}_{v_p})\To 1}
\] if it is referred to corollary 4.3.3.
\Ei
\end{proof}
\vskip 11pt


\subsection{Inverse elliptic quantum deformations and the Goldbach conjecture}

\subsubsection{Definition (Inverse elliptic quantum deformations)} 
Let $\DD\RL^{[p+j]\to[p+j+k]}(n)$ denote the equivalence class representative of lifting of the global elliptic subbisemimodule $\ellip\RL(n,[p+j],m_{p+j})$,  of class $[p+j]$, towards the global elliptic subbisemimodule $\ellip\RL(n,[p+j+k],m_{p+j+k})$ of class $[p+j+k]$ as considered in proposition 5.1.3.

Then, the {\bbf inverse deformation $\DD\RL^{[p+j+k]\to[p+j]}(n)$ can be introduced by the  surjective mapping\/}:
\begin{align*}
\DD\RL^{[p+j+k]\to[p+j]}(n) :\quad &\ellip\RL(n,[p+j+k],m_{p+j+k})\\
&\qquad \To  \ellip\RL(n,[p+j],m_{p+j})\end{align*}
which is associated with the  exact sequence:
\[
1\to \ellip\RL(n,[k]) \to \ellip\RL(n,[p+j+k],m_{p+j+k})
 \to \ellip\RL(n,[p+j],m_{p+j})\to 1\]
{\bbf in such a way that:
\[
\ellip\RL(n,[p+j],m_{p+j}) \simeq
\ellip\RL(n,[p+j+k],m_{p+j+k})- \ellip\RL(n,[k])\]
corresponds to the endomorphism\/}:
\begin{align*}
\End\RL^{[p+j+k]\to[p+j]}(n) : & \quad \ellip\RL(n,[p+j+k],m_{p+j+k})\\
&\qquad \To \ellip\RL(n,[p+j],m_{p+j})+ \ellip\RL(n,[k])\;.\end{align*} \vskip 11pt

\subsubsection{Proposition}  

{\em The set of inverse quantum deformations 
$\{\DD^{[p+j+k]\to[p+j]}\RL(n)\}^{k^{\rm{up}}}_{k=1}$ are artinian deformations.} \vskip 11pt

\begin{proof} 
Indeed, the set of  class representatives of decreasing global elliptic subbisemimodules:
\begin{multline*}
\ellip\RL(n,[p+j+k^{\rm up}],m_{p+j+k^{\rm up}}) \subset \cdots 
\subset \ellip\RL(n,[p+j+k],m_{p+j+k})\\
\qquad \subset \cdots \subset \ellip\RL(n,[p+j],m_{p+j})\;,\qquad 1\le k\le k^{\rm{up}}\le\infty \;,\end{multline*}
generated under the set of inverse quantum deformations $\{\DD\RL^{[p+j+k]\to[p+j]}(n)\}^{k^{\rm up}}_{k=1}$, forms an artinian sequence.\epr 
\vskip 11pt

\subsubsection{Inverse quantum deformations of one-dimensional tori}  
Consider the inverse quantum deformation: 
\begin{align*}
 \DD\RL^{[p+j+k]\to[p+j]}(1) : \quad   &\ellip\RL(1,[p+j+k],m_{p+j+k}) \\
&\qquad \To
\ellip\RL(1,[p+j],m_{p+j})\end{align*}
of a  global elliptic subbisemimodule
$\ellip\RL(1,[p+j+k],m_{p+j+k})$
from a class $[p+j+k]$ to a class $[p+j[$,\newline
where
\Bi
\item $\ellip\RL(1,[p+j],m_{p+j})\simeq T^1_R[p+j,m_{p+j}] \times T^1_L[p+j,m_{p+j}]$ 
is in bijection with the product, right by left, of the analytic developments of semitori of 
class $[p+j]$, i.e. having a rank $r^{(1)}_{v_{p+j}}=(p+j)\cdot N$ such that $(p+j)$ is a global  residue degree;
\item $\ellip\RL(1,[p+j+k],m_{p+j+k})\simeq T^1_R[p+j+k,m_{p+j+k}]\times 
T^1_L[p+j+k,m_{p+j+k}]$ of which semitori have a rank $r^{(1)}_{v_{p+j+k}}=(p+j+k)\cdot N$.
\Ei

\vskip 11pt

The semitorus (or semicircle) $T^1_L[p+j,m_{p+j}]$ is defined in the upper half plane.  Consequently, if we have to consider an analytic continuation of it to the whole plane, we have to undouble it according to:
\[T^1_L[p+j,m_{p+j}] \To T^1_{2_L}[2(p+j),m_{2(p+j)}]\equiv T^1_{2_L}[p'+j',m_{p'+j'}]\]
such that:
\Bi
\item   $T^1_{2_L}[\cdot]$ has a rank $r^{(1)}_{v_{2(p+j)}}=2\cdot (p+j)\cdot N$;
\item  $T^1_{2_L}[p'+j',m_{p'+j'}]$ has been defined with respect to a new prime number $p'\neq p$.
\Ei
\vskip 11pt

{\bbf So, a general class representative of inverse lifting corresponding to an inverse quantum deformation of a $1D$-torus of class $[p'+j'+k']$ will be defined by the projective mapping\/}:
\[ \DD_L^{[p'+j'+k']\to[p'+j']}(1) : \quad   T^1_{2_L}[p'+j'+k',m_{p'+j'+k'}] \To
T^1_{2_L}[p'+j',m_{p'+j'}]\]
{\bbf associated with the endomorphism\/}:
\begin{align*}
 \End_L^{[p'+j'+k']\to[p'+j']}(1) : 
\quad    &T^1_{2_L}[p'+j'+k',m_{p'+j'+k'}] \\
&\qquad \To
T^1_{2_L}[p'+j',m_{p'+j'}]
+T^1_{2_L}[k',m_{k'}]\end{align*}
splitting the circle of class $[p'+j'+k']$ into two complementary portions of circles of classes $[p'+j']$ and $[k']$.
\vskip 11pt

If the global residue degrees $f_{v_{p'+j'}}=p'+j'$ and $f_{v_{k'}}=k'$ are not even integers, then $T^1_{2_L}[p'+j',m_{p'+j'}]$
 and $T^1_{2_L}[k',m_{k'}]$ are not ``closed'' circles.
\vskip 11pt

But, in any case, we have the following equality between the global residue degrees 
associated with the endormophism $\End_L^{[p'+j'+k']\to[p'+j']}(1) $ of the circle 
$T^1_{2_L}[p'+j'+k',m_{p'+j'+k'}]$:
\[ f_{v_{p'+j'+k'}}=f_{v_{p'+j'}}+f_{v_{k'}} \quad \Rightarrow \quad p'+j'+k' = (p'+j')+(k')\]
where $f_{v_{p'+j'+k'}}=p'+j'+k'$ is an even integer since it is assumed to be the global residue degree of a one-dimensional torus undoubled from the corresponding semicircle.
\vskip 11pt

\subsubsection{Proposition (Goldbach conjecture)}  

{\em Let $c^1_{2_L}[p'+j'+k',m_{p'+j'+k'}]$ denote a closed curve isomorphic 
to the one-dimensional torus $T^1_{2_L}[p'+j'+k',m_{p'+j'+k'}]$ of class $(p'+j'+k')$ generated 
from the corresponding semitorus localized in the upper half plane.

Let
\[ \DD_L^{[p'+j'+k']\to[p'+j']}(1) : \quad   c^1_{2_L}[p'+j'+k',m_{p'+j'+k'}] \To
c^1_{2_L}[p'+j',m_{p'+j'}]\]
be the projective inverse quantum deformation of the closed curve 
$c^1_{2_L}[p'+j'+k',m_{p'+j'+k'}] $ associated with the endomorphism
\begin{align*}
 \End_L^{[p'+j'+k']\to[p'+j']}(1) : \quad    &c^1_{2_L}[p'+j'+k',m_{p'+j'+k'}] \\
&\qquad \To
c^1_{2_L}[p'+j',m_{p'+j'}]
+c^1_{2_L}[k',m_{k'}]\end{align*}
splitting the closed curve $C^1_{2_L}[p'+j'+k',m_{p'+j'+k'}]$ into two complementary portions,\newline
in such a way that:
\Bi
\item the curves $c^1_{2_L}[p'+j',m_{p'+j'}]
$ and $c^1_{2_L}[k',m_{k'}]$ are not necessarily closed;
\item the global residue degrees of these curves verify:
\[ f_{v_{p'+j'+k'}}=f_{v_{p'+j'}}+f_{v_{k'}} \quad \Rightarrow \quad p'+j'+k' = (p'+j')+(k')\]
where $G_{\rm {even}}=2G=p'+j'+k'$ is an even integer.
\Ei 
Then, to the even integer $G_{\rm{even}}$, there corresponds at least one basic class representative of inverse lifting corresponding to the inverse deformation $\DD_L^{[p'+j'+k']\to[p'+j']}(1)$ of the closed curve $c^1_{2_L}[p'+j'+k',m_{p'+j'+k'}] $ such that the global residue degrees associated to
$\DD_L^{[p'+j'+k']\to[p'+j']}(1)$  verify:
\begin{align*}
G_{\rm{even}}=f_{v_{p'+j'}}+f_{v_{k'}}\;, && G_{\rm{even}}\le\infty \;, ,\end{align*}
where:
\[
f_{v_{p'+j'}}  = p'+j' \qquad
\text{and} \qquad f_{v_{k'}}  = q'+n' = k'\]
in such a way that the even integer $F_{\rm{even}}
=G_{\rm{even}}-j'-n'=p'+q'$, is the sum of two prime numbers $p'$ and $q'$, $4\le F_{\rm{even}}\le\infty $.
} \vskip 11pt

\bpr \Be
\item To the set of one-dimensional tori $\{T^1_{2_L}[p'+j'+k',m_{p'+j'+k'}]\}_{j',k'} $ is associated a classical Riemann $\zeta$-function whose trivial zeros are negative even integers which are proved to be in one-to-one correspondence with the global residue degrees
\[ f_{v_{p'+j'+k'}}=  p'+j'+k' = 2G\]
as developed in \cite{Pie2}.

It was also seen in \cite{Pie2} that the trivial zeros and the nontrivial zeros of the classical $\zeta$-function all on the line $\sigma=\half$ are in one-to-one correspondence under the action of the Lie algebra of the decomposition group.

So, if we admit that the Riemann hypothesis is verified, then, the prime number theorem follows:
\[\Pi(x)\sim \F x{\ln x}\qquad \text{when}\quad  x\to \infty\]
\cite{B-K}, \cite{Ing}, \cite{K-S},\newline
where $\Pi(x)$ is the number of primes not exceeding $x$.

\item To every even integer $2G$, there corresponds a pair of nontrivial zeros of the Riemann zeta function:
\[\zeta(s) =\sum^\infty_{n=1}n^{-s}=\prod_p\ (1-p^{-s})^{-1}\]
and a set of prime numbers inferior to $2G$.

So, according to the prime number theorem fixing the density of the prime numbers, it is likely that there exists at least a pair $\{p',q'\}$ of prime numbers whose sum is an even integer.

But, the question is now the following: does there always exist a pair $\{p',q'\}$ of prime numbers whose sum is an even integer?  The response is affirmative.

Indeed, the integers $G_{\rm{even}}$, $f_{v_{p'+j'}}$ and $f_{v_{k'}}$, are interpreted as global residue degrees associated with the inverse quantum deformation $\DD_L^{[p'+j'+k']\to[p'+j']}(1)$  of the closed curve\linebreak $c^1_{2_L}[p'+j'+k',m_{p'+j'+k'}]   $.  Now, the  global residue degree $G_{\rm{even}}=p'+j'+k'$ was defined in 2.1 and in \cite{Pie1} by:
\[[L^{nr}_{v_{p'+j'+k'}}:k] = f_{v_{p'+j'+k'}} = hp'+j''+k''=p'+j'+k'\;,\]
where $j''$ and $k''$ are integers, referring to congruence classes modulo $p'$, and where $p'$ is a prime number.

Similarly, we have that
\begin{align*}
[L^{nr}_{v_{p'+j'}}:k] &= f_{v_{p'+j'}} = \ell p'+j'' = p'+j'\\
\noalign{\qquad \qquad and}
[L^{nr}_{v_{k'}}:k] &= f_{v_{k'}} = m q'+n = q'+n'=k'\end{align*}
where $q'$ is a prime number.

Indeed, the global  residue degrees can be defined with respect to prime numbers in order that the corresponding local fields  be able to be handled by classical $p'$ (resp. $q'$)-adic methods as it was developed in section 3.3.

We thus have that
\[ F_{\rm{even}} = G_{\rm{even}}-j'-n'=p'+q'\]
where $j'$ and $n'$ are both either even or odd integers if $F_{\rm{even}}$ is an even integer.

Now, it is always possible to find pairs of even or odd integers $j'$ and $n'$ in order that $F_{\rm{even}}$ be an even integer as it results from section 3.2.1.

Finally, if $F_{\rm{even}}$ could not be defined as the sum $(p'+q')$ of two primes $p'$ and $q'$, that should mean that $G_{\rm{even}}$ could not be equal to the sum $(f_{v_{p'+j'}}+f_{v_{k'}})$ of two global residue degrees and, thus, that it could not be possible to define all class representatives of inverse liftings of every closed curve, which is absurd.
\qedhere 
\Ee
\end{proof}\vskip 11pt

As a result, we have the \vskip 11pt

\subsubsection{Goldbach's conjecture}  

{\em Every even integer superior or equal to 4 is the sum of two prime numbers.} 

\vfill
Department of Mathematics\\
Universit\'e de Louvain\\
Chemin du cyclotron, 2\\
B--1348 Louvain-la-Neuve, Belgium\\
pierre.math.be@gmail.com

\end{document}